\documentclass[11pt, a4paper, final]{article}
\usepackage[english]{babel}
\usepackage[utf8]{inputenc}
\usepackage{amsmath}
\usepackage{amssymb}
\usepackage{amsthm}
\usepackage{txfonts}
\usepackage{color}
\usepackage{enumerate}
\usepackage[inline]{enumitem}
\usepackage{xspace}										
\usepackage{bbm}

\usepackage[a4paper, left=3cm, right=3cm, top=2.5cm, bottom=2.5cm]{geometry}

\usepackage{standalone}									
\usepackage{xcolor}										
\usepackage{tikz}										
\usetikzlibrary{%
	arrows
}

\usepackage{ragged2e}
\usepackage[
sorting = nyt,
sortcites,
giveninits = true,
isbn = false,
doi = false, 
url = false, 
clearlang = true,
backend = biber
]{biblatex}
\AtEveryBibitem{\clearfield{pagetotal}}					
\AtEveryBibitem{\clearfield{note}}						
\usepackage{csquotes}	
\DeclareRedundantLanguages{English,english,german,french, eng}{english,german,ngerman,french}

\renewbibmacro{in:}{%
	\ifentrytype{article}{}{\printtext{\bibstring{in}\intitlepunct}}}		

\renewcommand{\textcite}{\cite}

\usepackage{aliascnt}		
\usepackage[pdfpagelabels=true,plainpages=false, hidelinks]{hyperref}
\usepackage[nameinlink]{cleveref}

\theoremstyle{plain}
\newtheorem{theorem}{Theorem}[section]
\newtheorem*{theorem*}{Theorem}

\newaliascnt{lemma}{theorem}
\newtheorem{lemma}[lemma]{Lemma}
\aliascntresetthe{lemma}
\crefname{lemma}{lemma}{lemmas}

\newaliascnt{cor}{theorem}
\newtheorem{cor}[cor]{Corollary}
\aliascntresetthe{cor}

\newaliascnt{prop}{theorem}
\newtheorem{prop}[prop]{Proposition}
\aliascntresetthe{prop}
\crefname{prop}{proposition}{propositions}

\theoremstyle{definition}
\newtheorem*{defi}{Definition}

\theoremstyle{remark}
\newaliascnt{remk}{theorem}
\newtheorem{remk}[remk]{Remark}
\aliascntresetthe{remk}
\crefname{remk}{remark}{remarks}

\newtheorem{ex}{Example}

\newcommand{\RR}{\mathbb{R}}
\newcommand{\OO}{\mathcal{O}}
\newcommand{\FF}{\mathcal F}
\newcommand{\LL}{\mathcal L}
\newcommand{\varell}{\mathfrak{l}}
\newcommand{\End}[1]{\operatorname{End}_\Sigma (#1)}
\newcommand{\gl}[1]{\mathfrak{gl}(#1)}
\newcommand{\gker}[1]{\operatorname{ker}_0 (#1)}
\newcommand{\XX}{\mathcal{X}}
\newcommand{\Xc}{\XX^c}
\newcommand{\Xh}{\XX^h}

\newcommand{\nole}{\trianglelefteq}
\newcommand{\noge}{\trianglerighteq}
\newcommand{\nol}{\vartriangleleft}
\newcommand{\nog}{\vartriangleright}

\newcommand{\ox}{\otimes}

\newcommand{\idv}{\mathbbm{1}_V}
\newcommand{\idvd}{\mathbbm{1}_{\Vd}}

\newcommand{\redim}[1]{\operatorname{im}_0 (#1)}
\newcommand{\musup}{\mu}
\newcommand{\musub}{\mu}
\newcommand{\Dsup}{D}
\newcommand{\Dsub}{D}
\newcommand{\Bsup}{B}

\newcommand{\od}{\textbf{\textsf{1D}}\xspace}
\newcommand{\dd}{\textbf{\textsf{DD}}\xspace}

\newcommand{\No}{\mathcal{N}^{\mathbf{\mathsf{1}}}}
\newcommand{\Nd}{\mathcal{N}^{\mathbf{\mathsf{D}}}}
\newcommand{\Ao}{A^{\mathbf{\mathsf{1}}}}
\newcommand{\Ad}{A^{\mathbf{\mathsf{D}}}}
\newcommand{\Bo}{B^{\mathbf{\mathsf{1}}}}
\newcommand{\Bd}{B^{\mathbf{\mathsf{D}}}}
\newcommand{\Xo}{\XX_{\mathbf{\mathsf{1}}}}
\newcommand{\Xd}{\XX_{\mathbf{\mathsf{D}}}}
\newcommand{\Po}{P_{\mathbf{\mathsf{1}}}}

\newcommand{\Vo}{\RR}
\newcommand{\Vd}{W}
\newcommand{\Gamo}{\Gamma^{\mathbf{\mathsf{1}}}}
\newcommand{\Gamd}{\Gamma^{\mathbf{\mathsf{D}}}}
\newcommand{\Delo}{\Delta^{\mathbf{\mathsf{1}}}}
\newcommand{\Deld}{\Delta^{\mathbf{\mathsf{D}}}}

\newcommand{\uu}[1]{\underline{#1}}
\newcommand{\vu}{\underline{v}}
\newcommand{\Gu}{\underline{\Gamma}}
\newcommand{\Xu}{\underline{\mathcal{X}}}

\newcommand{\pu}{\uu{\psi}}

\DeclareMathOperator{\Id}{Id}

\newcommand{\field}[1]{\mathbb{#1}}
\newcommand{\R}{\field{R}}

\numberwithin{equation}{section}

\definecolor{red}{RGB}{193,56,160}									
\definecolor{blue}{RGB}{0,37,170}

\begin{document}
	
\title{Homogeneous coupled cell systems with high-dimensional internal dynamics}
\author{Sören von der Gracht\thanks{Paderborn University, Paderborn, Germany, \href{mailto:soeren.von.der.gracht@uni-paderborn.de}{soeren.von.der.gracht@uni-paderborn.de}}, Eddie Nijholt\thanks{\mbox{University of São Paulo, São Carlos, Brazil, \href{mailto:eddie.nijholt@gmail.com}{eddie.nijholt@gmail.com}}}, Bob Rink\thanks{\mbox{Vrije Universiteit Amsterdam, Amsterdam, The Netherlands, \href{mailto:b.w.rink@vu.nl}{b.w.rink@vu.nl}}}}
\date{}
\maketitle

\begin{abstract}
	\noindent
	We investigate homogeneous coupled cell systems with high-dimensional internal dynamics. In many studies on network dynamics, the analysis is restricted to networks with one-dimensional internal dynamics. Here, we show how symmetry explains the relation between dynamical behavior of systems with one-dimensional internal dynamics and with higher dimensional internal dynamics, when the underlying network topology is the same. Fundamental networks of homogeneous coupled cell systems (compare to \fullcite{Rink.2014}) 
	can be expressed in terms of monoid representations, which uniquely decompose into indecomposable subrepresentations. 
	In the high-dimensional internal dynamics case, these subrepresentations are isomorphic to multiple copies of those one computes in the one-dimensional internal dynamics case. 
	This has interesting implications for possible center subspaces in bifurcation analysis. 
	We describe the effect on steady state and Hopf bifurcations in $l$-parameter families of network vector fields.  
	The main results in that regard are that (1) generic one-parameter steady state bifurcations are qualitatively independent of the dimension of the internal dynamics and that, (2) in order to observe all generic $l$-parameter bifurcations that may occur for internal dynamics of any dimension, the internal dynamics has to be at least $l$-dimensional for steady state bifurcations and $2l$-dimensional for Hopf bifurcations. 
	Furthermore, we illustrate how additional structure in the network can be exploited to obtain even greater understanding of bifurcation scenarios in the high-dimensional case beyond qualitative statements about the collective dynamics. 
	One-parameter steady state bifurcations in feedforward networks exhibit an unusual amplification in the asymptotic growth rates of individual cells, when these are one-dimensional (\fullcite{vonderGracht.2022}). As another main result, we prove that (3) the same cells exhibit this amplifying effect with the same growth rates when the internal dynamics is high-dimensional.
\end{abstract}

\vspace{.5cm}

\noindent {\em Keywords:} coupled cell systems, network dynamics, dimension reduction, bifurcation theory, symmetry, monoid representation theory

\vspace{.3cm}

\noindent {\em MSC:} 37G40, 37G10, 20M30, 37C81, 37C79
\vspace{1cm}

\section*{Introduction}
Dynamical systems as they arise in fields such as neuroscience (the workings of the brain), systems biology (metabolic systems), and electrical engineering (power grids), exhibit the structure of a network. 
That is, they consist of nodes (neurons, proteins, power stations) with connections between them and the behavior of one cell influences that of another. 
It usually does not suffice to understand the nature of the individual nodes to deduce the behavior of the network, and the specific interaction structure of a network can produce remarkable dynamics. One of the most staggering examples is synchronization (e.g. the simultaneous firing of neurons). The analysis of network dynamical systems is challenging as standard techniques are not tailored to the underlying structure. In recent years, numerous approaches and formalisms have been put forward to robustly encode network structure in dynamical systems. 
The \emph{groupoid formalism} (see \textcite{Golubitsky.2006,Golubitsky.2004} and an equivalent definition in \textcite{Field.2004}) among other uses allows for the classification of synchrony patterns in terms of balanced colorings of the network nodes. 
It has been generalized to \emph{hypernetworks} to specifically model groupwise and indirect interactions which have recently been identified in numerous applications as a main driver of collective dynamics \cite{Battiston.2020,Nijholt.2022b,Aguiar.2023,Bick.2023c,Boccaletti.2023,vonderGracht.2023d,Bick.2024b, vonderGracht.2024}. \emph{Asynchronous networks} \textcite{Bick.2016,Bick.2017,Bick.2017b} and \emph{open systems} \textcite{Lerman.2018,Schultz.2020} are aimed at flexibly and realistically modeling real-world applications.

Slightly more specialized -- that is, for the smaller class of \emph{homogeneous networks} -- is the approach via so-called \emph{homogeneous coupled cell systems} \textcite{Rink.2013,Rink.2014,Rink.2015,Nijholt.2016,Nijholt.2017,Nijholt.2017c} (therein, also generalizations to non-homogeneous networks are made). After the category-theoretic tool of \emph{graph fibrations} was introduced to network dynamics (see \textcite{DeVille.2015,DeVille.2015b}), the investigation of their implications on dynamical systems made it possible to relate the behavior of homogeneous networks to algebraic properties of a representation of a semigroup. In particular, the admissible vector fields of a lift of the network -- the \emph{fundamental network} -- are precisely the vector fields that are equivariant with respect to this representation. This is referred to as \emph{hidden symmetry}. The observation sheds light on the -- often anticipated -- connection between networks and symmetry in the context of homogeneous coupled cell systems. More precisely, equivariant dynamics has proven to be a powerful tool to study the dynamics of complex systems by structuring and reducing systems in terms of representation theory of classical symmetry groups. Even if many networks do not have classical symmetry groups (i. e., graph automorphisms) their dynamics exhibit features that resemble those of symmetric systems. It turns out that more general structures such as quivers, semigroups, and monoids adequately reflect the structural features of non-symmetric networks and their representations can be exploited in a similar manner \textcite{Rink.2014,Nijholt.2020} -- a strategy that we outline and pursue below.

Oftentimes, one is interested in \emph{generic bifurcations} in network dynamical systems. These can be interpreted as a prediction of bifurcation behavior that is dictated only by the network structure and independent of the specific system. Investigations range from small examples to qualitative statements for entire classes of networks. The book \cite{Golubitsky.2023} and references therein give an excellent overview of the current state of the theory and some newer results can for example be found in \cite{vonderGracht.2024,Aguiar.2025}. A major issue in determining generic bifurcations, say in a given network, is the computational complexity stemming from high-dimensional phase spaces. The \emph{total phase space} -- the phase space of the entire network dynamical system -- has at least one spatial direction for each cell of the network. In order to reduce this difficulty to its minimum, one often restricts to the case where the \emph{internal phase space} -- the phase space of a single cell -- is one-dimensional. However, this often requires additional work either motivating this restriction from a modeling point of view or demonstrating its meaning in more general systems, for instance as reduced dynamics.

In the following example, we illustrate in a very simple network an effect high-dimensional internal dynamics may have on the dynamical analysis compared to one-dimensional internal dynamics.
\begin{figure}[h]
	\begin{center}
		\resizebox{.4\linewidth}{!}{	
			\centering
\begin{tikzpicture}[->,
	>=stealth',
	shorten >=1pt,
	auto,
	node distance=1cm,
	main node/.style={line width=1.5pt, circle, scale = 3, draw, font=\sffamily\tiny, inner sep=1pt}]
	\node[main node] (1) {$1$};
	\node[main node, left of=1] (2){$2$};
	\node[main node, left of=2] (3) {$3$};
	\path[every node/.style={font=\sffamily\small}, line width =1.5pt]
	(2) edge [color = {red}] node {} (1)
	(3) edge [color = {red}] node {} (2)
	(3) edge [color = {red}, in = 190, out = 170, looseness = 8] node {} (3)
;
\end{tikzpicture}%
		}%
	\end{center}%
	\caption{A $3$-cell homogeneous feedforward chain.}
	\label{fig:3cellff}
\end{figure}%
\begin{ex}
	\label{ex:int}
	Consider the $3$-cell homogeneous feedforward chain in \Cref{fig:3cellff}. Its dynamics is governed by the system of ordinary differential equations
	\begin{align*}
		\dot{v}_1	&= f(v_1, \textcolor{red}{v_2}) \\
		\dot{v}_2	&= f(v_2, \textcolor{red}{v_3}) \\
		\dot{v}_3	&= f(v_3, \textcolor{red}{v_3}),
	\end{align*}
	where $v_i \in V$ is the state variable of cell $i$ in the internal phase space $V$. In order to investigate dynamical phenomena -- in particular stability -- one analyzes spectral properties of linearizations of the right hand side of such systems. Especially for the investigation of generic steady state bifurcations one is interested in spectral properties of a generic linear right hand side -- that is, a linear admissible map --, which is of the form
	\begin{equation*}
		L = \begin{pmatrix} 
			A & \textcolor{red}{B} & 0 \\
			0 & A & \textcolor{red}{B} \\
			0 & 0 & A+\textcolor{red}{B}
		\end{pmatrix} 
	\end{equation*}
	where $A, \textcolor{red}{B} \in \gl{V}$ are generic linear maps on $V$. The spectrum of $L$ is made up of the eigenvalues of $A$ and those of $A+\textcolor{red}{B}$, where the eigenvalues of $A$ occur with algebraic multiplicity $2$, even though they are generically simple as eigenvalues of $A$. This spectral degeneracy -- a double eigenvalue is unheard of in a generic linear map without any additional structure -- is independent of the dimension of the underlying space $V$.
	
	However, the investigation of generic steady state bifurcations also relies on information about the generalized \emph{eigenspaces} of the linearization at a bifurcation point. Let us, for simplicity, assume that $A$ has an eigenvalue $0$. In the case $V=\RR$, this is equivalent to the assumption $A=0$. In that case, we also have $\textcolor{red}{B}\ne 0$ generically. 
    The generalized eigenspace of the eigenvalue $0$ is spanned by an eigenvector and a generalized eigenvector as 
	\[ E_0 = \left\langle \begin{pmatrix} 1 \\ 0 \\ 0  \end{pmatrix}, \begin{pmatrix} 0 \\ \frac{1}{\textcolor{red}{B}} \\ 0 \end{pmatrix} \right\rangle. \]
	
    If, on the other hand, $V = \RR^d$ for some $d > 1$, the eigenvalue $0$ of $A$ is generically simple. Hence, there is an eigenvector $v \in V$ such that $Av=0$ and no other (generalized) eigenvector. 
    As we are assuming $0$ is a simple eigenvalue of $A$, we have a direct sum decomposition of $V$ into the image of $A$ and the span of $v$. Thus there exists a scalar $c \in \R$ and an element $w \in V$ such that \[-\textcolor{red}{B}v = Aw-cv \, \text{ or }  \, Aw = \left( c\idv - \textcolor{red}{B} \right) v. \]
	Then the generalized eigenspace is spanned by the eigenvector with corresponding generalized eigenvector
	\[ E_0 = \left\langle \begin{pmatrix} v \\ 0 \\ 0  \end{pmatrix}, \begin{pmatrix} w \\ v \\ 0 \end{pmatrix} \right\rangle. \]
	This structure significantly differs from the one in the case $V=\RR$. The generalized eigenvector depends not only on $\textcolor{red}{B}$ but also on $A$. As $A$ and $\textcolor{red}{B}$ do not necessarily commute,
	\[ \begin{pmatrix} 0 \\ \textcolor{red}{B}^{-1}v \\ 0 \end{pmatrix} \]
	is in general not a generalized eigenvector. Summarizing, already this simple $3$-cell feedforward chain produces significantly different spectral properties when the internal dynamics is high-dimensional.
	\hspace*{\fill}$\triangle$	
\end{ex}

In this paper, we show that issues, such as the one illustrated in \Cref{ex:int}, only have a `controllable' qualitative impact on generic steady state bifurcations in homogeneous coupled cell systems. After briefly summarizing the formalism in \Cref{sec:pr}, we use techniques from representation theory to show that critical eigenspaces in bifurcation analysis of networks with high-dimensional internal dynamics are in some sense the same as the ones we encounter in the one-dimensional case. In particular, this allows for dimension-reduction. The main result of the first part of this paper is
\begin{theorem*}
    The decomposition of the monoid representation that is equivalent to the fundamental network with $d$-dimensional internal dynamics is isomorphic to $d$ copies of the decomposition of the total phase space of the same fundamental network with one\nobreakdash-dimensional internal dynamics. As a result, generic $1$-parameter steady state bifurcations are (qualitatively) the same as in the one-dimensional internal dynamics case. Moreover, any generic $l$-parameter bifurcation in a fundamental network with one\nobreakdash-dimensional internal dynamics also occurs generically in the same network with $d$-dimensional internal dynamics. Even more so, every generic $l$-parameter bifurcation that a given fundamental network supports can be observed when the internal dynamics is of a minimal dimension $\overline{d}$ -- , $\overline{d}=l$ for steady state bifurcations and $\overline{d}=2l$ for Hopf bifurcations.
\end{theorem*}
\noindent
Throughout this paper, `the same', `qualitatively the same', `qualitatively equivalent' in terms of bifurcations always means that two equivariant bifurcation problems occur in phase spaces that are equivariantly isomorphic with respect to a representation of a monoid. Hence, the class of \emph{generic} bifurcation problems in both cases are literally the same in a suitably chosen basis. In particular, this implies that the number of solution branches, their branching directions, the leading growth rates, and their stability properties are locally the same.
The main result is proved in multiple theorems. In \Cref{sec:hd}, we provide the algebraic details by comparing the decompositions of the monoid representations that are equivalent to the fundamental network structure with internal dynamics of varying dimensions (\Cref{thm:dechd}). Implications for possible center subspaces in steady state or Hopf bifurcation analysis are discussed in \Cref{sec:bi}. We compare $1$-parameter bifurcations to the one-dimensional internal dynamics case first (\Cref{thm:ddss,thm:ddhopf}) and turn to arbitrary $l$-parameter bifurcations after (\Cref{thm:ddlparam,thm:ddlparamgen,thm:ddlparamall}). The second part of this paper (\Cref{sec:ff}), is devoted to the application of these techniques and results to the class of feedforward networks. Generic $1$-parameter steady state bifurcations for such systems are described in full detail in \textcite{vonderGracht.2022} for one-dimensional internal dynamics. Here we use this knowledge to investigate generic steady state bifurcations for high-dimensional internal dynamics. Finally, we close this article with a discussion and outlook in~\Cref{sec:discussion}.

\section{Preliminaries: Hidden symmetries in homogeneous coupled cell systems}
\label{sec:pr}
In this section we briefly summarize the underlying theory of \emph{homogeneous coupled cell systems} in the language described in \textcite{Nijholt.2017d, Nijholt.2017c, Nijholt.2016, Rink.2013, Rink.2014, Rink.2015}.
To define a homogeneous network, we label the set of \emph{nodes} (or \emph{cells}) as $C=\lbrace p_1,\dotsc,p_N \rbrace$ and denote the interaction structure in the form of \emph{input maps} $\Sigma = \lbrace \sigma_1, \dotsc , \sigma_n \rbrace$. Then each $\sigma_i \colon C \to C$ characterizes one specific input type, i.e., cell $p$ receives an input from cell $\sigma(p)$ via an arrow of color $\sigma$.

Each cell has a state variable in the \emph{internal state space} which is the same finite dimensional real vector space for all cells: $v_p \in V$. 
The \emph{internal dynamics} of a cell is governed by the same function $f\colon V^n \to V$ with arguments given by its inputs (see \eqref{eq:netvf} below). 
Using the same internal phase space and the same governing function for each cell reflects homogeneity. The \emph{total phase space} is $\bigoplus_{p \in C} V \cong V^N$, for which we choose coordinates according to the cells of the network: \mbox{$v=(v_p)_{p\in C}=(v_{p_1},\dotsc,v_{p_N})^T$}. The dynamics is governed by the ordinary differential equations
\begin{equation}
	\label{eq:netvf}
	\dot{v} = \gamma_f(v) = \begin{pmatrix}
		f (v_{\sigma_1(p_1)}, \dotsc, v_{\sigma_n(p_1)}) \\
		f (v_{\sigma_1(p_2)}, \dotsc, v_{\sigma_n(p_2)}) \\
		\vdots \\
		f (v_{\sigma_1(p_N)}, \dotsc, v_{\sigma_n(p_N)})
	\end{pmatrix}.
\end{equation}
Indicating the inputs of $f$ for each cell by the input maps has the effect that a cell receives precisely one input of each type without imposing (symmetry) relations between the arguments of $f$, hence the term \emph{asymmetric inputs}. The \emph{network vector fields} are also referred to as \emph{admissible maps} or \emph{admissible vector fields}.

Linear admissible maps are of particular interest for the investigation of network dynamics. They appear as linearizations of network vector fields of the form \eqref{eq:netvf} at a steady state and are necessary for the determination of stability properties.
It was shown in Proposition 1.1 in \cite{vonderGracht.2022} that the space of linear admissible maps is spanned by linear maps $B_\sigma\colon V^N \to V^N$ with \mbox{$(B_{\sigma}(v))_p = v_{\sigma(p)}$} in the following sense: If $L\colon V^N \to V^N$ is linear and admissible, i.e., defined as in \eqref{eq:netvf} for some linear internal dynamics, then $L$ is of the form
\begin{equation}
	\label{eq:linadm}
	(Lv)_p = \sum_{\sigma \in \Sigma} b_\sigma (B_\sigma (v))_p,
\end{equation}
where $b_\sigma \in \gl{V}$ are linear maps on $V$ independent of $p$. In particular, if $V=\RR$ any linear admissible map $L$ is a linear combination of the $B_\sigma$, i.e.,
\begin{equation}
    \label{eq:linadmod}
    L = \sum_{\sigma \in \Sigma} b_\sigma B_\sigma
\end{equation}
for some $b_\sigma \in \RR$ using the identification $\gl{\RR}\cong\RR$.
\begin{remk}
	\label{rem:adjacency}
	After choosing a labeling of the cells, the maps $B_\sigma$ in the case $V=\RR$ can be interpreted as matrices that are also known as the \emph{adjacency matrices} of the network. Then $B_\sigma$ encodes the interaction structure of the input map $\sigma$. On the other hand, the linear maps can be interpreted as matrices with entries in $\gl{V}$ which have the same structure independent of the dimension of the internal dynamics. Consequently, we refer to them as (generalized) adjacency matrices.
	\hspace*{\fill}$\triangle$
\end{remk}

We make two additional assumptions on the set of input maps $\Sigma$. We want it to include the identity map $\sigma_1 = \Id \colon C \to C$ which is natural in the sense that we require each cell's dynamics to depend on its own state. Furthermore, we want $\Sigma$ to be closed under composition of maps. This means that every indirect input is also a direct input. Note that closedness can always be achieved by considering the closure -- the smallest set of input maps that is closed under composition and contains the original maps -- or, put differently by including concatenations of arrows in $\Sigma$. This is not a restriction as it only leads to an extension of the set of admissible vector fields $\gamma_f$. For more details on these assumptions, consult the aforementioned references.

Summarizing, we assume that $\Sigma$ has the structure of a monoid, which is to say a group without inverses or a semigroup with an identity. This algebraic property introduces \emph{hidden symmetry} to homogeneous coupled cell systems via the following two constructions. On the one hand, we may define a second network which has nodes labeled by the elements of $\Sigma$ and the same input maps $\Sigma$ as before. The inputs are now given by multiplication from the left: $\sigma \colon \Sigma \to \Sigma, \tau \mapsto \sigma\tau$. This network is called the \emph{fundamental network}. Its construction is essentially the same as for the left Cayley graph of $\Sigma$. The total phase space is $\bigoplus_{\sigma \in \Sigma} V \cong V^n$ with coordinates chosen accordingly again: $X = (X_\sigma)_{\sigma\in\Sigma} = (X_{\sigma_1}, \dotsc, X_{\sigma_n})$. The dynamics is governed by
\begin{equation*}
	\dot{X} = \Gamma_f(X) = \begin{pmatrix}
		f (X_{\sigma_1\sigma_1}, \dotsc, X_{\sigma_n\sigma_1}) \\
		f (X_{\sigma_1\sigma_2}, \dotsc, X_{\sigma_n\sigma_2}) \\
		\vdots \\
		f (X_{\sigma_1\sigma_n}, \dotsc, X_{\sigma_n\sigma_n})
	\end{pmatrix}.
\end{equation*}
Note that we use the same governing function $f$ and the same internal state space $V$ for both networks.

On the other hand, we can construct the regular representation $\sigma \mapsto A_\sigma$ of $\Sigma$ on $V^n$, denoted by $(V^n, A_\sigma)$, where we identify $\bigoplus_{\sigma \in \Sigma}V$ with $V^n$. The -- in general non-invertible -- linear maps for the action of $\Sigma$ are defined by multiplication from the right via $(A_\sigma X)_\tau= X_{\tau \sigma}$ for $\sigma, \tau \in \Sigma$. They satisfy the standard properties of representation maps, i.e., $A_{\Id} = \mathbbm{1}_{V^n}$ -- the identity on $V^n$ -- and $A_\sigma A_\tau = A_{\sigma \tau}$ for all $\sigma,\tau\in\Sigma$. The relation between the two constructions is that it can be seen that the class of \emph{equivariant} vector fields on $(V^n, A_\sigma)$ is precisely the same as the class of admissible vector fields $\Gamma_f$ for the fundamental network (Theorem 3.11 in \cite{Rink.2014}):
\begin{equation}
	\label{eq:equivariance}
	\left\{ \Gamma_f\ \left|\ f \in \mathcal{C}^\infty \left(\bigoplus_{\sigma\in\Sigma}V, V\right) \right. \right\} = \left\{ F \in \mathcal{C}^\infty (V^n, V^n) \mid F \circ A_\sigma = A_\sigma \circ F \text{ for all } \sigma \in \Sigma \right\}.
\end{equation}
Hence, the investigation of fundamental network vector fields can be performed using symmetry properties.
\begin{cor}
	\label{cor:adjacencyfund}
	Reformulating the definitions of adjacency matrices in terms of fundamental networks, we see that these are defined via multiplication from the left in $\Sigma$, i.e., $(B_{\sigma}(X))_\tau = X_{\sigma\tau}$. Additionally, we may compute that the $B_\sigma$ respect the multiplicative structure of $\Sigma$, now in the sense that $B_{\Id} = \mathbbm{1}_{V^n}$ and $B_\sigma B_\tau = B_{\tau\sigma}$. Due to \eqref{eq:equivariance}, these adjacency matrices span the space of all linear equivariant maps -- also called \emph{endomorphisms} -- as in \eqref{eq:linadm} and \eqref{eq:linadmod}.
\end{cor}

The contribution of the fundamental network to the analysis of the dynamics of the original network stems from the fact that, using technical tools called graph fibrations \textcite{DeVille.2015,DeVille.2015b}, one can show that the admissible vector fields of the two networks are semi-conjugate. In particular there exist linear maps $\pi_p\colon V^n \to V^N$ for each $p \in C$ -- non-invertible in general -- such that $\pi_p \circ \Gamma_f = \gamma_f \circ \pi_p$. Even more specifically, under the mild additional assumption that the original network contains a cell that receives an input from every other cell, the original network can be retrieved as a quotient network of the fundamental network. This means that there exists a \emph{balanced coloring} (see \textcite{Golubitsky.2006}) of the nodes of the fundamental network, such that identification of nodes of the same color gives rise to the original network. This has the effect that the total phase space of the original network -- and with that its dynamics -- is a synchrony subspace of the total phase space of the fundamental network. The condition to be fulfilled for this to be true -- sometimes also referred to as \emph{backward connectivity} \textcite{Aguiar.2019} -- is not very restrictive, as indirect inputs are to be considered as direct inputs in this setting. The semi-conjugacy between the dynamics of the two networks -- or the quotient relation -- together with the equivariance of the fundamental network vector fields \eqref{eq:equivariance} is what coins the term \emph{hidden symmetry}. For full details on these constructions see \textcite{Rink.2014,Nijholt.2016,Nijholt.2017c}.

In order to investigate dynamical properties such as bifurcations, we may now analyze the fundamental network using techniques from equivariant dynamics. In fact, many techniques from dynamics with underlying compact Lie group symmetry can be applied in a similar way (see \textcite{Rink.2014,Rink.2015,Nijholt.2017c}). In order for a bifurcation to occur, the steady state of a system of the form \eqref{eq:netvf} has to change its stability properties when one or multiple parameters vary. More precisely, one or multiple eigenvalues of the linearization $L$ have to cross the imaginary axis. 
The bifurcation then occurs along the \emph{generalized kernel} $\gker{L}$ or the \emph{center subspace} $\XX^c$ -- i.e., the generalized eigenspace to the critical eigenvalues -- of $L$. As the system is equivariant, it can be seen that the generalized kernel and center subspace are invariant under the action of $\Sigma$, they are so-called \emph{subrepresentations}.

On the other hand, it is known that any representation of a monoid uniquely (up to isomorphism) decomposes as a direct sum of \emph{indecomposable subrepresentations}. That means
\begin{equation}
	\label{eq:dec}
	V^n = W_1 \oplus \dotsb \oplus W_k,
\end{equation}
where $W_i$ for all $1 \le i \le k$ is a subspace satisfying $A_\sigma W_i \subset W_i$ for all $\sigma \in \Sigma$ and that cannot be decomposed further into non-trivial subspaces satisfying the same property. These indecomposable subrepresentations can be classified as being of \emph{real type} -- also called \emph{absolutely indecomposable} --, \emph{complex type} or \emph{quaternionic type}. This terminology is essentially identical to the one in the context of group representations. It indicates that the space of endomorphisms of the respective subrepresentation contains a real, complex, or quaternionic structure. However, these specific details are of no particular importance in the remainder of this article.

Summarizing, this means that the generalized kernel and center subspace are isomorphic to the direct sum of some of the indecomposable components: $\gker{L} \cong W_{i_1} \oplus \dotsb \oplus W_{i_s}$ with $1 \le i_1 < \dotso < i_s \le k$ and accordingly for $\XX^c$. Furthermore, one can determine which configurations of these components and their types can generically be present in an $l$-parameter bifurcation (see \textcite{Schwenker.2018} for the $1$-parameter case and \textcite{Nijholt.2017} for the $l$-parameter case). For example, in a $1$-parameter steady state bifurcation the generalized kernel is generically exactly one absolutely indecomposable subrepresentation. In a generic $1$-parameter Hopf bifurcation the center subspace is either the direct sum of two isomorphic components of real type or one component of complex or of quaternionic type. In general, for every fundamental network for which the decomposition of its regular representation into indecomposable components is known, we now have a systematic way of classifying the generic bifurcations and mode interactions in (multi-parameter) bifurcation problems.

In the remainder of this paper we frequently compare properties of dynamical systems with the same underlying network structure with one- and $d$-dimensional internal dynamics, where $d>1$. To avoid confusion when setting them side by side we introduce some notational conventions. In general, we distinguish these two settings by referring to them as the cases \od and \dd, respectively. Oftentimes, the distinction is highlighted in the notation by a super- or subscript $\mathbf{\mathsf{1}}$ (\od) or $\mathbf{\mathsf{D}}$ (\dd). We denote the internal phase spaces by $V=\Vo$ or $V=W \cong \RR^d$. If we do not want to specify one of the two cases -- i.e., if an observation holds for both cases simultaneously -- we keep on denoting the internal phase space by $V$. In particular, we have coordinates in the total phase spaces given by $(x_p)_{p\in C} \in \bigoplus_{p \in C} \Vo, (w_p)_{p\in C} \in \bigoplus_{p \in C} \Vd$ and $(v_p)_{p \in C} \in \bigoplus_{p \in C}V$. Note that in the case \od we commonly use the basis
\begin{equation}
\label{eq:basisod}
\left\{ (\delta_{p,q})_{p\in C} \right\}_{q\in C}
\end{equation}
for the total phase space $\bigoplus_{i = 1}^N \RR$, where $\delta_{p,q}$ is the Kronecker delta, which equals $1$ if $p=q$ and $0$ otherwise. The majority of the investigations in this paper focuses on fundamental networks. For these, we abbreviate the total phase spaces by $\No = \bigoplus_{\sigma\in\Sigma}\Vo$ and $\Nd = \bigoplus_{\sigma\in\Sigma} \Vd$. 
The coordinates on these spaces are denoted  by $x=(x_\sigma)_{\sigma\in\Sigma} \in \No, \omega=(w_\sigma)_{\sigma\in\Sigma}\in \Nd$ and $v=(v_\sigma)_{\sigma\in\Sigma}\in \bigoplus_{\sigma\in\Sigma}V$, instead of the capital letters $X_{\sigma}$ used above.
Finally, recall that the total phase space of the fundamental network is the representation space of the right regular representation of the monoid of input maps $\Sigma$. We denote the linear maps by which its elements act on $\No$ and $\Nd$ by $\Ao_\sigma$ and $\Ad_\sigma$, respectively. They are defined by
\[ (\Ao_\sigma x)_\tau = x_{\tau\sigma}, \qquad (\Ad_\sigma \omega)_\tau = w_{\tau\sigma} \]
in accordance to \eqref{eq:equivariance}. After choosing coordinates for both spaces that respect the network structure -- as above --, the representation maps can be interpreted as matrices. They have a very similar structure in both cases. The maps $\Ao_\sigma$ have one entry $1$ per row and all the other entries equal to $0$. The matrices $\Ad_\sigma$ have the same structure, but now the entries come from the ring $\mathfrak{gl}(\Vd)$ of linear maps on $\Vd$. 
Hence, whenever $\Ao_\sigma$ has an entry $1$, $\Ad_\sigma$ has an entry $\idvd$. Accordingly an entry $0$ in $\Ao_\sigma$ corresponds to an entry $0 \in \mathfrak{gl}(\Vd)$ in $\Ad_\sigma$. Compare also to the generalized adjacency matrices (\Cref{rem:adjacency}).

\section{Dimension reduction in fundamental networks}
\label{sec:hd}
In this section, we investigate the relation between networks with one-dimensional internal dynamics, i.e., $V \cong \RR$, and those with arbitrary finite-dimensional ($d$-dimensional), real vector spaces as internal phase spaces. This serves as a motivation for the restriction to the former case which is often done in bifurcation analysis of steady states. As was laid out in \Cref{sec:pr}, one typically restricts to the analysis of generic vector fields on the (generalized) kernel of the linearization at the bifurcation point, in order to qualitatively investigate steady state bifurcations. The technical tool for this is Lyapunov-Schmidt reduction. 
For more detailed results (e.g. stability properties of branching solutions) one computes the center manifold, which is a graph over the center subspace. These methods were thoroughly described in \textcite{Rink.2014,Nijholt.2017c,Nijholt.2019} for homogeneous coupled cell systems \eqref{eq:netvf}. 
In this section, we mostly restrict the investigations to fundamental networks. As was pointed out in \Cref{sec:pr}, this has the advantage that the reduction methods can be performed to respect symmetry: the generalized kernel and center subspace are subrepresentations and the reduced vector fields are precisely the equivariant ones on the respective subrepresentation. The bifurcation results for the original network can then be regained by restriction to a synchrony subspace. It turns out, that the generalized kernels and center subspaces in both cases (\od and \dd) are strongly related due to the monoid symmetry. In particular, for $1$-parameter steady state bifurcations, the case \dd leads to `the same' generalized kernels and center subspaces as one would obtain in the case \od. The bifurcation patterns are therefore qualitatively equivalent.

One of the main tools in this section is a representation of the total phase space in the case \dd in tensor product notation. To that end, we may identify
\begin{equation}
	\label{eq:tensor}
	\bigoplus_{p \in C} \Vd \cong \left( \bigoplus_{p \in C} \RR \right) \ox \Vd.
\end{equation}
In particular, we may split each $(w_p)_{p \in C} \in \bigoplus_{p \in C} \Vd$ as a sum of vectors that have precisely one non-vanishing coordinate entry: $(w_p)_{p \in C} = \sum_{q \in C} (\delta_{p,q} w_p)_{p \in C}$. This can be represented as a sum of pure tensors
\begin{eqnarray}
	\label{eq:tensorelt}
	\sum_{p \in C} (\delta_{p,q})_{q \in C} \ox w_p
\end{eqnarray}
using the basis for $\bigoplus_{p \in C} \RR$ as in \eqref{eq:basisod}. Note that this \emph{notation} as the sum of pure tensors is not unique, since multiple sums of pure tensors can represent the same element $(w_p)_{p\in C}\in \bigoplus_{p \in C}\Vd$. For example, 
\[ (\delta_{p,q})_{p \in C} \ox s w_p = s (\delta_{p,q})_{p \in C} \ox w_p \]
for a scalar $s \in \RR$. However, in critical cases like this the tensors are identified via an equivalence relation. Hence, the \emph{representation} via tensors in \eqref{eq:tensorelt} is unique. For more on this see the details and definitions of tensor products of vector spaces (e.g., in \cite{Halmos.1958}). This formalism allows for a graphical interpretation of attaching a vector space to each cell instead of using vectors of vectors and block matrices. It has been used, for example, in \textcite{Dias.2006,Golubitsky.2009,Leite.2006,Aguiar.2012} to describe adjacency matrices and their spectral properties for networks with higher-dimensional internal dynamics. Note that the isomorphism in \eqref{eq:tensor} is a mere identification of notations.
\begin{remk}
	\label{remk:tensorbasis}
	The representation \eqref{eq:tensorelt} relies on the fact that for two finite-dimensional vector spaces $A$ and $B$ with bases $\{ a_i \}_{i\in I}$ and $\{ b_j \}_{j\in J}$ the elements $\{ a_i \ox b_j \}_{i \in I, j\in J}$ form a basis for the tensor product $A \ox B$. In particular, given a basis $\{ b_1, \dotsc, b_d \}$ for $\Vd$, we can represent each element $(w_p)_{p\in C} \in \bigoplus_{p\in C} \Vd$ as
	\[ \sum_{p \in C} (\delta_{p,q})_{q \in C} \ox w_p = \sum_{p \in C} \sum_{i = 1}^d (\delta_{p,q})_{q \in C} \ox \alpha_i^p \cdot b_i = \sum_{p \in C} \sum_{i = 1}^d \alpha_i^p \cdot (\delta_{p,q})_{q \in C} \ox b_i, \]
	where $w_p = \alpha_1^p b_1 + \dotsb + \alpha_d^p b_d$.
	\hspace*{\fill}$\triangle$
\end{remk}

The tensor product formalism allows us to give a precise characterization of linear admissible maps in the case \dd. Recall from \eqref{eq:linadm} that any linear admissible map $L \colon \bigoplus_{p\in C}V \to \bigoplus_{p\in C}V$ can be represented using the generalized adjacency matrices $B_\sigma$ defined by $(B_\sigma (v_p)_{p\in C})_q = v_{\sigma(q)}$. Using the tensor notation, we see that these matrices from the case \od are sufficient to represent linear admissible maps also in the case \dd.
\begin{prop}
	\label{prop:linadmtensor}
	Define the \od linear map $\Bo_\sigma \colon \bigoplus_{p \in C} \RR \to \bigoplus_{p \in C} \RR$ through $(\Bo_\sigma (x_p)_{p\in C})_q = x_{\sigma(q)}$ for all $\sigma \in \Sigma$. Then any linear admissible map in the case \emph{\dd} $L \colon \bigoplus_{p \in C} \Vd \to \bigoplus_{p \in C} \Vd$ is of the form
	\begin{equation}
		\label{eq:linadmtensor}
		L = \sum_{\sigma \in \Sigma} \Bo_\sigma \ox b_\sigma
	\end{equation}
	for suitable linear maps $b_\sigma \in \gl{\Vd}$.
\end{prop}
\begin{proof}
	Let $L\colon \bigoplus_{p\in C}\Vd \to \bigoplus_{p\in C} \Vd$ be a linear admissible map. According to \eqref{eq:linadm}
    there are linear maps $b_\sigma \in \gl{\Vd}$ for each $\sigma \in \Sigma$ such that
	\[ (L((w_q)_{q\in C}))_p = \sum_{\sigma \in \Sigma} b_\sigma (\Bd_\sigma (w_q)_{q\in C})_p = \sum_{\sigma \in \Sigma} b_\sigma (w_{\sigma(p)}) \]
	in the non-tensor notation, where $(\Bd_\sigma (w_q)_{q\in C})_p = w_{\sigma(p)}$. In the tensor notation \eqref{eq:tensorelt}, $L((w_q)_{q\in C})$ can therefore be represented as
	\begin{equation}
		\label{eq:linadm1}
		L((w_q)_{q\in C}) = \sum_{p \in C} (\delta_{p,q})_{q \in C} \ox \sum_{\sigma \in \Sigma} b_\sigma (w_{\sigma(p)}) = \sum_{p \in C} \sum_{\sigma \in \Sigma} (\delta_{p,q})_{q \in C} \ox b_\sigma (w_{\sigma(p)}).
	\end{equation}
	On the other hand $\left( \Bo_\sigma (\delta_{p,q})_{q\in C} \right)_r = \delta_{p,\sigma(r)}$, which equals $1$ if $r \in \sigma^{-1}(p)$ and $0$ otherwise. In particular,
	\[ \Bo_\sigma (\delta_{p,q})_{q\in C} = \sum_{r \in \sigma^{-1}(p)} (\delta_{r,q})_{q\in C}. \]
	Using the tensor notation \eqref{eq:tensorelt}, i.e., representing $(w_q)_{q\in C}$ as
	\[ \sum_{p \in C} (\delta_{p,q})_{q \in C} \ox w_p, \]
	we compute
	\begin{align}
		\left[ \sum_{\sigma \in \Sigma} \Bo_\sigma \ox b_\sigma \right] \left( \sum_{p \in C} (\delta_{p,q})_{q \in C} \ox w_p \right)	&= \sum_{\sigma \in \Sigma} \sum_{p \in C} \Bo_\sigma (\delta_{p,q})_{q \in C} \ox  b_\sigma (w_p) \notag \\
			&= \sum_{\sigma \in \Sigma} \sum_{p \in C} \sum_{r \in \sigma^{-1}(p)} (\delta_{r,q})_{q\in C} \ox \left( b_\sigma w_{\sigma(r)} \right) \notag \\
			&= \sum_{\sigma \in \Sigma} \sum_{r \in C} (\delta_{r,q})_{q\in C} \ox \left( b_\sigma w_{\sigma(r)} \right). \label{eq:linadm2}
	\end{align}
	Therein the last equation holds since $\left\{ \sigma^{-1}(p) \mid p \in C \right\}$ forms a partition of $C$ for all $\sigma \in \Sigma$. As \eqref{eq:linadm1} and \eqref{eq:linadm2} agree, this completes the proof.
\end{proof}
\begin{remk}
	Note that the tensor notation is also applicable in the case \od. Then in \eqref{eq:tensorelt} and \eqref{eq:linadmtensor} we tensor with a scalar. Furthermore, $\gl{\RR}\cong\RR$ so that application of a linear map can be identified with scalar multiplication. This only plays a role in cases where we do not explicitly distinguish between \od and \dd.
	\hspace*{\fill}$\triangle$
\end{remk}

Furthermore, we may also describe robust synchrony subspaces in the tensor notation. It is well known that these are in one-to-one correspondence to balanced partitions of the cells of the network (see for example \textcite{Field.2004,Golubitsky.2006} for an exact definition).
\begin{prop}
	\label{prop:tensorsynch}
	Let $P = \{P_1, \dotsc, P_r \}$ be a balanced partition of the cells $C$ and let
	\begin{align*}
		\Delo_P	&= \left\{ \left. (x_p)_{p\in C} \in \bigoplus_{p \in C} \Vo\ \right|\ x_p = x_q, \quad \text{if} \quad p,q \in P_i \quad \text{for some } 1 \le i \le r \right\}, \\
		\Deld_P	&= \left\{ \left. (w_p)_{p\in C} \in \bigoplus_{p \in C} \Vd\ \right|\ w_p = w_q, \quad \text{if} \quad p,q \in P_i \quad \text{for some } 1 \le i \le r \right\}
	\end{align*}
	be the corresponding synchrony subspaces for one- and for high-dimensional internal dynamics, respectively. Then
	\[ \Deld_P \cong \Delo_P \ox \Vd . \]
\end{prop}
\begin{proof}
	The result follows almost directly from the characterization of bases of tensor products in \Cref{remk:tensorbasis}. The synchrony subspace $\Delo_P$ is spanned by elements $\{ (x^1_p)_{p\in C}, \dotsc, (x^r_p)_{p\in C} \}$, where \mbox{$x^i_p=1$} if $p \in P_i$ and $x^i_p=0$ otherwise. Hence,
	\begin{equation}
		\label{eq:tensorsynch1}
		\left\{ \left. (x^i_p)_{p\in C} \ox b_j\ \right|\ i=1, \dotsc,r \text{ and } j=1,\dotsc, d  \right\}
	\end{equation}
	is a basis of $\Delo_P \ox \Vd$, when $\{ b_1, \dotsc,b_d \}$ is a basis of $\Vd$.
	
	On the other hand, let $p_1, \dotsc, p_r \in C$ be a set of representatives of the partition $P$, i.e., $p_i \in P_i$. Then every element $(w_p)_{p\in C} \in \Deld_P$ can be represented as
	\[ (w_p)_{p\in C} = \sum_{i = 1}^r (x^i_{p} \cdot w_{p_i})_{p\in C} \]
	using the basis of $\Delo_P$, since $w_p = w_{p_i}$ if $p \in P_i$. Furthermore, every element $w_{p_i}$ is of the form
	\[ w_{p_i} = \sum_{j = 1}^d \alpha^i_j \cdot b_j \]
	using the basis of $\Vd$. Hence, we obtain
	\[ (w_p)_{p\in C} = \sum_{i = 1}^r \sum_{j = 1}^d \alpha^i_j \cdot (x^i_{p} \cdot b_j)_{p\in C}. \]
	In particular, we see that
	\begin{equation}
		\label{eq:tensorsynch2}
		\left\{ \left. (x^i_p \cdot b_j)_{p\in C}\ \right|\ i=1, \dotsc,r \text{ and } j=1,\dotsc, d  \right\}
	\end{equation}
	is a basis of $\Deld_P$. Representing these basis elements in their respective tensor notation \eqref{eq:tensorelt} shows that \eqref{eq:tensorsynch2} agrees with \eqref{eq:tensorsynch1} which completes the proof.
\end{proof}

From now on, we focus on fundamental networks. We want to understand the structure of the regular representation and its decomposition into indecomposable subrepresentations \eqref{eq:dec} to determine all possible generalized kernels when the internal phase space has dimension greater than $1$. As it turns out, this is strongly related to the case \od.

The main goal of this part is to prove \Cref{thm:dechd}, which states that the indecomposable subrepresentations in the decomposition of the total phase space \eqref{eq:dec} are independent of the dimension of the internal phase space. In fact, only their multiplicities are varied by increasing this dimension. The proof of this result relies on multiple technical steps in the upcoming propositions and lemmas that aim at making the intuition that the tensor representation suggests rigorous. We first prove that the right regular representations of the monoid $\Sigma$ on $\No$ and on $\Nd$ are isomorphic. Then, we need to ensure that the tensor representation also behaves well when decomposing $\No$ into indecomposable components. The result for $\Nd$ then follows by simultaneously decomposing $\No$ into indecomposable components and $\Vd$ into one-dimensional subspaces spanned by single basis vectors.
\begin{prop}
	The $\Sigma$-representation $\lbrace \Ad_\sigma \rbrace_{\sigma\in\Sigma}$ on $\Nd$ is isomorphic to $\No \ox \Vd$ on which $\sigma \in \Sigma$ acts as $\Ao_\sigma \ox \idvd$.
\end{prop}
\begin{proof}
	The main idea for the proof is the interpretation of the total phase space $\Nd$ as having the vector space $\Vd$ attached to each cell of the network. These are in 1-to-1 correspondence to the coordinates in $(x_\sigma)_{\sigma\in\Sigma} \in \No = \bigoplus_{\sigma \in \Sigma} \RR$. Hence, we assign a vector $w_\sigma \in \Vd$ to each coordinate $x_\sigma$ which is reflected in the tensor notation $\No \ox \Vd$. 
	
	First, note that both vector spaces have dimension $n \cdot \dim \Vd$. Hence, they are isomorphic as such. An isomorphism can be defined as
	\begin{align*}
		\Phi \colon \No \ox \Vd &\to \Nd \\
		(x_\sigma)_{\sigma \in \Sigma} \ox w &\mapsto (x_\sigma w)_{\sigma \in \Sigma},
	\end{align*}
	which is linearly extended to non-pure tensors (sums of elements $(x_\sigma)_{\sigma \in \Sigma} \ox w$). Note that this, except for the transposition, coincides with the often used identification of $(x_\sigma)_{\sigma \in \Sigma} \ox w$ with the outer product $(x_\sigma)_{\sigma \in \Sigma} w^T = (x_\sigma w^T)_{\sigma \in \Sigma}$.
	
	As in \eqref{eq:tensorelt}, we may uniquely split each $\omega = (w_\sigma)_{\sigma\in\Sigma} \in \Nd$ as a sum of vectors that have precisely one non-vanishing coordinate entry: $(w_\sigma)_{\sigma\in\Sigma} = \sum_{\tau \in \Sigma} (\delta_{\sigma,\tau} w_\sigma)_{\sigma\in\Sigma}$. Hence, the map 
	\[ \Psi \colon (w_\sigma)_{\sigma\in\Sigma} \mapsto \sum_{\tau \in \Sigma} (\delta_{\sigma,\tau})_{\sigma\in\Sigma} \ox w_\tau \]
	is inverse to $\Phi$. In particular,
	\begin{align*}
		\Psi \left( \Phi ((x_\sigma)_{\sigma \in \Sigma} \ox w) \right) &= \Psi \left( (x_\sigma w)_{\sigma \in \Sigma} \right) \\
			&= \Psi \left( \sum_{\tau \in \Sigma} (\delta_{\sigma,\tau} x_\sigma w)_{\sigma\in\Sigma} \right) \\
			&= \sum_{\tau \in \Sigma} (\delta_{\sigma,\tau})_{\sigma\in\Sigma} \ox (x_\sigma w) \\
			&= \sum_{\tau \in \Sigma} (\delta_{\sigma,\tau} v_\sigma)_{\sigma\in\Sigma} \ox w \\
			&= (x_\sigma)_{\sigma \in \Sigma} \ox w,
	\end{align*}
	which extends linearly to non-pure tensors.	Recall that the basis $\{ (\delta_{\sigma,\tau})_{\sigma\in\Sigma} \}_{\tau\in\Sigma} \subset \No$ corresponds to the cells of the network. Therefore, we may interpret $\Psi$ as the map that picks the vector $w_\sigma$ in the $\sigma$-entry and attaches it to cell $\sigma$ via the tensor product.
	
	It remains to be checked that $\Phi$ intertwines the two $\Sigma$-representations. In order to do so, we compute
	\[ \Phi\left( [\Ao_\tau \ox \idvd] ((x_\sigma)_{\sigma\in\Sigma} \ox w) \right) = \Phi ((x_{\sigma\tau})_{\sigma\in\Sigma} \ox w) = (x_{\sigma\tau} w)_{\sigma\in\Sigma} \]
	but also
	\[ \Ad_\tau \Phi((x_\sigma)_{\sigma\in\Sigma} \ox w) = \Ad_\tau(x_\sigma w)_{\sigma\in\Sigma} = (x_{\sigma\tau} w)_{\sigma\in\Sigma}. \]
	Equivariance on non-pure tensors follows from linearity of the representation matrices. This proves equivalence of the representations.
\end{proof}
\begin{remk}
	\label{rem:linadmfundtensor}
	Combining \Cref{cor:adjacencyfund} and \Cref{prop:linadmtensor} we see that any endomorphism \mbox{$L \in \End{\No \ox \Vd}$} of the regular representation that characterizes the fundamental network is of the form \eqref{eq:linadmtensor} in the tensor notation, i.e.,
	\[ L = \sum_{\sigma \in \Sigma} B_\sigma \ox b_\sigma \]
	for linear maps $b_\sigma \in \gl{W}$.
	\hspace*{\fill}$\triangle$
\end{remk}

The tensor notation relates the representation $\Nd$ in a straightforward way to the representation $\No$ of the fundamental networks in the cases \dd and \od, respectively. Understanding the structure of the representation, especially its decomposition into subrepresentations, is essential for the investigation of generic dynamics. 
We now relate the decomposition of $\No \ox \Vd$ to that of $\No$. It is well known that the tensor product commutes with direct sums: if $A$ and $B$ are vector spaces with $A=A_1 \oplus A_2$, then $A \ox B = (A_1 \ox B) \oplus (A_2 \ox B)$ (see for example Theorem 17 in \cite{Dummit.2004}). In particular, projections $\pi_1, \pi_2 \in \mathfrak{gl}(A)$ onto $A_1$ and $A_2$, respectively, can be extended to projections $\pi_1 \ox \mathbbm{1}_{B}, \pi_2 \ox \mathbbm{1}_{B} \in \mathfrak{gl}(A\ox B)$ onto $(A_1 \ox B), (A_2 \ox B)$, respectively. Applied to the above setting, this even applies to decompositions into subrepresentations: If $\No = Y_1 \oplus Y_2$ where $Y_1$ and $Y_2$ are subrepresentations, the projections $\pi_1$ and $\pi_2$ are equivariant with respect to $\lbrace \Ao_\sigma \rbrace_{\sigma\in\Sigma}$. Then $\pi_1 \ox \idvd$ and $\pi_2 \ox \idvd$ are equivariant with respect to $\lbrace \Ao_\sigma \ox \idvd \rbrace_{\sigma\in\Sigma}$. Thus, $\No \ox \Vd = (Y_1 \ox \Vd) \oplus (Y_2 \ox \Vd)$ as a decomposition into subrepresentations. The same constructions works for a decomposition of $\Vd$.
\begin{lemma}
	\label{lem:subrepDD}
	Suppose $Y \subset \No$ is a subrepresentation with respect to $\lbrace \Ao_\sigma \rbrace_{\sigma\in\Sigma}$ and let $w \in \Vd \setminus \{0\}$ be an arbitrary element. Then $Y \ox \langle w \rangle \subset \No \ox V$ is a subrepresentation with respect to $\lbrace \Ao_\sigma \ox \idvd \rbrace_{\sigma\in\Sigma}$. Furthermore, $Y \ox \langle w \rangle \cong Y$ as subrepresentations.
\end{lemma}
\begin{proof}
	We first note
	\begin{align*}
	y_1 \ox r_1 w + \dotso + y_k \ox r_k w	&= r_1 y_1 \ox w + \dotso + r_k y_k \ox w \\	
											&= (r_1 y_1 + \dotso + r_k y_k) \ox w
	\end{align*}
	for an arbitrary formal sum in $Y \ox \langle w \rangle$, i.e., $y_1, \dotsc, y_k \in Y$ and $r_1, \dotsc, r_k \in \RR$. Hence, every element in $Y \ox \langle w \rangle$ can uniquely be expressed as a pure tensor $y \ox w$ with $y = (r_1 y_1 + \dotso + r_k y_k) \in Y$. Therefore, we may identify $y \ox w$ with $y$ which is equivariant by definition of the representation maps. This proves the claim.
\end{proof}
\begin{remk}
	In particular, if $Y$ is indecomposable of a specific type, the same holds for $Y \ox \langle w \rangle$, due to equivalence of the representations.
	\hspace*{\fill}$\triangle$
\end{remk}

We obtain the main result of this section as a corollary of the above.
\begin{theorem}
	\label{thm:dechd}
	Suppose $\No = Y_1 \oplus \dotsb \oplus Y_s$ as a decomposition into indecomposable subrepresentations with respect to $\lbrace \Ao_\sigma \rbrace_{\sigma\in\Sigma}$. Let $\lbrace b_1, \dotsc, b_d \rbrace$ be a basis for $\Vd$. Then
	\begin{subequations}
		\label{eq:1DDD}
		\begin{align}
			\Nd \cong \No \ox \Vd	&= \bigoplus_{i = 1}^s \bigoplus_{j=1}^d Y_i \ox \langle b_j \rangle \label{eq:1DDDa} \\
									&\cong \bigoplus_{i=1}^s Y_i^d	\label{eq:1DDDb}
		\end{align}
	\end{subequations}	
	as a decomposition into indecomposable subrepresentations with respect to $\lbrace A_\sigma \ox \idvd \rbrace_{\sigma\in\Sigma}$.
\end{theorem}
\begin{remk}
	Note that the first isomorphism in \eqref{eq:1DDDa} is only an identification of different notations. The second relation however is an equality. Hence, we may identify the fundamental network representation space $\Nd$ in the case \dd with the decomposition given in that equation without changing the coordinates. In particular the network structure is preserved by this identification. On the other hand, not every decomposition of $\Nd$ is of the form \eqref{eq:1DDDa}. However, as the decomposition into indecomposable subrepresentations of a monoid representation is unique up to isomorphisms, every indecomposable component $W_i \subset \Nd$ is isomorphic to one of $\No$, i.e., $W_i \cong Y_i \ox \langle w \rangle \cong Y_i$ after relabeling the indices, but they are not necessarily equal. Nonetheless, every decomposition of $\No$ gives rise to a decomposition of $\Nd$ as in \Cref{thm:dechd}.
	\hspace*{\fill}$\triangle$
\end{remk}
\begin{cor}
	\label{cor:hddecomp}
	Let $W^1$ and $W^2$ be two finite-dimensional real vector spaces that we choose as internal phase spaces of a fundamental network. Furthermore, assume $\dim{W^1} \le \dim{W^2}$. Then there is a subrepresentation $U \subset \bigoplus_{\sigma\in\Sigma} W^2$ such that
	\[ U \cong \bigoplus_{\sigma\in\Sigma}W^1. \]
\end{cor}
\begin{proof}
	This follows directly from \eqref{eq:1DDDb}.
\end{proof}

Finally, we show that the identification of indecomposable subrepresentations in the cases \od and \dd respects synchrony subspaces. This is particularly relevant in the analysis of bifurcations, as it implies that the restriction to the case \od does not change patterns of synchrony. In particular, this holds true for those patterns of synchrony that provide the original network as a quotient of the fundamental network. We need the following technical result on subspaces of tensor product spaces. The proof is elementary. It is included here, because we could not find a suitable reference in the literature.
\begin{lemma}
	\label{lem:intersectionsubsp}
	Let $A$ and $B$ be finite-dimensional real vector spaces and let $A_1, A_2 \subset A$ and \mbox{$B_1, B_2 \subset B$} be subspaces. Then $(A_1 \ox B_1) \cap (A_2 \ox B_2) = (A_1 \cap A_2) \ox (B_1 \cap B_2)$ as a subspace of $A \ox B$.
\end{lemma}
\begin{proof}
	This result follows from the representation of a basis of the tensor product space in terms of bases of the components in \Cref{remk:tensorbasis}. Let $\{ a_i \}_{i \in I}$ and $\{ b_j \}_{j \in J}$ be bases of $A$ and $B$, respectively, and let $I_1, I_2 \subset I$ and $J_1, J_2 \subset J$ be subsets such that
	\begin{alignat*}{4}
		& A_1	&&= \left\langle a_i\ \left|\ i \in I_1 \right. \right\rangle, \qquad	&& A_2	&&= \left\langle a_i\ \left|\ i \in I_2 \right. \right\rangle \\
		& B_1	&&= \left\langle b_j\ \left|\ j \in J_1 \right. \right\rangle, \qquad	&& B_2	&&= \left\langle b_j\ \left|\ j \in J_2 \right. \right\rangle.
	\end{alignat*}
	Note that $I, I_1$ and $I_2$ can be constructed by completing a basis of $A_1 \cap A_2$ to bases of $A_1$ and $A_2$. Then the set of all basis elements of $A_1$ and $A_2$ is completed to a basis of $A$. Accordingly we construct $J, J_1$, and $J_2$ for $B$ from a basis of $B_1\cap B_2$. In particular, we obtain
	\[ A_1 \cap A_2	= \left\langle a_i\ \left|\ i \in I_1 \cap I_2 \right. \right\rangle, \qquad B_1 \cap B_2 = \left\langle b_j\ \left|\ j \in J_1 \cap J_2 \right. \right\rangle. \]
	Hence, using \Cref{remk:tensorbasis} we see
	\[ (A_1 \cap A_2) \ox (B_1 \cap B_2) = \left\langle a_i \ox b_j\ \left|\ i \in I_1 \cap I_2 \text{ and } j \in J_1 \cap J_2 \right. \right\rangle. \]
	On the other hand
	\[ A_1 \ox B_1 = \left\langle a_i \ox b_j\ \left|\ i \in I_1 \text{ and } j \in J_1 \right. \right\rangle, \qquad (A_2 \ox B_2) = \left\langle a_i \ox b_j\ \left|\ i \in I_2 \text{ and } j \in J_2 \right. \right\rangle. \]
	Thus,
	\[ (A_1 \ox B_1) \cap (A_2 \ox B_2) = \left\langle a_i \ox b_j\ \left|\ i \in I_1 \cap I_2 \text{ and } j \in J_1 \cap J_2 \right. \right\rangle \]
	which completes the proof.
\end{proof}
\noindent
As a corollary of \Cref{prop:tensorsynch} and \Cref{lem:intersectionsubsp} we obtain
\begin{prop}
	\label{prop:ddodsynch}
	Let $P = \{P_1, \dotsc,P_r\}$ be a balanced partition of the cells $\Sigma$ of the fundamental network and let $\Delo_P$ and $\Deld_P$ denote the corresponding robust synchrony subspaces in the case \od and \dd, respectively. Furthermore, let $Y\subset \No$ be an indecomposable component and $U \subset \Nd$ such that \mbox{$U \cong Y \ox \langle w \rangle \cong Y$}. Then
	\[ \Deld_P \cap U \cong (\Delo_P \ox \Vd) \cap (Y \ox \langle w \rangle) = (\Delo_P \cap Y) \ox \langle w \rangle. \]
\end{prop}

\section{Implications for bifurcations of steady states}
\label{sec:bi}
In \Cref{thm:dechd}, we describe the relation between the algebraic structures of a given fundamental network in the cases \od and \dd. In particular, decomposing the regular representation in the case \od provides a decomposition in the case \dd by choosing a basis for $\Vd$. In this section, we want to investigate how this allows us to reduce the investigation of bifurcations in fundamental networks with high-dimensional internal dynamics to that in fundamental networks with one-dimensional internal dynamics. We provide the general setting first and discuss the implications of \Cref{thm:dechd} on generic steady state and generic Hopf bifurcations in \Cref{subsec:is,subsec:ih,subsec:il,subsec:in}. Throughout this section, we still focus on fundamental networks.

In the investigation of bifurcations of steady states, one is interested in qualitative changes of the set of steady state or periodic solutions when a given steady state changes its stability properties as a parameter is varied. We assume that the admissible vector fields \eqref{eq:netvf} depend on a real parameter $\lambda \in \RR^l$:
\begin{equation}
	\label{eq:netvfparam}
	\Gamma_f(v,\lambda) = \begin{pmatrix}
		f (v_{\sigma_1\sigma_1}, \dotsc, v_{\sigma_n\sigma_1}, \lambda) \\
		\vdots \\
		f (v_{\sigma_1\sigma_n}, \dotsc, v_{\sigma_n\sigma_n}, \lambda)
	\end{pmatrix}.
\end{equation}
The case $l>1$ is also referred to as an $l$-parameter bifurcation in which we interpret each component of $\lambda$ as one parameter. We aim at describing generic bifurcations from a fully synchronous steady state. Without loss of generality, we may assume this to be the origin and the bifurcation to occur for $\lambda=0$. Hence, we assume 
\[ \Gamma_f(0,0) = 0, \]
which implies $f(0,0) = 0$. Furthermore, the assumption that the steady state changes its stability can be translated to certain technical conditions on the partial derivatives of $f$. We are interested in steady states and periodic solutions close to this bifurcation point for a generic smooth function $f$ satisfying these conditions.

More precisely, due to the implicit function theorem and the Hopf bifurcation theorem the linearization $D_v\Gamma_f(0,0)$ needs to have eigenvalues on the imaginary axis in order for non-trivial (steady state or periodic) solutions to exist close to the bifurcation point. A steady state bifurcation requires vanishing eigenvalues and a Hopf bifurcation requires purely imaginary eigenvalues of the linearization at the bifurcation point. As the vector field \eqref{eq:netvfparam} is equivariant in its $v$-component and the steady state is fully synchronous -- which is equivalent to fully symmetric --, so is the linearization $D_v\Gamma_f(0,0)$. Hence, in either case it induces a decomposition of the right regular representation into subrepresentations given by its generalized kernel and reduced image or by its center and hyperbolic subspaces, respectively, i.e.,
\begin{align*} 
	\bigoplus_{\sigma\in\Sigma} V &= \gker{D_v\Gamma_f(0,0)} \oplus \redim{D_v\Gamma_f(0,0)} \quad \text{or} \\
	\bigoplus_{\sigma\in\Sigma} V &= \Xc \oplus \Xh
\end{align*}
Therein $\gker{D_v\Gamma_f(0,0)}$ is the generalized eigenspace of the eigenvalue $0$ and $\Xc$ the direct sum of the generalized eigenspaces of the eigenvalues on the imaginary axis. The respective complements $\redim{D_v\Gamma_f(0,0)}$ and $\Xh$ are given by the direct sum of the remaining generalized eigenspaces.

The bifurcation problem can be reduced to an equivalent one on the generalized kernel or the center subspace. In particular, branching steady state or periodic solutions of the original system lie on the center manifold which is a graph over the corresponding subspace. We refer to this fact by saying that the bifurcation occurs along that generalized kernel or center subspace. Furthermore, the reduced bifurcation problem is fully characterized by symmetry of the subrepresentation. Hence, in order to classify the bifurcations of steady states that may occur in the fundamental network, one has to determine all possible subrepresentations of $\bigoplus_{\sigma\in\Sigma} V$ that can form center subspaces. As the decomposition into indecomposable subrepresentations \eqref{eq:dec} is unique up to isomorphisms, that means one needs to decompose the regular representation
\[ \bigoplus_{\sigma\in\Sigma} V = W_1 \oplus \dotsb \oplus W_k. \]
The center subspace is then isomorphic to the direct sum of a suitable collection of the components $W_i$. \Cref{thm:dechd} shows that it is sufficient to determine this decomposition in the case \od\ -- i.e., to decompose $\No$ --, as each indecomposable component of $\Nd$ is also an indecomposable component of $\No$. In particular, if we know the decomposition of $\No$ and if the dimension of the internal phase space $\Vd$ (in the case \dd) is $d$, we immediately obtain the decomposition of $\Nd$ in the form of $d$ copies of the components of $\No$.

Furthermore, recall that indecomposable representations come in three types: real (also called absolutely indecomposable), complex and quaternionic. This is determined by technical properties of the linear equivariant maps of that representation. We make the decomposition more precise by writing
\[ \bigoplus_{\sigma\in\Sigma} V = V_1^\mathrm{R} \oplus \dotsb \oplus V_{m_\mathrm{R}}^\mathrm{R}\oplus V_1^\mathrm{C} \oplus \dotsb \oplus V_{m_\mathrm{C}}^\mathrm{C} \oplus V_1^\mathrm{H} \oplus \dotsb \oplus V_{m_\mathrm{H}}^\mathrm{H}, \]
where
\[ V_i^\mathrm{R} \cong \left( W_i^\mathrm{R} \right)^{s_i^\mathrm{R}}, \quad V_i^\mathrm{C} \cong \left( W_i^\mathrm{C} \right)^{s_i^\mathrm{C}}, \quad V_i^\mathrm{H} \cong \left( W_i^\mathrm{H} \right)^{s_i^\mathrm{H}} \]
collect isomorphic subrepresentations. The $W_i^\mathrm{R}, W_i^\mathrm{C}$, and $W_i^\mathrm{H}$ are pairwise non-isomorphic indecomposable subrepresentations of real, complex, and quaternionic type, respectively. Then we may determine, which configurations of indecomposable components are possible as generalized kernels or center subspaces for a generic $l$-parameter family of equivariant vector fields. In particular, let
\[ U \cong \bigoplus_{i = 1}^{m_\mathrm{R}} \left( W_{i}^\mathrm{R} \right)^{\rho_i} \oplus \bigoplus_{i = 1}^{m_\mathrm{C}} \left( W_{i}^\mathrm{C} \right)^{\gamma_i} \oplus \bigoplus_{i = 1}^{m_\mathrm{H}} \left( W_{i}^\mathrm{H} \right)^{\iota_i} \]
with $0 \le \rho_i \le s_i^\mathrm{R}, 0 \le \gamma_i \le s_i^\mathrm{C}, 0 \le \iota_i \le s_i^\mathrm{H}$ for every $i$. Then $U$ can only occur as a generalized kernel -- i.e., $\gker{D_v\Gamma_f(0,0)} \cong U$ --, if
\begin{equation}
	\label{eq:condker}
	K_U = \sum_{i = 1}^{m_\mathrm{R}} \rho_i + 2 \cdot \sum_{i = 1}^{m_\mathrm{C}} \gamma_i + 4 \cdot \sum_{i = 1}^{m_\mathrm{H}} \iota_i \le l,
\end{equation}
and as a center subspace -- i.e., $\Xc \cong U$ --, if
\begin{equation}
	\label{eq:condcent}
	C_U = \sum_{i = 1}^{m_\mathrm{R}} \lceil \rho_i/2 \rceil + \sum_{i = 1}^{m_\mathrm{C}} \gamma_i + \sum_{i = 1}^{m_\mathrm{H}} \iota_i \le l.
\end{equation}
Here $\lceil a \rceil$ denotes the nearest larger integer. For details on this see \textcite{Schwenker.2018} for the case $l=1$ and \textcite{Nijholt.2017} for the general case. The conditions \eqref{eq:condker} and \eqref{eq:condcent} have simpler interpretations for the case of $1$-parameter families.

In the upcoming subsections, we investigate the interplay of \Cref{thm:dechd} with the classification of generalized kernels and center subspaces in generic bifurcation problems. This allows to relate generic bifurcations of a fixed fundamental network with high-dimensional internal dynamics to those of the same network with one-dimensional internal dynamics. We begin with a discussion of the $1$-parameter case, before turning to general $l$-parameter families.

\subsection[Generic \texorpdfstring{$1$}{1}-parameter steady state bifurcations]{Generic $1$-parameter steady state bifurcations}
\label{subsec:is}
When focusing on steady state bifurcations, we want to characterize solutions to
\[ \Gamma_f (v,\lambda) = 0 \] 
close to the bifurcation point $(v_0,\lambda_0) = (0,0)$. We may restrict the analysis to the generalized kernel $\gker{D_v\Gamma_f(0,0)}$ in order to qualitatively determine branching steady states. By Lyapunov-Schmidt reduction, all generic branches of steady states can be found in that kernel. The actual branches (and their stability properties) still require investigation of the center manifold, the qualitative picture including branching directions, however, does not. As the generalized kernel is a complementable subrepresentation (i.e., there exists an invariant complement) and the Lyapunov-Schmidt reduction can be performed to respect equivariance (see \textcite{Rink.2014}), this allows us to restrict to a lower-dimensional equivariant bifurcation problem on the generalized kernel. From \eqref{eq:condker}, we know that this generalized kernel is generically an absolutely indecomposable subrepresentation.

Let us turn to the case \dd, i.e., to a bifurcation problem on $\Nd$. Applying \Cref{thm:dechd} (especially \eqref{eq:1DDDb}), we obtain that the generalized kernel in a given network in the case \dd is isomorphic as a subrepresentation to one of the indecomposable subrepresentations one computes for $\No$ in the case \od. That is
\[ \gker{D_\omega\Gamma_f(0,0)} \cong Y_i \subset \No \]
in the notation of \Cref{thm:dechd}. As the dynamics -- in particular the generic steady state bifurcations -- on these subrepresentations is entirely classified by monoid symmetry, the reduced bifurcation problem in the case is equivalent to one on the subrepresentation $Y_i$. Hence, the generic steady state bifurcations in the case \dd occur generically also in the case \od. On the other hand, as every component $Y_i$ can occur as a generalized kernel $\gker{D_\omega\Gamma_f(0,0)}$ in this way, any generic steady state bifurcation in the case \od occurs generically in the case \dd as well. Summarizing, we have shown
\begin{theorem}
	\label{thm:ddss}
	The generic $1$-parameter steady state bifurcations in a fundamental network with $d$-dimensional internal dynamics are qualitatively the same as those for the same network with $1$\nobreakdash-dimensional internal dynamics in the sense that the reduced bifurcation problems are equivalent in both cases.
\end{theorem}
\noindent
We cannot expect a more precise comparative result. The generalized kernels are the same in both cases so that the reduced bifurcation problems are the same. However, the full system has different dimensions and requires different coordinate systems. Therefore, the branching solutions for the full systems in general cannot be `equal' in a stricter sense.

\subsection[Generic \texorpdfstring{$1$}{1}-parameter Hopf bifurcations]{Generic $1$-parameter Hopf bifurcations}
\label{subsec:ih}
Similar to \Cref{subsec:is}, we investigate Hopf bifurcations in a generic $1$-parameter family of fundamental network vector fields. In order to do so, we have to investigate possible center subspaces corresponding to non-vanishing purely imaginary eigenvalues. Condition \eqref{eq:condcent} shows that only three cases can occur generically:
\[ \Xc \cong W_i^\mathrm{C}, \quad \Xc \cong W_i^\mathrm{H}, \quad \Xc \cong \left( W_i^\mathrm{R} \right)^2. \]
That is, either $\Xc$ is isomorphic to precisely one indecomposable component of complex or of quaternionic type or it is the direct sum of two isomorphic components of real type. Note that $\Xc \cong W_i^\mathrm{R}$ -- i.e., the center subspace is isomorphic to one indecomposable component of real type -- satisfies \eqref{eq:condcent} as well. However, this situation can only occur for an eigenvalue $0$ and cannot provide a pair of purely imaginary eigenvalues (for more details see \textcite{Nijholt.2017}). In the case \dd, \Cref{thm:dechd} shows that the existence of two isomorphic components of real type in $\Nd$ can occur in two different situations. Either $\No$ contains two isomorphic components of real type -- i.e., \mbox{$\Xc \cong Y_i^2 \subset \No$} in the notation of \Cref{thm:dechd} -- or $d \ge 2$ and we obtain two copies of the same component of $\No$ due to the high-dimensional internal dynamics -- i.e., $\Xc \cong Y_i^2 \not\subset \No$. Moreover, the latter choice is the only center subspace for a generic $1$-parameter Hopf bifurcation in the case \dd that does not occur in the case \od. The other three cases are possible independent of $d$.

The solutions to the reduced bifurcation problem on $\Xc$ are entirely classified by symmetry. In particular, in the cases that are independent of $d$ the reduced bifurcation problem is equivalent for each choice of $d$. The bifurcations are qualitatively the same as in the case \od. Furthermore, note that these cases describe all possible center subspaces in a generic $1$-parameter Hopf bifurcation in the case \od. Hence, all generic $1$-parameter Hopf bifurcations in the case \od can also be observed generically in the case \dd. Conversely, the center subspace that is due to the high-dimensional internal dynamics -- i.e., $\Xc \cong Y_i^2 \not\subset \No$ -- can only occur as a generic center subspace for \mbox{$d \ge 2$}. Hence, in general there is no equivalent reduced bifurcation problem in the case \od and the corresponding Hopf bifurcations can only be observed in the case \dd. The discussion of this subsection can be summarized as
\begin{theorem}
	\label{thm:ddhopf}
	\begin{enumerate}[label=(\roman*)]
		\item All center subspaces in generic $1$-parameter Hopf bifurcations in the case \od are generic as center subspaces in a $1$-parameter Hopf bifurcation in the case \dd for any $d$. The branching periodic solutions corresponding to one center subspace are qualitatively the same for all values of $d$ in the sense that the reduced bifurcation problems are equivalent.
		\item Let $\No = Y_1 \oplus \dotsb \oplus Y_s$ be a decomposition into indecomposable subrepresentations and assume $Y_i$ to be of real type such that $Y_i \not\cong Y_j$ for all $j \ne i$. If the internal dynamics is at least two-dimensional, i.e., $d \ge 2$, $Y_i$ yields a center subspace $\Xc \cong Y_i^2$ of a generic $1$-parameter Hopf bifurcation in the case $\dd$. The corresponding branching periodic solutions cannot be observed in the case \od. All remaining generic $1$-parameter Hopf bifurcations in the case \dd are as described in (\romannumeral 1).
	\end{enumerate}
\end{theorem}
\begin{remk}
	Note that in both situations the indecomposable component of $\No$ can be high-dimensional due to symmetry. Hence, in general \Cref{thm:ddhopf} does not describe a standard Hopf bifurcation with $2$-dimensional center manifold.
	\hspace*{\fill}$\triangle$
\end{remk}
\begin{cor}
	\label{cor:ddhopf}
	Assume that the specific network structure forces the fundamental network to decompose into only components of real type that are pairwise non-isomorphic. Then Hopf bifurcations in generic $1$-parameter families are only possible in networks with internal dynamics of dimension greater or equal to $2$. 
\end{cor}
\begin{remk}
	In particular, both conditions of \Cref{cor:ddhopf} hold true if the network structure forces the linear admissible maps to only have real eigenvalues. An example of this phenomenon is the class of feedforward networks that were introduced in \textcite{vonderGracht.2022}. We investigate $1$-parameter steady state bifurcations in feedforward networks with high-dimensional internal dynamics in \Cref{sec:ff}.
	\hspace*{\fill}$\triangle$
\end{remk}

\subsection[Generic \texorpdfstring{$l$}{l}-parameter bifurcations]{Generic $l$-parameter bifurcations}
\label{subsec:il}
The situation for $l$-parameter bifurcations is a lot more involved than in \Cref{subsec:is,subsec:ih}. Similar precise statements relating the case \dd to \od are not possible in full generality, as conditions \eqref{eq:condker} and \eqref{eq:condcent} allow for greater flexibility in the composition of generalized kernel or the center subspace the larger the value of $l$ is. Nevertheless, the underlying mechanism that made the two different characterizations in \Cref{thm:ddhopf} possible, applies to $l$-parameter bifurcations (with $l >1$) as well. In the case \dd the representation space $\Nd$ decomposes into the indecomposable subrepresentations of $\No$, each occurring $d$ times. These subrepresentations can be components of generalized kernels or center subspaces. Hence, there are potentially numerous possibilities to find suitable combinations of components that satisfy \eqref{eq:condker} and \eqref{eq:condcent}. Nonetheless, any combination of subrepresentations that occurs in a generic $l$-parameter bifurcation in the case \od\ -- i.e., one that does not make use of extra copies -- also occurs as a generalized kernel or center subspace in a generic $l$-parameter bifurcation in the case \dd. Once again, the reduced bifurcation problems are equivalent due to symmetry. They can be seen as the ones that are inherent to the network structure and independent of the internal dynamics. On the other hand, in general a generalized kernel or center subspace in a generic $l$-parameter bifurcation in the case \dd with $d \ge 2$ contains multiple copies of the same indecomposable component of $\No$. Then there is no equivalent reduced bifurcation problem in the case \od. We summarize these results as
\begin{theorem}
	\label{thm:ddlparam}
	Generic $l$-parameter bifurcations in a fundamental network with $1$\nobreakdash-dimensional internal dynamics are also generic in the same network with $d$-dimensional internal dynamics.
\end{theorem}

More generally, \Cref{thm:ddlparam} follows almost directly from \Cref{cor:hddecomp}. The total phase space $\No$ is a subrepresentation of $\Nd$. Hence, any combination of indecomposable components of $\No$ that make up a generalized kernel in a generic $l$\nobreakdash-parameter bifurcation in the case \od also occur in a generic $l$\nobreakdash-parameter bifurcation problem in the case \dd. This yields the previous result. Even more so, it can be generalized to compare bifurcations in the same network with internal dynamics of arbitrary dimension. If $\dim{W^1} \le \dim{W^2}$ the total phase space $\bigoplus_{\sigma\in\Sigma}W^1$ is isomorphic to a subrepresentation of $\bigoplus_{\sigma\in\Sigma}W^2$. Hence, every generalized kernel or center subspace in a generic $l$\nobreakdash-parameter bifurcation in $\bigoplus_{\sigma\in\Sigma}W^1$ also occurs as a generalized kernel or center subspace in a generic $l$\nobreakdash-parameter bifurcation in $\bigoplus_{\sigma\in\Sigma}W^2$.
\begin{theorem}
	\label{thm:ddlparamgen}
	Generic $l$-parameter bifurcations in a fundamental network with $d_1$\nobreakdash-dimensional internal dynamics are also generic in the same network with $d_2$-dimensional internal dynamics, whenever $d_1 \le d_2$.
\end{theorem}

On the other hand, for a fixed number of parameters $l$, conditions \eqref{eq:condker} and \eqref{eq:condcent} impose restrictions on the maximal number of indecomposable components of $\No$ that can occur as a generalized kernel or as a center subspace in a generic $l$-parameter bifurcation for any value of $d$. More precisely the generalized kernel can at most be composed of $l$ components of real type, of $\lfloor l/2 \rfloor$ components of complex type, or of $\lfloor l/4 \rfloor$ components of quaternionic type. Here $\lfloor a \rfloor$ denotes the nearest smaller or equal integer. Likewise, the center subspace can at most be composed of $2l$ components of real type, of $l$ components of complex type, or of $l$ components of quaternionic type. In particular, the total number of indecomposable components is always less than or equal to $l$ for generalized kernels and less than or equal to $2l$ for center subspaces. Recall that increasing the dimension of the internal dynamics $d$ yields additional copies of the indecomposable components of $\No$ in the decomposition of $\Nd$. In particular, we find all possible combinations of $l$ or $2l$ indecomposable components in the case that $d = l$ or $d=2l$, respectively. Increasing $d$ beyond these critical values does not provide any further solutions to the combinatorial problems \eqref{eq:condker} and \eqref{eq:condcent}. As a result, all possible generalized kernels in a generic $l$-parameter bifurcation problem with internal dynamics of dimension $d' > l$ can also be observed in the case $d=l$. Likewise, all possible center subspaces in a generic $l$-parameter bifurcation problem with internal dynamics of dimension $d' > 2l$ can also be observed in the case $d=2l$. Once again, the reduced bifurcation problems are therefore equivalent to those in the cases $d=l$ and $d=2l$, respectively. In combination with \Cref{thm:ddlparamgen}, we obtain
\begin{theorem}
	\label{thm:ddlparamall}
	\begin{enumerate}[label=(\roman*)]
		\item All generic $l$-parameter steady state bifurcations that can occur in a fundamental network can be investigated in the case of an internal phase space of dimension $d=l$.
		\item All generic $l$-parameter Hopf bifurcations that can occur in a fundamental network can be investigated in the case of an internal phase space of dimension $d=2l$.
	\end{enumerate}
\end{theorem}
\begin{remk}
	The minimal values of the dimension of the internal phase space $d$ stated in \Cref{thm:ddlparamall} are optimal in the sense that for any fundamental network there is indeed a bifurcation scenario that can only occur in a generic $l$-parameter bifurcation problem if $d \ge l$ or $d \ge 2l$, respectively. In the remainder of this remark, we construct the corresponding decomposition of the center subspace into indecomposable subrepresentations of $\No$ which specifies this bifurcation scenario.
    
    To that end, let
	\[ \Delta_0 = \left\{ \left. (x_\sigma)_{\sigma\in\Sigma} \in \No\ \right|\ x_\sigma = x_\tau, \quad \text{for all} \quad \sigma,\tau \in \Sigma \right\} \subset \No \]
	be the fully synchronous subspace in the case \od. In particular, $\Delta_0$ is a (not necessarily complementable) subrepresentation on which each representation map $\Ao_\sigma$ acts as the identity. This assertion holds true independent of the specific fundamental network at hand. Assume $Y \subset \No$ is a subrepresentation with $Y \cong \Delta_0$. Then all $\Ao_\sigma$ act as the identity on $Y$ as well. That is, for all $y = (y_\sigma)_{\sigma\in\Sigma} \in Y$ we have $(\Ao_\tau y)_\sigma = y_{\sigma\tau} = y_\sigma$ for all $\sigma,\tau\in\Sigma$ by definition of the right regular representation. In particular, for $\sigma=\Id$ we see
	\[ y_{\Id} = y_\tau \]
	for all $\tau\in\Sigma$. As $\Delta_0$ is one-dimensional, we obtain $y \in \Delta_0$ and therefore $Y = \Delta_0$. In particular, there is no subrepresentation in $\No$ that is isomorphic but not equal to $\Delta_0$. This observation depends crucially on the fact that the internal dynamics are one-dimensional. In general, there is an indecomposable component of $\No$ that contains the fully synchronous subspace $Y \supset \Delta_0$. This follows from the fact that $\Delta_0$ is spanned by $(1, \dotsc, 1)^T$, which is an eigenvector of all the equivariant projections onto the indecomposable components, as these are network maps. If a second subrepresentation $Y' \subset \No$ is isomorphic to $Y$ the consideration above implies that $Y$ and $Y'$ contain $\Delta_0$. Hence, $Y \cap Y' \supset\Delta_0 \ne \emptyset$.
	
	As a result, any decomposition of $\No$ into indecomposable components contains precisely one component isomorphic to $Y$ where $\Delta_0 \subset Y$. Consequently, in the case of $d$-dimensional internal dynamics there is a subrepresentation $U \subset \Nd$ with $U \cong Y^{l}$ only if $d \ge l$ by \Cref{thm:dechd}. Likewise, a subrepresentation $U \subset \Nd$ with $U \cong Y^{2l}$ exists only if $d \ge 2l$.
	\hspace*{\fill}$\triangle$
\end{remk}
\begin{remk}
    The argument in the previous remark uses only the property of the subrepresentation $\Delta_0$ that it occurs precisely once in any fundamental network. It can thus be applied in the same way to any other subrepresentation with that property in any other fundamental network. On the other hand, if a subrepresentation $Y\subset\No$ occurs with multiplicity $m$, a subrepresentation $U\subset\Nd$ with $U\cong Y^l$ (or $U\cong Y^{2l}$) exists only if the internal dimension $d$ satisfies $dm\ge l$ (or $dm\ge 2l$, respectively).
	\hspace*{\fill}$\triangle$
\end{remk}
\begin{remk}
	In the case of $1$-parameter bifurcations \Cref{thm:ddlparamall} shows that all generic steady state bifurcations can be observed in the network with $1$-dimensional internal dynamics and all Hopf bifurcations can be observed in the case of $2$-dimensional internal dynamics. This matches the results in \Cref{subsec:is,subsec:ih}.
	\hspace*{\fill}$\triangle$
\end{remk}

\begin{figure}[h]
	\begin{center}
		\resizebox{.55\linewidth}{!}{	\centering
        \begin{tikzpicture}[->,
        	>=stealth',
        	shorten >=1pt,
        	auto,
        	node distance=1cm,
        	main node/.style={line width=1.5pt, circle, scale = 3, draw, font=\sffamily\tiny, inner sep=1pt}
        ]
        	\node[main node] (1) at (-4,0) {$1$};
        	\node[main node] (2) at (-2,0) {$2$};
        	\node[main node] (3) at (0,0) {$3$};
            \node[main node] (4) at (1.7321,1) {$4$};
            \node[main node] (5) at (1.7321,-1) {$5$};
        
        	\path[every node/.style={font=\sffamily\small}, line width =1.5pt]
            	(2) edge [color = {red}] node {} (1)
            	(3) edge [color = {red}, bend left = 0] node {} (2)
                (4) edge [color = {red}, bend left = -15] node {} (3)
                (5) edge [color = {red}, bend left = -15] node {} (4)
                (3) edge [color = {red}, bend left = -15] node {} (5)
            ;
        \end{tikzpicture}
		}%
	\end{center}%
	\caption{A $5$-cell fundamental network.}
	\label{fig:tadpole}
\end{figure}%

\begin{ex}
    Consider the class of network ODEs given by
    \begin{align}\label{eq:ex5celll}
		\dot{X}	= 
        \begin{pmatrix}
        f(X_1, X_2, X_3, X_4, X_5) \\
		f(X_2, X_3, X_4, X_5, X_3) \\
        f(X_3, X_4, X_5, X_3, X_4) \\
		f(X_4, X_5, X_3, X_4, X_5) \\
        f(X_5, X_3, X_4, X_5, X_3) 
        \end{pmatrix}.
	\end{align}
    \Cref{fig:tadpole} shows the corresponding network structure, which is isomorphic to its own fundamental network. 
    For convenience, we have only drawn one type of interaction, with the others found by tracing the given arrows back two, three, and four times. Self-loops corresponding to internal dynamics have likewise been left out. 
    Systems of the form $\eqref{eq:ex5celll}$ are precisely those for which the linear map
    \begin{align*}
		A	= 
        \begin{pmatrix}
        0 & 1 & 0 & 0 & 0 \\
        0 & 0 & 1 & 0 & 0 \\
        0 & 0 & 0 & 1 & 0 \\
        0 & 0 & 0 & 0 & 1 \\
        0 & 0 & 1 & 0 & 0 
        \end{pmatrix}
	\end{align*}
    sends solutions to solutions. This map satisfies $A^2 = A^5$ and generates, together with the identity map, a monoid of $5$ elements.

    It is shown in \cite{Nijholt.2017d} that a decomposition into indecomposable representations for the \od case  is given by
    \begin{align*}
        \No=\mathbb{R}^5 = W_1 \oplus W_2 \oplus U\, ,
    \end{align*}
    where
    \begin{align*}
        W_1 &= \{(X,X,X,X,X)\mid X \in \mathbb{R}    \} \\ 
        W_2 &= \{(X,Y,0,0,0)\mid X,Y \in \mathbb{R}    \} \\
        U &= \{(X,Y,Z,X,Y)\mid X,Y,Z \in \mathbb{R}\, , X+Y+Z=0    \}\, .
    \end{align*}
    These are all non-isomorphic, with $W_1$ and $W_2$ of real type and $U$ of complex type. 
    Thus, for the case \dd we have $d$ copies of each of these indecomposables.
    
    Since $U$ is of complex type, an eigenvalue movement through the origin is a co-dimension $2$ event for this subrepresentation.
    Thus, in the case of a generic $1$-parameter steady state bifurcation, the generalized kernel may either be a single copy of $W_1$ or of $W_2$. 
    Similarly in a $2$-parameter steady state bifurcation, the generalized kernel may be isomorphic to either $U$, $W_i$ or $W_i\oplus W_j$, for any choice of $i,j \in \{1,2\}$. 
    It therefore suffices to set $d=2$ in order to observe all $1$- and $2$-parameter steady state bifurcations that the network structure in \Cref{fig:tadpole} supports generically, which illustrates the mechanism and the result of \Cref{thm:ddlparamall}.

    The center subspace in a generic $1$-parameter bifurcation (i.e., potentially also caused by complex eigenvalues crossing the imaginary axis) can be given by 
    \[ U\, , \,\,\, W_1 \oplus W_1 \, \text{ or } \, W_2 \oplus W_2\, ,\] 
    as well as the aforementioned generalized kernels $W_1$ and $W_2$.
    In a generic $2$-parameter bifurcation, the center subspace will be isomorphic to one of 
    \begin{alignat*}{3}
        &W_1\oplus W_1 \oplus W_2 \oplus W_2\, ;  &\qquad&W_1\oplus W_1 \oplus U\, ; &\qquad& \,W_2\oplus W_2 \oplus U\, ; \\ \nonumber
        & U\oplus U\, ; &&\bigoplus^4 W_1 \, ; &&\bigoplus^4 W_2\, 
    \end{alignat*}
    or any complementable subrepresentation thereof.
    Again in accordance with \Cref{thm:ddlparamall}, the last two show that we need $d\geq 4$ to observe all generic local $1$- and $2$-parameter (steady state or Hopf) bifurcations.\hspace*{\fill}$\triangle$
\end{ex}

\subsection{Beyond qualitative statements using center manifold reduction}
\label{subsec:in}
\Cref{subsec:is,subsec:ih,subsec:il} describe how to determine qualitative bifurcations in homogeneous coupled cell systems with (possibly) high-dimensional internal dynamics. In particular, how (parts of) the branching pattern in the case \dd can be observed in the case \od. The reason why the restriction to qualitative statements needs to be made is the fact that the relation between the two cases is made in terms of reduced bifurcation problems. The reduction methods (namely Lyapunov-Schmidt and center manifold reduction) require coordinate changes so that whatever information is stored in the precise choice of coordinates is lost when applying \Cref{thm:ddss,thm:ddhopf,thm:ddlparam}. Most importantly, in the investigation of network dynamics this includes the possibility to distinguish individual cells from the coordinates of the total phase space variables as they were chosen in \Cref{sec:pr}. However, \Cref{thm:dechd} allows for more precise statements about the case \dd, if more knowledge about bifurcations in the case \od is available. We can explicitly construct bifurcating branches in the case \dd from those in the case \od, capturing the spirit of \Cref{thm:ddlparam}, without losing all information about each cell's behavior.

\subsubsection{Technicalities on the center manifold reduction}
\label{subsubsec:cmfr}
Let us briefly recapture the center manifold reduction and its usage for bifurcation analysis of systems governed by fundamental network vector fields \eqref{eq:netvfparam} as it was introduced in \textcite{Nijholt.2017c,Nijholt.2019}. 
The main goal of this section is to introduce some notation. 
We note that the remainder of this paper can be understood without reading the proofs of \Cref{lem:1,lem:2,lem:Qprime,lem:2.5} below. We also add some additional technicalities that were not explicitly proved in \textcite{Nijholt.2017c,Nijholt.2019}. 
As the techniques presented here can be applied to steady state bifurcations as well as to Hopf bifurcations, 
we slightly abuse notation from now on and denote the generalized kernel or the center subspace and the respective complement by $\XX^c$ and $\XX^h$. We obtain equivariant projections along $\XX^c$ and $\XX^h$, respectively,
\[ P^h \colon \bigoplus_{\sigma\in\Sigma}V \to \XX^h, \quad P^c \colon \bigoplus_{\sigma\in\Sigma}V \to \XX^c. \]
These allow to uniquely split any element $v \in \bigoplus_{\sigma\in\Sigma}V$ as \mbox{$v=v^c+v^h$} with $v^c = P^c(v)$ and \mbox{$v^h=P^h(v)$}. 

In order to apply the center manifold reduction to bifurcation problems, one extends the system \eqref{eq:netvfparam} to $\uu{v} = (v, \lambda) \in \bigoplus_{\sigma\in\Sigma} V \times \RR^l$ to include the parameter as a dynamic variable:
\begin{equation}
	\label{eq:extended}
	\dot{\vu} = \begin{pmatrix} \dot{v} \\ \dot{\lambda} \end{pmatrix} = \begin{pmatrix} \Gamma_f (v, \lambda) \\ 0 \end{pmatrix} = \Gu_f(\vu).
\end{equation}
Solutions to this system are in one-to-one correspondence with solutions to the original system. In particular, $\uu{\Gamma}_f(0,0)=0$. If we extend the representation matrices accordingly
\[ \uu{A}_\sigma \colon \uu{v} = (v,\lambda) \mapsto (A_\sigma v, \lambda), \]
this system, furthermore, remains equivariant $\uu{\Gamma}_f \circ \uu{A}_\sigma = \uu{A}_\sigma \circ \uu{\Gamma}_f$. The following lemma explores the consequences of equivariance with respect to the extended monoid representation.
\begin{lemma}
	\label{lem:1}
	Let
	\begin{align*}
		F \colon \bigoplus_{\sigma\in\Sigma} V \times \RR^l	&\to \bigoplus_{\sigma\in\Sigma} V \times \RR^l \\
		(v, \lambda)										&\mapsto \begin{pmatrix} F_v(v, \lambda) \\ F_\lambda (x,\lambda) \end{pmatrix}
	\end{align*}
	be equivariant with respect to the extended monoid representation: $F \circ \uu{A}_\sigma = \uu{A}_\sigma \circ F$ for all $\sigma \in \Sigma$. Then the component functions $F_v$ and $F_\lambda$ are parameter dependent and equivariant or invariant with respect to the non-extended monoid representation, respectively. That is
	\[ F_v (A_\sigma v, \lambda) = A_\sigma F_v(v,\lambda), \qquad F_\lambda(A_\sigma v, \lambda) = F_\lambda(v,\lambda). \]
\end{lemma}
\begin{proof}
	This can be seen directly from the equivariance condition
	\begin{equation*} 
		\left(F\circ \uu{A}_\sigma \right) (v, \lambda) = \begin{pmatrix} F_v (A_\sigma v, \lambda) \\ F_\lambda (A_\sigma v, \lambda) \end{pmatrix} = \begin{pmatrix} A_\sigma F_v (v, \lambda) \\ F_\lambda (v,\lambda) \end{pmatrix} = \left(\uu{A}_\sigma \circ F \right) (v,\lambda).
	\end{equation*}
\end{proof}

Under the bifurcation assumption -- eigenvalue $0$ or purely imaginary eigenvalues -- the linearization of the extended system \eqref{eq:extended} induces a splitting
\[ \bigoplus_{\sigma\in\Sigma}V \times \RR^l = \Xu^c \oplus \Xu^h \]
into center subspace -- generalized eigenspace to purely imaginary eigenvalues -- and its hyperbolic complement of $D\uu{\Gamma}_f(0,0)$. These are subrepresentations of the extended monoid representation $\bigoplus_{\sigma\in\Sigma}V \times \RR^l$. Furthermore, they are also related to the corresponding subspaces of the non-extended system. In particular,
\begin{equation}
	\label{eq:cspextended}
	\begin{split}
		\uu{\XX}^h	&= \XX^h \times \{ 0 \}, \\
		\uu{\XX}^c	&= \XX^c \times \{ 0 \} \oplus \left\langle (v_1^h, \lambda_1), \dotsc, (v_l^h, \lambda_l) \right\rangle,
	\end{split}
\end{equation}
where $\lbrace \lambda_1, \dotsc, \lambda_l \rbrace$ is a basis of $\RR^l$ and $v_1^h, \dotsc, v_l^h \in \XX^h$ are chosen suitably. The extended projections along $\uu{\XX}^c$ and along $\uu{\XX}^h$, respectively, are denoted accordingly
\[ \uu{P}^h \colon \bigoplus_{\sigma\in\Sigma}V\times\RR^l \to \uu{\XX}^h, \quad \uu{P}^c \colon \bigoplus_{\sigma\in\Sigma}V\times\RR^l \to \uu{\XX}^c. \]
They allow us to uniquely split every $\uu{v} \in \bigoplus_{\sigma\in\Sigma}V \times \RR^l$ as $\uu{v} = \uu{v}^c+\uu{v}^h$ with $\uu{v}^c = \uu{P}^c(\uu{v})$ and $\uu{v}^h=\uu{P}^h(\uu{v})$.
\begin{lemma}[see also the proof of Lemma 5.3 in \cite{Nijholt.2019}]
	\label{lem:2}
	Let 
	\begin{equation}
		\label{eq:uuPc}
		\uu{P}^c \colon \bigoplus_{\sigma\in\Sigma} V \times \RR^l \to \uu{\XX}^c
	\end{equation}
	be the equivariant projection onto the center subspace of the extended system. Then $\uu{P}^c$ acts as the identity on the $\lambda$-component, i.e., $\uu{P}^c(v,\lambda) = \left(\uu{P}^c_v(v, \lambda), \lambda \right)$.
\end{lemma}
\begin{proof}
	Fix $(v,\lambda) \in \bigoplus_{\sigma\in\Sigma} V \times \RR^l$. Recall that
	\[ \bigoplus_{\sigma\in\Sigma} V = \XX^c \oplus \XX^h. \]
	Hence we find a unique representation $v = v^c + v^h$ with $v^c \in \XX^c$ and $v^h \in \XX^h$ using projections $v^c = P^c (v), v^h = P^h(v)$. From the representation in \eqref{eq:cspextended}, we see that $(v^c, 0) \in \uu{\XX}^c$. Furthermore, we may represent the parameter as $\lambda = \alpha_1 \lambda_1 + \dotsb + \alpha_l \lambda_l$ and it holds that
	\[ \alpha_1 (v_1^h, \lambda_1) + \dotsb + \alpha_l (v_l^h, \lambda_l) \in \uu{\XX}^c. \]
	Lastly, as $v_1^h, \dotsc, v_l^h \in \XX^h$ we obtain
	\[ \left(v^h - (\alpha_1 v_1^h + \dotsb + \alpha_l v_l^h), 0 \right) \in \uu{\XX}^h. \]
	Since $\uu{P}^c$ is the projection onto $\uu{\XX}^c$, it acts as the identity on $\uu{\XX}^c$ and maps $\uu{\XX}^h$ to $0$. We obtain
	\begin{align*}
	(v,\lambda) = (v^c+v^h,\lambda) &= (v^c,0) + (v^h - (\alpha_1 v_1^h + \dotsb + \alpha_l v_l^h), 0) + \alpha_1 (v_1^h, \lambda_1) + \dotsb + \alpha_l (v_l^h, \lambda_l) \\
									&\mapsto (v^c,0) + 0 + \alpha_1 (v_1^h, \lambda_1) + \dotsb + \alpha_l (v_l^h, \lambda_l) \\
									&= \left( v^c + \alpha_1 v_1^h + \dotsb + \alpha_l v_l^h, \lambda \right).
	\end{align*} 
\end{proof}

We define another projection that `forgets' about the skewed directions of the extended center subspace. Every element in $\RR^l$ can uniquely be represented in terms of its basis as \mbox{$\lambda = \alpha_1 \lambda_1 + \dotsb + \alpha_l \lambda_l$} with coefficients $\alpha_1, \dotsc, \alpha_l \in \RR$. Thus, every element in $\uu{\XX}^c$ has a unique representation as
\[ \uu{v}^c = (v^c,0) + \alpha_1 (v_1^h, \lambda_1) + \dotsb + \alpha_l (v_l^h, \lambda_l) \]
where $v^c = P^c(v) \in \XX^c$. Hence, the map
\begin{equation}
	\label{eq:Pprime}
	\begin{split}
		P' \colon \uu{\XX}^c	&\to \XX^c \times \RR^l \\
		\uu{v}^c 				&\mapsto (v^c, \alpha_1 \lambda_1 + \dotsb + \alpha_l \lambda_l) = (v^c,\lambda)
	\end{split}
\end{equation}
is a projection. Note that this map satisfies $P' = \left. (P^c \times \mathbbm{1}_{\RR^l}) \right|_{\uu{\XX}^c}$.
\begin{lemma}[see also the proof of Theorem 5.4 in \cite{Nijholt.2019}]
	\label{lem:Qprime}
	The projection $P' \colon \uu{\XX}^c \to \XX^c \times \RR^l$ is invertible with inverse given by
	\begin{equation}
		\label{eq:Qprime}
		\begin{split}
			Q' \colon \XX^c \times \RR^l	&\to \uu{\XX}^c \\
			(v^c,\lambda)					&\mapsto (v^c,0) + \alpha_1 (v_1^h, \lambda_1) + \dotsb + \alpha_l (v_l^h, \lambda_l) = \uu{v}^c,
		\end{split} 
	\end{equation}
	where $\lambda = \alpha_1 \lambda_1 + \dotsb + \alpha_l \lambda_l$ is a representation in the basis $\RR^l = \langle \lambda_1, \dotsc, \lambda_l \rangle$ and $v_1^h, \dotsc,v_l^h \in \XX^h$ suitable as in \eqref{eq:cspextended}.
\end{lemma}
\begin{proof}
	This follows directly from the definition of $P'$ in \eqref{eq:Pprime}.
\end{proof}

After restriction to a suitably small neighborhood around the bifurcation point $(0,0)$ the extended system \eqref{eq:extended} admits a unique center manifold $M^c$, invariant under the dynamics, which locally contains all solutions whose $\Xu^h$-component is bounded. In the formalism of the extended system this includes all bifurcating steady state or periodic solutions. The center manifold can be realized as the graph of a function $\pu \colon \Xu^c \to \Xu^h$. More precisely
\begin{equation}
	\label{eq:cmfr}
	M^c = \left\{ Q'(v^c,\lambda) + \uu{\psi}(Q'(v^c, \lambda)) \mid (v^c,\lambda) \in \XX^c \times \Lambda \right\} \subset \bigoplus_{\sigma\in\Sigma}V \times \RR^l
\end{equation}
In particular, the dynamics of a generic system of the form \eqref{eq:extended} restricted to the center manifold is bijectively conjugate to that of a generic system of the form 
\begin{equation}
	\label{eq:cmfred}
	\begin{split}
		\dot{v}^c		&= r(v^c,\lambda), \\
		\dot{\lambda}	&= 0
	\end{split}
\end{equation}
on $\XX^c \times \RR^l$, where
\begin{enumerate}[label=(\roman*)]
	\item $r(0,0)=0$.
	\item The center subspace of $D_{v^c}r(0,0)$ is the full space $\XX^c$.
	\item It is equivariant with respect to the (non-extended) monoid representation restricted to $\XX^c$.
\end{enumerate}
The conjugation is realized by the maps $P = P' \circ \uu{P}^c = P^c \times \mathbbm{1}_{\RR^l}$ and $\uu{\psi} \circ Q'$. As all maps that are needed to define the center manifold -- in particular $\uu{\psi}, P$ and the projections onto $\uu{\XX}^c$ and $\uu{\XX}^h$ -- are equivariant as well, the entire center manifold $M^c$ is invariant under the extended monoid representation.

Lastly, we state a technical result that explicitly describes the parameter-dependence of the center manifold $M^c$ of the extended system and the maps that define the interconnection between the center manifold and the center subspace.
\begin{lemma}
	\label{lem:2.5}
	Let $\uu{\psi} \colon \uu{\XX}^c \to \uu{\XX}^h$ be the map whose graph is the center manifold $M^c$ as in \eqref{eq:cmfr}. Then $\uu{\psi}$ has a trivial $\lambda$-component, i.e., $\uu{\psi}(v,\lambda) = \left( \uu{\psi}_v(v,\lambda), 0 \right)$.
\end{lemma}
\begin{proof}
	This follows directly from the representation of the subspaces in \eqref{eq:cspextended}: $\uu{\XX}^h = \XX^h \times \{ 0 \}$.
\end{proof}
\begin{remk}
	\Cref{lem:2,lem:2.5} can be summarized by stating that the center manifold and the center subspace of the extended system share the same parameter component. Furthermore, the same holds for the reduced system on $\XX^c \times \RR^l$ governed by $r$. The conjugation maps leave the $\lambda$-component unchanged. Introducing the parameter dependence into the system for bifurcation analysis yields a skewed center subspace (see \eqref{eq:cspextended}) that extends into parameter space. The same holds for the center manifold. This can be regarded as a continuation of the center manifold of the non-extended system for $\lambda=0$ to varying parameter values. Nevertheless, no dynamic effects occur in the $\lambda$-component of the extended system.
	\hspace*{\fill}$\triangle$
\end{remk}

\subsubsection{\texorpdfstring{Determining $\dd$-branches from \od-branches}{Determining DD-branches from 1D-branches}}
\label{subsubsec:ddbranch}
Assume now that we know the branching behavior of each cell in a branch of a generic bifurcation of steady states of the system \eqref{eq:netvfparam} in the case \od. That is, we have a smooth curve of steady states or of initial conditions for periodic solutions $(x_\sigma (\lambda))_{\sigma\in\Sigma}$ for small absolute values of $\lambda$. We denote the subrepresentation that forms the corresponding generalized kernel or center subspace by $\Xo^c \subset \No$. Due to the conjugation of the dynamics on the center manifold and the reduced system \eqref{eq:cmfred}, we obtain that
\[ \Po^c((x_\sigma (\lambda))_{\sigma\in\Sigma}) \]
-- where $\Po^c$ is the projection onto the generalized kernel or center subspace $\Xo^c$ in the case \od\ -- is the branching solution of a generic reduced bifurcation problem corresponding to $(x_\sigma (\lambda))_{\sigma\in\Sigma}$. As the center subspace, in general, is a proper subspace, not every cell's coordinate entry $x_\sigma(\lambda)$ can be found in the projected solution branch. For example, the projection might map coordinates to zero or sum up multiple coordinate entries. Nevertheless, as long as the projection $\Po^c$ is known, this method provides a method to represent the qualitative branching solutions of the reduced system while respecting the coordinates that are chosen according to the cells of the network. This projection operator, however, is often computed while classifying the generic bifurcations in the case \od as a byproduct.

Furthermore, \Cref{thm:ddlparam} shows that any bifurcating branch in the case \od occurs generically in the same network in the case \dd as well. \Cref{thm:dechd} allows us to transform a generic \mbox{\od-branch} into a generic \dd-branch. In particular, there is a generic $l$-parameter bifurcation problem in $\Nd$ whose generalized kernel or center subspace is isomorphic to $\Xo^c$ in the case $\od$. Furthermore, from \eqref{eq:1DDDa} we see that this generalized kernel or center subspace $\Xd^c \subset \Nd$ is of the form
\begin{equation}
	\label{eq:genkerDD}
	\Xd^c \cong \Xo^c \ox \langle w \rangle,
\end{equation}
where $w \in \Vd \setminus \{0\}$. Denote the equivariant isomorphism in \eqref{eq:genkerDD} by $\Psi \colon \Xo^c \ox \langle w \rangle \to \Xd^c$. Consequently, the branch $(x_\sigma (\lambda))_{\sigma\in\Sigma}$ can be represented as a generic branch on the center subspace $\Xd^c$ as 
\[ \Psi(\Po^c (x_\sigma (\lambda))_{\sigma \in \Sigma} \ox w). \]
In more general terms, there exists a direction $w \in \Vd \setminus \{0\}$ such that the generic bifurcation pattern -- that is, all generic bifurcating solutions -- restricted to the center subspace in the $1$-dimensional case is reflected in the $d$-dimensional case, where internal dynamics is restricted to this direction $w$. This interpretation, however, is only fully accurate in the case that $\Psi$ is the identity -- after identifying tensor and non-tensor notation. In general, it yields only qualitatively the same bifurcation diagram in the case \dd. Nevertheless, as the coordinates for the center subspace in the case \dd reflect the cells of the network, we can read off cell-by-cell information from this representation. Finally, using the conjugacy between the center manifold and the reduced system once more, this time in the case $\dd$, we find the representation of the branch $(x_\sigma (\lambda))_{\sigma\in\Sigma}$ in the center manifold of the \dd-system as
\[ \uu{\psi}^{\mathbf{\mathsf{D}}} \left( Q'_{\mathbf{\mathsf{D}}} \left( \Psi(\Po^c (x_\sigma (\lambda))_{\sigma \in \Sigma} \ox w), \lambda \right) \right). \]
Due to \Cref{thm:ddlparam} this branch occurs in a generic bifurcation problem with generalized kernel or center subspace isomorphic to $\Xo^c \ox \langle w \rangle$.

In theory, this procedure provides a mechanism to translate generic branching solutions in the case \od to generic solutions in the case \dd. As long as information about the maps $\Po^c, Q'_{\mathbf{\mathsf{D}}}, \uu{\psi}^{\mathbf{\mathsf{D}}}$ and $\Psi$ is available, it also transforms information about the branching behavior of each individual cell. However, this latter part is not to be expected in general. The results in this chapter aim at simplifying the investigations of generic steady state bifurcations in the case \dd by restricting to the case \od. In particular, we do not want to determine a generic center manifold in the high-dimensional case. As a result, knowledge about these maps is not available in general. Nevertheless, the observations in this subsection provide a theoretical tool to relate generic branching solutions in the case \dd to those in the case \od. In \Cref{sec:ff}, we investigate how additional structure in the network -- in particular feedforward structure -- provides information about the maps $\Po^c, Q'_{\mathbf{\mathsf{D}}}, \uu{\psi}^{\mathbf{\mathsf{D}}}$ and $\Psi$ without explicitly computing them. This suffices to exploit this mechanism to characterize generic branching solutions in the case \dd from those in the case \od including the behavior of each cell without computing the center manifolds.
\begin{remk}
	The method presented in this section can also be used in the spirit of \Cref{thm:ddlparamall} to translate a bifurcating branch of steady state or periodic solutions in a fundamental network with $d_1$-dimensional internal dynamics into one for the same network with $d_2$-dimensional internal dynamics if $d_1 \le d_2$ and whenever the corresponding generalized kernel or center subspace occurs in both cases generically. In particular, when a branch exhibits a robust pattern of synchrony -- i.e., coordinate entries corresponding to a balanced partition of the cells coincide -- then the representation of that branch exhibits the same synchrony for internal dynamics of any dimensions for which it exists generically. This follows from \Cref{prop:ddodsynch}.
	\hspace*{\fill}$\triangle$
\end{remk}

\section{Steady state bifurcations in feedforward networks}
\label{sec:ff}
In this section, we investigate how additional structure in the network helps to gain more insight into the generic steady state bifurcations of networks with high-dimensional internal dynamics for the example of feedforward structure. In \Cref{sec:hd,sec:bi}, we show that symmetry of the fundamental network forces the qualitative bifurcation picture in the high-dimensional case to be the same as that in the case \od. Exploiting the feedforward structure, we can even understand bifurcations in individual cells, which was already hinted at in \Cref{subsec:in}. We begin by recalling different equivalent characterizations of feedforward networks as they are introduced in \textcite{vonderGracht.2022}. The generic steady state bifurcation result in the case \od is summarized in \Cref{subsec:ffod}. Finally, we discuss implications for generic steady state bifurcations in the case \dd. Note that the results in \cite{vonderGracht.2022} hold for general homogeneous networks with asymmetric inputs. Here, we restate only the required bits in the context of fundamental networks.

In general terms, the feedforward structure can be interpreted as the absence of any feedback. In that sense, a network is called a feedforward network if it does not contain any closed, directed loops of length $2$ or larger. Note that this allows for self-loops. We do not exclude the situation that a cell influences itself, which is not considered as feedback in this setting.

Secondly, we view feedforward structure as a notion of regularity in the network: all arrows are either self-loops or `oriented in the same direction'. We introduce a preorder on the monoid of input maps $\Sigma$ which are the cells of the fundamental network by saying $\sigma \nole \sigma'$ if there is a path from $\sigma'$ to $\sigma$. For a fundamental network, this means $\sigma \nole \sigma'$ if there exists $\tau \in \Sigma$ with $\tau\sigma = \sigma'$. Then, a network does not contain any closed directed loops of length $2$ or larger, if and only if this preorder is a partial order which reflects the aforementioned regularity. We denote the situation that $\sigma \nole \sigma'$ but $\sigma \ne \sigma'$ by $\sigma \nol \sigma'$. 

In \Cref{subsec:ffod}, we see that we may relate the asymptotics in a generic feedforward bifurcation to the partial order of the cells in the case \od. For a fundamental network, we want to investigate the implications of this for the case \dd.
\begin{remk}
    It can be shown, that a homogeneous coupled cell network is a feedforward network if and only if its fundamental network is. In that sense, the upcoming analysis is also relevant for networks that are not fundamental networks.
\end{remk}
\begin{remk}
	\label{rem:pot}
	Note that labeling the nodes $\{\sigma_1, \dotsc, \sigma_N\}$ of a \emph{feedforward} fundamental network according to $\nole$, such that $\sigma_i \nole \sigma_j$ implies $i \le j$, yields that all the linear admissible maps are upper triangular with entries in $\gl{V}$ -- or, put differently, upper triangular block matrices. If $V = \RR$ they are truly upper triangular. In the tensor notation (see \Cref{prop:linadmtensor}) this yields that linear admissible maps are of the form
	\[ \sum_{\sigma\in\Sigma} B_\sigma \ox b_\sigma, \]
	where $b_\sigma \in \gl{V}$ and the $B_\sigma$ are upper triangular. In particular, this implies that all the equivariant linear maps of the regular representation are upper triangular due to \eqref{eq:equivariance}. In fact, this property can be shown to be equivalent to the network having feedforward structure. However, for most of our considerations it is more convenient not to use the tensor notation of linear maps.
	\hspace*{\fill}$\triangle$
\end{remk}
\noindent
The structure of linear admissible maps immediately implies the following result on generic $1$-parameter bifurcations, such as the Hopf bifurcation, in feedforward networks.
\begin{theorem}[\cite{vonderGracht.2022}, Rem. 2.12]
	\label{prop:ffhopf}
	In a $1$-parameter bifurcation in a feedforward network with one-dimensional internal dynamics there cannot be a pair of conjugate imaginary eigenvalues at the synchronous bifurcation point. On the other hand, if the internal dynamics is at least $2$-dimensional, a $1$-parameter bifurcation in which a pair of complex eigenvalues crosses the imaginary axis is possible.
\end{theorem}
\begin{proof}
	In order for a bifurcation to occur, a $1$-parameter family of linear admissible maps has to have an eigenvalue/a pair of complex conjugate eigenvalues that crosses/cross the imaginary axis at the bifurcation point. In the case $V=\RR$ all linear admissible maps are real upper triangular matrices. Their eigenvalues are the diagonal elements which are real. Hence, only real eigenvalues can cross the imaginary axis. On the other hand, when the internal dynamics is in $V\cong \RR^d$ with $d\ge 2$, linear admissible maps are upper triangular with entries in $\gl{V}$. Thus the eigenvalues of a linear admissible map are the union of the eigenvalues of all diagonal elements which are arbitrary elements in $\gl{V}$ (some of which might be related). In particular, there are possible diagonal elements with complex eigenvalues.
\end{proof}
\begin{remk}
	In particular, the emergence (or collapse) of periodic solutions in a bifurcation -- as in classical and non-classical Hopf bifurcations -- requires a pair of complex conjugate eigenvalues to cross the imaginary axis. The previous theorem shows that this can only occur in feedforward networks, if the internal dynamics is at least $2$-dimensional.
\hspace*{\fill}$\triangle$
\end{remk}

Finally, we restate the following definition from \textcite{vonderGracht.2022}. It proves to be useful for the investigation of bifurcations. It induces an equivalence relation on the cells of a feedforward network whose equivalence classes are in one-to-one correspondence with the eigenvalues of linear admissible maps.
\begin{defi}
	Given a homogeneous coupled cell network, we denote $\LL_\sigma = \{\tau \in \Sigma\, |\,  \tau\sigma = \sigma\}$ for all $\sigma \in \Sigma$ and define an equivalence relation $\multimapboth$ on the input maps as follows:
	\begin{equation}
		\label{eq:deflooptype}
		\sigma \multimapboth \sigma' \quad \iff \quad \LL_\sigma = \LL_{\sigma'} \, .
	\end{equation}
	If $\sigma \multimapboth \sigma'$ we say that $\sigma$ and $\sigma'$ have the same \textit{loop-type}. In a network, two nodes have the same loop-type if and only if they have the same self-loops (of the same color).
	\hspace*{\fill}$\triangle$
\end{defi}
\begin{remk}
	\label{rem:maximal}
	As the network consists of only finitely many cells, there are well-defined maximal elements with respect to $\nole$. These are all of the same loop-type.
	\hspace*{\fill}$\triangle$
\end{remk}

\subsection[Summary of the case \texorpdfstring{\od}{1D}]{Summary of the case \od}
\label{subsec:ffod}
In \textcite{vonderGracht.2022}, the generic steady state bifurcations in a $1$-parameter family of vector fields for a feedforward network are thoroughly investigated. It is shown that feedforward structure induces a so-called \emph{amplification effect} in bifurcating branches, i.e., the branching solution in a given cell has a steeper slope the `lower' the cell is in the network with respect to $\nole$. These results in particular apply to fundamental networks. The bifurcation problem is as laid out in \Cref{sec:bi}. We investigate parameter-dependent dynamics governed by
\begin{equation}
	\label{eq:ffparam}
	\dot{v}_\sigma = f (v_{\tau_1\sigma}, \dotsc, v_{\tau_n\sigma}, \lambda),
\end{equation}
for all $\sigma \in \Sigma$, where $\lambda \in \RR$. Due to the feedforward structure, the equation for component $v_\sigma$ depends only on those states $v_{\sigma'}$ for $\sigma' \noge \sigma$. Furthermore, we assume $f(0,0) = 0$ and the linearization $D_v\Gamma_f (0,0)$ at the bifurcation point to have an eigenvalue $0$. In accordance to \Cref{rem:pot}, the linearization takes the form
\begin{equation}
	\label{eq:lin}
	D_v \Gamma_f (0,0) = \begin{pmatrix}
		\sum_{\tau \in \LL_{\sigma_1}} a_\tau	& \bullet								& \dots		& \bullet	\\
		0										& \sum_{\tau \in \LL_{\sigma_2}} a_\tau	& \ddots	& \vdots 	\\
		\vdots									& \ddots								& \ddots	& \bullet	\\
		0										& \dots									& 0			& \sum_{\tau \in \LL_{\sigma_N}} a_\tau
	\end{pmatrix},
\end{equation}
where $a_\tau = \partial_\tau f(0,0)$ is the partial derivative in the direction of the $\tau$-input. 
\begin{remk}
	In the tensor notation, this linear map has the form
	\[ D_v \Gamma_f (0,0) = \sum_{\sigma \in \Sigma} B_\sigma \ox a_\sigma, \]
	as in \Cref{prop:linadmtensor} and \Cref{rem:pot}. Here the maps $B_\sigma$ are upper triangular. However, as we are only interested in the diagonal elements, the non-tensor notation is more convenient.
	\hspace*{\fill}$\triangle$
\end{remk}
The eigenvalues of the linearization \eqref{eq:lin} are determined per cell from the diagonal elements $\sum_{\sigma \in \LL_{\overline{\sigma}}} a_\sigma$. Their multiplicities are given by the number of elements of the same loop-type, as loop-type equivalent elements yield the same diagonal elements of the linearization. Hence, $D_v \Gamma_f(0,0)$ is non-invertible, if and only if there is a node $\overline{\sigma} \in \Sigma$ such that
\begin{equation}
	\label{eq:critical}
	0 \in \operatorname{spec}\left( \sum_{\sigma \in \LL_{\overline{\sigma}}} a_\sigma \right).
\end{equation}
Then the same holds for all nodes $\hat{\sigma}\in \Sigma$ such that $\hat{\sigma}\multimapboth \overline{\sigma}$. We call these nodes, as well as their loop-type, \emph{critical}. All other cells and their loop-types are referred to as \emph{non-critical}. Note that generically precisely one loop-type is critical while all others are not. In particular, in the light of \Cref{rem:maximal} either all maximal cells are critical (and no other cell is) or no maximal cell is critical. The bifurcation problem is to find solutions to
\[ \Gamma_f(v,\lambda) = 0 \]
locally around the bifurcation point $(0,0)$.

In the remainder of this subsection, we summarize the results for the bifurcation problem stated above in the case \od as is provided in \textcite{vonderGracht.2022}. Once again, we replace the coordinates describing the internal state of each cell $v_\sigma \in V$ by coordinates $x_\sigma \in \RR$. Note that in this case the linear admissible maps are truly upper triangular and the eigenvalue condition \eqref{eq:critical} becomes $\sum_{\sigma \in \LL_{\overline{\sigma}}} a_\sigma = 0$ for critical cells $\overline{\sigma}$. Generically, four different types of branching steady states may occur depending on whether the maximal cells are critical or not and whether the branching occurs for $\lambda >0$ or $\lambda<0$. The following results are a summary of corresponding ones in \textcite{vonderGracht.2022}.

\begin{prop}[\textcite{vonderGracht.2022}, Thm. 3.7]
	\label{prop:maxcrit}
	Branching steady state solutions of a feedforward fundamental network with maximal critical cells, generically are as in one of the following two cases:
	\begin{enumerate}[label=(\roman*)]
		\item \hfill (supercritical) \[ x_\sigma(\lambda) = D_\sigma \cdot \sqrt{\lambda} + \OO(|\lambda|) \quad \text{for small} \quad \lambda>0; \]
		\item \hfill (subcritical) \[ x_\sigma(\lambda) = D_\sigma \cdot \sqrt{- \lambda} + \OO(|\lambda|) \quad \text{for small} \quad \lambda<0; \]
	\end{enumerate}
	where $D_\sigma\in\RR\setminus\{0\}$. The direction of branching is the same for all $\sigma \in \Sigma$.
\end{prop}
\begin{prop}[\textcite{vonderGracht.2022}, Thms. 3.21, 3.23 \& 3.24]
	\label{prop:maxnoncrit}
	Branching steady state solutions in a feedforward fundamental network with non-maximal critical cells are as follows: There is a \emph{root subnetwork}, i.e., a subset of nodes $\Bsup \subset \Sigma$ containing all maximal cells that is not influenced from outside of $\Bsup$ meaning $\sigma \Bsup \subset \Bsup$ for all $\sigma\in\Sigma$. The state variables of all $\sigma\in\Bsup$ bifurcate as
	\[ x_\sigma ( \lambda ) = X(\lambda) = D \cdot \lambda + \OO (|\lambda|^2) \]
	for $|\lambda|$ small. The remaining cells bifurcate as in one of the following two cases
	\begin{enumerate}[label=(\roman*)]
		\item \hfill (supercritical) \[ x_\sigma ( \lambda ) = \Dsup_\sigma \cdot \lambda^{2^{-\musup_\sigma}} + \OO \left( |\lambda|^{2^{-(\musup_\sigma-1)}} \right) \quad \text{for small} \quad \lambda>0; \]
		\item \hfill (subcritical) \[ x_\sigma ( \lambda ) = \Dsub_\sigma \cdot (-\lambda)^{2^{-\musub_\sigma}} + \OO \left( |\lambda|^{2^{-(\musub_\sigma-1)}} \right) \quad \text{for small} \quad \lambda<0; \]
	\end{enumerate}
	where $\mu_\sigma$ is defined inductively setting $\mu_\sigma=0$ for all $\sigma \in B$ and
	\begin{equation}
		\label{eq:defmu}
		\musup_\sigma = \begin{cases}
			\max_{\tau \nog \sigma} \musup_\tau		&\text{for} \quad \sigma \text{ non-critical}, \\
			\max_{\tau \nog \sigma} \musup_\tau		&\text{for} \quad \sigma \text{ critical and } \tau\in B \text{ for all } \tau\nog p, \\
			\max_{\tau \nog \sigma} \musup_\tau + 1	&\text{for} \quad \sigma \text{ critical and there exists } \tau\nog p \text{ such that } \tau\notin B.
		\end{cases}
	\end{equation}		
	That is, $\mu_\sigma$ is the maximal number of critical cells $\tau \notin \Bsup$ along paths from any cell in $\Bsup$ to $\sigma$. Furthermore, $D,\Dsub_\sigma\in\RR\setminus\{0\}$ for all $\sigma \in \Sigma$. The direction of branching is the same for all $\sigma \notin \Bsup$.
\end{prop}
\begin{remk}
	\label{remk:muextended}
	We may extend the definition of $\musup_\sigma$ to the case of critical maximal cells. Therein the state variables of all cells bifurcate as
	\begin{align*}
		x_\sigma(\lambda)	&= D_\sigma \cdot \lambda^{2^{-\musup_\sigma}} + \OO(|\lambda|)	\quad \text{or} \\
		x_\sigma(\lambda)	&= D_\sigma \cdot (-\lambda)^{2^{-\musup_\sigma}} + \OO(|\lambda|)
	\end{align*}
	depending on sub- or supercriticality of the solution branch, where $\musup_\sigma=1$ for all $\sigma\in\Sigma$. This also satisfies the growing condition \eqref{eq:defmu}, as only the maximal cells are critical and these receive all their inputs from themselves.
\end{remk}
\begin{remk}
	\label{remk:incomplete}
	The results in \textcite{vonderGracht.2022} are more complete than what we state here. They also include formulas to explicitly compute the lowest order coefficients. Furthermore, we can deduce a condition on the system parameters -- the partial derivatives of the governing function $f$ -- that determines whether a specific branch exists. In particular, the branching pattern may differ throughout parameter space. We refer to a bifurcating branch as being \emph{generic} if it occurs for an open set of system parameters. \Cref{prop:maxcrit,prop:maxnoncrit} describe all branching solutions for an open and dense subset of system parameters. The details are omitted here because we only need the general form of all branching solutions.
	\hspace*{\fill}$\triangle$
\end{remk}

Even if we do not state the full results here, we can see the main characteristic of the bifurcating solutions, which is the amplification effect. When we investigate the dynamics cell-by-cell starting with the maximal cells, we observe that in each critical cell -- if the maximal cells are critical, no other cell is -- the square root order of the lowest order terms in $\lambda$ increases. Hence, the branch exits $0$ with a steeper slope. This is independent of the branching direction. As indicated in \Cref{remk:incomplete}, in bifurcation analysis one is often not interested in the explicit expression describing how to compute the branching state for every parameter value but rather in the qualitative behavior. To that end, we consider the asymptotics of the state variable for each cell separately in a generic steady state bifurcation which encodes qualitative bifurcation information of each individual cell. \Cref{prop:maxcrit,prop:maxnoncrit} allow to describe these asymptotics for each cell separately in a generic steady state bifurcation in the case \od. We introduce some notation to state these ideas more precisely and to make them practical for the case \dd as well. The following definition is independent of the feedforward structure and could also have been formulated for non-fundamental networks. We state it for arbitrary internal dynamics on $V$ from which we may deduce the case \od as a special case.
\begin{defi}
	Let $(v_\sigma (\lambda))_{\sigma\in\Sigma} \subset \bigoplus_{\sigma\in\Sigma}V$ be a branch of steady states for a parameter-dependent feedforward fundamental network vector field for small $\lambda>0, \lambda<0$ or $|\lambda|$ small, i.e.,
	\begin{equation}
		\label{eq:bif}
		\Gamma_f((v_\sigma (\lambda))_{\sigma\in\Sigma}, \lambda) = 0.
	\end{equation}
	We say that a cell $\sigma \in \Sigma$ \emph{has the square root order $\xi_\sigma \ge 0$ in $\lambda$} or that it \emph{grows asymptotically as $\lambda^{2^{-\xi_\sigma}}$} if
	\[ \begin{cases}
		\|v_\sigma(\lambda) \| = d_\sigma \cdot \lambda^{2^{-\xi_\sigma}} + \OO \left( |\lambda|^{2^{-(\xi_\sigma+1)}} \right) \quad \text{with } d_\sigma \ne 0 	&\text{for} \quad \xi_\sigma \ge 1, \\[5pt]
		\|v_\sigma(\lambda) \| = \OO \left( |\lambda| \right)	&\text{for} \quad \xi_\sigma = 0.
	\end{cases} \]
	Therein $\| \cdot \|$ denotes the euclidean norm on $V$. If the branch only exists for $\lambda<0$ replace $\lambda$ by $-\lambda$. We denote this situation by $v_\sigma \sim \lambda^{2^{-\xi_\sigma}}$.
	\hspace*{\fill}$\triangle$
\end{defi}
\begin{remk}
	The situation $v_\sigma \sim \lambda^{2^{-\xi_\sigma}}$ with $\xi_\sigma >0$ is equivalent to
	\[ v_\sigma(\lambda) = \lambda^{2^{-\mu_\sigma}} \cdot \vartheta_\sigma + R_\sigma (\lambda), \]
	where $\vartheta_\sigma \in V\setminus \{0\}$ suitable and $R_\sigma\colon\RR \to V$ -- restricted to the suitable neighborhood of $\lambda_0=0$ -- such that $\|R_\tau(\lambda)\| = \OO \left( |\lambda|^{2^{-(\mu_\tau-1)}} \right)$.
	\hspace*{\fill}$\triangle$
\end{remk}

The key observation for the investigation of implications of the steady state bifurcation results in the case \od for the case \dd is the following:
\begin{prop}
	\label{prop:asorder}
	Let $(x_\sigma(\lambda))_{\sigma\in\Sigma}$ be a branch of steady states for a generic steady state bifurcation of a feedforward fundamental network as described in \Cref{prop:maxcrit,prop:maxnoncrit} with cell-by-cell asymptotics $x_\sigma \sim \lambda^{2^{-\mu_\sigma}}$ for all $\sigma \in \Sigma$. Then the asymptotic orders of the cells are partially ordered with respect to $\nole$, i.e.,
	\begin{equation}
		\label{eq:asorder}
		\sigma \noge \tau \qquad \text{implies} \qquad \mu_\sigma \le \mu_\tau.
	\end{equation}
	Hence, if $\tau \nole \sigma$, then $x_\tau$ has at most the same asymptotic order as $x_\sigma$.
\end{prop}
\begin{proof}
	In the case of critical maximal cells, all cells have square root order $1$, i.e., $x_\sigma \sim \sqrt{\lambda}$. Hence, nothing is to be shown. If, on the other hand, the maximal cells are non-critical, we obtain a subset of cells $B \subset \Sigma$ that is not influenced by any cells outside of $B$ and such that all cells in $B$ have square root order $0$: $x_\sigma \sim \lambda$ -- or $\mu_\sigma = 0$ -- for all $\sigma \in B$. In particular, $\sigma B \subset B$ for all $\sigma\in\Sigma$ implies that $B$ contains all maximal cells. All $\sigma\notin B$ have square root order $\mu_\sigma$, where $\mu_\sigma$ is defined as in \eqref{eq:defmu}. By definition $\mu_\sigma \le \mu_\tau$ if $\sigma \nog \tau$.
\end{proof}

\subsection{Feedforward networks with high-dimensional internal dynamics}
\label{subsec:ffdd}
In this section, we show that feedforward structure can be used to observe the amplification effect also in generic steady state bifurcations when the internal dynamics is high-dimensional. In particular, for a fixed feedforward fundamental network, the amplification effect we observe in a generic $1$-parameter steady state bifurcation is the same independent of the dimension of the internal dynamics. The most important tools are the partial order on the cells and the fact that the admissible maps are upper triangular -- or more precisely, they respect the partial order. As the center manifold reduction for fundamental networks respects monoid symmetry, this in particular holds for the maps involved in the procedure that allows to translate a generic branch in the case \od into a generic branch in the case \dd presented in \Cref{subsubsec:ddbranch}. As equivariance is equivalent to admissibility \eqref{eq:equivariance}, this has the convenient effect that whenever we compute some property of the state variable $v_\sigma$ of cell $\sigma$ it only depends on the state variables of cells $\tau \noge \sigma$. This observation proves to be powerful enough to translate the amplification effect into the high-dimensional case without having to determine center manifolds explicitly. In particular, we show that critical cells have the same square root order in any dimension of the internal dynamics. In a special case we obtain the same cell-by-cell asymptotics for all cells.

A crucial part of this section contains the technical analysis of the maps that are needed to define the center manifold reduction and their interaction with branching steady state solutions. In particular, we will make heavy use of the technicalities proved in \Cref{subsubsec:cmfr}. To that end, we reuse the notation from \Cref{subsec:in}. We investigate steady state solutions of
\[ \Gamma_f \colon \bigoplus_{\sigma\in\Sigma}V \times \RR \to \bigoplus_{\sigma\in\Sigma}V, \]
which is an admissible vector field for a feedforward fundamental network that depends on a real parameter $\lambda\in\RR$ as in \eqref{eq:ffparam}, close to the bifurcation point $(v_0, \lambda_0)=(0,0)$, as in the beginning of \Cref{sec:ff}. That is, we are interested in solutions to
\[ \Gamma_f(v,\lambda) = 0, \]
close to $(0,0)$, where
\[ D_v \Gamma_f (0,0) \]
is non-invertible. As stated before, for technical reasons, we focus on the extended system as in \eqref{eq:extended}
\begin{equation}
	\label{eq:ffext}
	\dot{\vu} = \begin{pmatrix} \dot{v} \\ \dot{\lambda} \end{pmatrix} = \begin{pmatrix} \Gamma_f (v, \lambda) \\ 0 \end{pmatrix} = \Gu_f(\vu)
\end{equation}
on $\bigoplus_{\sigma\in\Sigma}V \times \RR$. The extended system is equivariant with respect to the extended right regular representation $\sigma \mapsto \uu{A}_\sigma$ where
\[ \uu{A}_\sigma \colon \uu{v} = (v,\lambda) \mapsto (A_\sigma v, \lambda). \]
Under the bifurcation assumption the linearization of the extended system induces a splitting of $\bigoplus_{\sigma\in\Sigma}V \times \RR$ into subrepresentations
\begin{equation}
	\label{eq:cspextendedff}
	\begin{split}
		\uu{\XX}^h	&= \redim{D_v \Gamma_f(0,0)} \times \{ 0 \}, \\
		\uu{\XX}^c	&= \gker{D_v \Gamma_f(0,0)} \times \{ 0 \} \oplus \left\langle (v^h, 1) \right\rangle,
	\end{split}
\end{equation}
where $v^h\in \redim{D_v \Gamma_f(0,0)}$ suitable, as in \eqref{eq:cspextended}. These subspaces are generalized kernel and reduced image of $D\uu{\Gamma}_f(0) = D_{(v,\lambda)}\uu{\Gamma}_f(0,0)$, respectively. The corresponding equivariant projections are again denoted by $\uu{P}^c$ and $\uu{P}^h$, respectively.
In particular, every element $\uu{v}^c \in \uu{\XX}^c$ can be represented as
\[ \uu{v}^c = (v^c,0) + \lambda (v^h, 1), \]
where $v^c \in \gker{D_v \Gamma_f(0,0)}$ and $\lambda\in\RR$. Furthermore, the isomorphism of monoid representations \mbox{$P' \colon \uu{\XX}^c \to \gker{D_v \Gamma_f(0,0)} \times \RR$} is given by
\begin{equation}
	\label{eq:Pprimeff}
	P' (\uu{v}^c) = (P^c(v^c+\lambda v^h), \lambda) = (v^c,\lambda)
\end{equation}
and its inverse is given by
\begin{equation}
	\label{eq:Qprimeff}
	Q' (v^c,\lambda) = (v^c,0) + \lambda (v^h, 1) = \uu{v}^c
\end{equation}
(compare to \eqref{eq:Pprime} and \Cref{lem:Qprime}). Then the extended system admits the local center manifold
\begin{equation}
	\label{eq:cmf}
	M^c = \left\{ \left( \mathbbm{1}_{\uu{\XX}^c} + \uu{\psi} \right) (\uu{v}^c) \mid \uu{v}^c \in \uu{\XX}^c \right\} = \left\{ \left( \mathbbm{1}_{\uu{\XX}^c} + \uu{\psi} \right) (Q'(v^c,\lambda)) \mid (v^c,\lambda) \in \gker{D_v \Gamma_f(0,0)} \times \RR \right\},
\end{equation}
which is represented as a graph over either $\uu{\XX}^c$ or equivalently over $\gker{D_v \Gamma_f(0,0)} \times \RR$ by the equivariant map $\uu{\psi}$.

For any branch of steady states of the bifurcation problem we have $\left( (v_\sigma (\lambda))_{\sigma\in\Sigma}, \lambda \right) \subset M^c$. Furthermore, the dynamics restricted to the center manifold is bijectively conjugate to dynamics on the center subspace $\uu{\XX}^c$ or on $\gker{D_v \Gamma_f(0,0)} \times \RR$ (compare to \eqref{eq:cmfred}). Hence, any branch can uniquely be represented in coordinates of these subspaces using the projections $\uu{P}^c$ or \mbox{$P = P' \circ \uu{P}^c$} such that
\begin{align*}
	\uu{P}^c \left( (v_\sigma (\lambda))_{\sigma\in\Sigma}, \lambda \right)						&\subset \uu{\XX}^c, \\
	P' \left( \uu{P}^c \left( (v_\sigma (\lambda))_{\sigma\in\Sigma}, \lambda \right) \right) 	&\subset \gker{D_v \Gamma_f(0,0)} \times \RR
\end{align*}
(compare to \Cref{subsubsec:ddbranch}). Recall from \eqref{eq:Pprime} and \Cref{lem:1,lem:2,lem:2.5,lem:Qprime} that all maps required to switch between the three representation of a branch of steady states -- i.e., $\uu{P}^c, P', Q'$ and $\uu{\psi}$ leave the $\lambda$\nobreakdash-coordinate unchanged (the third representation is the original representation \mbox{$\left( (v_\sigma (\lambda))_{\sigma\in\Sigma}, \lambda \right)$}). Hence, we see that the parameter of each representation remains the same and we may write
\begin{equation}
	\label{eq:repr}
	\begin{split}
		\left( (y_\sigma (\lambda))_{\sigma\in\Sigma}, \lambda \right)	&= \uu{P}^c \left( (v_\sigma (\lambda))_{\sigma\in\Sigma}, \lambda \right) \subset \uu{\XX}^c, \\
		\left( (z_\sigma (\lambda))_{\sigma\in\Sigma}, \lambda \right)	&= P' \left( \uu{P}^c \left( (v_\sigma (\lambda))_{\sigma\in\Sigma}, \lambda \right) \right) \subset \gker{D_v \Gamma_f(0,0)} \times \RR.
	\end{split}
\end{equation}
Note that the second representation of \eqref{eq:repr} is of particular interest as it describes the branch of steady states in the corresponding reduced bifurcation problem on $\gker{D_v \Gamma_f(0,0)} \times \RR$ (see \eqref{eq:cmfred}).

We investigate to what extent cell-by-cell asymptotics are respected when switching between the different representations. In particular, we answer the following question: Under the assumption of square root orders of individual cells that are partially ordered as in \Cref{prop:asorder} for feedforward networks -- more precisely as in \eqref{eq:asorder} -- in either one of the representations, can we recover the square root order of each cell in any of the other representations? As the distinction of cells is reflected in the choice of coordinates, this is done by analyzing the projection operators in the original coordinates.

\begin{lemma}
	\label{lem:3}
	\begin{enumerate}[label=(\roman*)]
		\item Let $\left( (v_\sigma (\lambda))_{\sigma\in\Sigma}, \lambda \right) \subset M^c$ be a branch of steady states of \eqref{eq:ffext}. Assume cell-by-cell asymptotics $v_\sigma \sim \lambda^{2^{-\mu_\sigma}}$ that are partially ordered as in \eqref{eq:asorder}, i.e., $\sigma \noge \tau$ implies $\mu_\sigma \le \mu_\tau$. Then the representation in the generalized kernel of the extended system 
		\[ \left( (y_\sigma (\lambda))_{\sigma\in\Sigma}, \lambda \right) = \uu{P}^c \left( (v_\sigma (\lambda))_{\sigma\in\Sigma}, \lambda \right) \subset \uu{\XX}^c \]
		exhibits the same cell-by-cell asymptotics, i.e.,
		\[ y_\sigma \sim \lambda^{2^{-\mu_\sigma}} \quad \text{for all} \quad \sigma \in \Sigma. \]
		\item Let $\left( (y_\sigma (\lambda))_{\sigma\in\Sigma}, \lambda \right) \subset \uu{\XX}^c$ be the representation of a branch of steady states of \eqref{eq:ffext} in the generalized kernel of the extended system. Assume cell-by-cell asymptotics $y_\sigma \sim \lambda^{2^{-\mu_\sigma}}$ that are partially ordered as in \eqref{eq:asorder}, i.e., $\sigma \noge \tau$ implies $\mu_\sigma \le \mu_\tau$. Then the representation in the full space
		\[ \left( (v_\sigma (\lambda))_{\sigma\in\Sigma}, \lambda \right) = \left( \mathbbm{1}_{\uu{\XX}^c} + \uu{\psi} \right) \left( (y_\sigma (\lambda))_{\sigma\in\Sigma}, \lambda \right) \subset M^c \]
		exhibits the same cell-by-cell asymptotics, i.e.,
		\[ v_\sigma \sim \lambda^{2^{-\mu_\sigma}} \quad \text{for all} \quad \sigma \in \Sigma. \]
	\end{enumerate}
\end{lemma}
The proof of the lemma is quite technical but depends on a rather simple observation:  The maps used to change from one representation of a given branch to another are (in their spatial component) equivariant with respect to the monoid representation. As equivariance is equivalent to being admissible with respect to the feedforward network structure, this pairs very nicely with the the assumed ordering of the cell-by-cell asymptotics. In particular, the $\tau$-component of the target representation depends only on the components of the branch which are above $\tau$ with respect to $\nole$. By assumption, these grow with at most the same asymptotic square root order as the $\tau$-component. As the maps are further linear to leading order, the leading order in the target cell cannot be higher than that of the original cell.
The bulk part of the proof consists of formalizing this argument for both parts of the lemma and to ensure that the leading order is in fact as stated by proving that leading coefficients do not vanish.
\begin{proof}
	The central observation to prove both parts of the lemma is the fact that the maps $\left. \uu{P}^c \right|_{M^c} \colon M^c \to \uu{\XX}^c$ and $\left( \mathbbm{1}_{\uu{\XX}^c} + \uu{\psi} \right) \colon \uu{\XX}^c \to M^c$ are mutually inverse bijections between the center manifold $M^c$ and the generalized kernel of the extended system $\uu{\XX}^c$. The properties of these maps are well-understood from \Cref{lem:1,lem:2,lem:2.5}. We know that 
	\[ \uu{P}^c(v,\lambda) = \left( \uu{P}^c_v (v,\lambda), \lambda \right),\]
	where $\uu{P}^c_v$ can be interpreted as a parameter-dependent linear map on $\bigoplus_{\sigma\in\Sigma}V$ that is equivariant with respect to the non-extended monoid representation, i.e., $\uu{P}^c_v (A_\sigma v,\lambda) = A_\sigma \uu{P}^c_v (v,\lambda)$ for all \mbox{$v\in\bigoplus_{\sigma\in\Sigma}V ,\lambda\in\RR$} and $\sigma\in\Sigma$. Due to the equivalence of equivariance and admissibility, it is therefore a parameter-dependent linear admissible map of the fundamental network vector field. In particular, it respects the partial order $\nole$
	\begin{equation}
		\label{eq:brcmfgker-linform}
		\left( \uu{P}^c_v (v,\lambda) \right)_\tau = \ell_\tau \left( (v_\sigma \mid \sigma \noge \tau), \lambda \right),
	\end{equation} 
	where $\ell_\tau$ is a linear map that depends only on the entries $v_\sigma$ for $\sigma \noge \tau$ and the parameter $\lambda$. The parameter-dependence is linear as well.
	
	A similar observation holds for the other direction. Denote elements of the generalized kernel of the extended system by $\left( y,\lambda \right) = \left( (y_\sigma)_{\sigma\in\Sigma}, \lambda \right) \in \uu{\XX}^c \subset \bigoplus_{\sigma\in\Sigma} V \times \RR$. Then we know \mbox{$\uu{\psi} ( y, \lambda ) = ( \uu{\psi}_v (y,\lambda), 0 )$}, where $\uu{\psi}_v$ is a parameter-dependent admissible map of the fundamental network. That is
	\begin{equation}
		\label{eq:brcmfgker-quadform}
		\left( \uu{\psi}_v (y, \lambda) \right)_\tau = q_\tau \left( (y_\sigma \mid \sigma \noge \tau), \lambda \right).
	\end{equation}
	As $\uu{\psi}(0,0) = 0$ and $D \uu{\psi}(0,0) = 0$, we also obtain $q_\tau (0,0) = 0$ and $d q_\tau (0,0) = 0$.
	
	We may now prove both cases separately beginning with (\romannumeral 2). Let $\left( (y_\sigma (\lambda))_{\sigma\in\Sigma}, \lambda \right) \subset \uu{\XX}^c$ be the representation of a branch of steady states in the generalized kernel of the extended system. Assume cell-by-cell asymptotics $y_\sigma \sim \lambda^{2^{-\mu_\sigma}}$ that satisfy the ordering assumption \eqref{eq:asorder}, i.e., $\sigma \noge \tau$ implies $\mu_\sigma \le \mu_\tau$. First, we fix a cell $\tau \in \Sigma$ with $\mu_\tau \ge 1$. By definition $y_\tau \sim \lambda^{2^{-\mu_\tau}}$ is equivalent to
	\[ y_\tau (\lambda) = \lambda^{2^{-\mu_\tau}} \cdot \vartheta_\tau + R_\tau (\lambda), \]
	where $\vartheta_\tau \in V\setminus \{0\}$ suitable and $\|R_\tau(\lambda)\| = \OO \left( |\lambda|^{2^{-(\mu_\tau-1)}} \right)$. Combining the partial order of the cell-by-cell asymptotics with the form of the map whose graph is the center manifold \eqref{eq:brcmfgker-quadform}, and the fact that $q_\tau$ is at least quadratic to leading order, we compute
	\begin{equation}
		\label{eq:brcmfgker-iiquad}
		\left( \uu{\psi}_v \left( (y_\sigma (\lambda)_{\sigma\in\Sigma}, \lambda) \right) \right)_\tau = q_\tau \left( (y_\sigma(\lambda) \mid \sigma \noge \tau), \lambda \right) = \lambda^{2^{-(\mu_\tau-1)}} \cdot \eta_\tau + \rho_\tau (\lambda),
	\end{equation}
	where $\eta_\tau \in V$ suitable and $\|\rho_\tau (\lambda)\| = \OO \left( |\lambda|^{2^{-(\mu_\tau-1)}} \right)$.
	
	As we obtain the representation of the branch on the center manifold $\left( (v_\sigma(\lambda))_{\sigma\in\Sigma}, \lambda \right)$ from
	\[ \left( (v_\sigma(\lambda))_{\sigma\in\Sigma}, \lambda \right) = \left( \mathbbm{1}_{\uu{\XX}^c} + \uu{\psi} \right) \left( (y_\sigma (\lambda)_{\sigma\in\Sigma}, \lambda) \right) = \left( (y_\sigma (\lambda))_{\sigma\in\Sigma} + \uu{\psi}_v ((y_\sigma (\lambda))_{\sigma\in\Sigma}, \lambda), \lambda \right), \]
	we directly see
	\[ v_\tau(\lambda) = \lambda^{2^{-\mu_\tau}} \cdot \vartheta_\tau + R_\tau (\lambda) + \lambda^{2^{-(\mu_\tau-1)}} \cdot \eta_\tau + \rho_\tau (\lambda) = \lambda^{2^{-\mu_\sigma}} \cdot \vartheta_\tau + \varkappa_\tau (\lambda), \]
	where $\|\varkappa_\tau(\lambda)\| = \OO \left( |\lambda|^{2^{-(\mu_\tau-1)}} \right)$. Hence, $v_\tau \sim \lambda^{2^{-\mu_\tau}}$, since $\vartheta_\tau \ne 0$.
	
	On the other hand, consider a cell $\tau\in\Sigma$ with $\mu_\tau=0$. By the partial order of the square root orders \eqref{eq:asorder} we also obtain $\mu_\sigma=0$ for all $\sigma \noge \tau$. By definition this is equivalent to $\|y_\sigma(\lambda)\| = \OO (|\lambda|)$ for all $\sigma \noge \tau$ including $\tau$. Thus, using \eqref{eq:brcmfgker-quadform} as before, we compute 
	\[ \| v_\tau(\lambda) \| = \left\| y_\tau(\lambda) + q_\tau \left( (y_\sigma (\lambda) \mid \sigma \noge \tau), \lambda \right) \right\| = \OO (|\lambda|), \]
	which proves $v_\tau \sim \lambda$ and, therefore, completes the proof for (\romannumeral 2).
	
	Next we turn to the proof of the first statement (\romannumeral 1). This is slightly more complicated as we have to make sure that the projection onto the generalized kernel of the extended system does not `loose' information about the asymptotics of a given cell. Assume that $\left( (v_\sigma (\lambda))_{\sigma\in\Sigma}, \lambda \right) \subset M^c$ is a branch of steady states whose cell-by-cell asymptotics satisfy the ordering assumption \eqref{eq:asorder}. Then, fixing a cell $\tau \in \Sigma$ with $\mu_\tau \ge 1$, we obtain
	\[ v_\tau (\lambda) = \lambda^{2^{-\mu_\tau}} \cdot \vartheta_\tau + R_\tau (\lambda), \]
	with $\vartheta_\tau \in V\setminus \{0\}$ suitable and $\|R_\tau(\lambda)\| = \OO \left( |\lambda|^{2^{-(\mu_\tau-1)}} \right)$. Combining the partial order of the cell-by-cell asymptotics with the specific form of the projection onto $\uu{\XX}^c$ \eqref{eq:brcmfgker-linform}, in which $\ell_\tau$ is linear, we obtain
	\begin{equation}
		\label{eq:brcmfgker-iproj}
		\left( \uu{P}_v^c \left( (v_\sigma (\lambda)_{\sigma\in\Sigma}, \lambda) \right) \right)_\tau = \ell_\tau \left( (v_\sigma (\lambda) \mid \sigma \noge \tau), \lambda \right) = \lambda^{2^{-\mu_\tau}} \cdot \eta_\tau + \rho_\tau (\lambda),
	\end{equation}
	where $\eta_\tau \in V$ suitable and $\|\rho_\tau (\lambda)\| = \OO \left( |\lambda|^{2^{-(\mu_\tau-1)}} \right)$. The representation of the branch in the generalized kernel of the extended system $\left( (y_\sigma(\lambda))_{\sigma\in\Sigma}, \lambda \right)$ is computed as
	\[ \left( (y_\sigma(\lambda))_{\sigma\in\Sigma}, \lambda \right) = \uu{P}^c \left( (v_\sigma(\lambda))_{\sigma\in\Sigma}, \lambda \right) = \left( \uu{P}^c_v \left( (v_\sigma(\lambda))_{\sigma\in\Sigma}, \lambda \right), \lambda \right). \]
	Thus, we see that $y_\tau \sim \lambda^{2^{-\mu_\tau}}$, if $\eta_\tau \ne 0$.\footnote{The opposite case $\eta_\tau = 0$ is what we refer to as `loosing' information about the asymptotics in cell $\tau$.}
	
	However, as $\left( (y_\sigma(\lambda))_{\sigma\in\Sigma}, \lambda \right) \subset \uu{\XX}^c$ is the representation of a branch of steady states in the generalized kernel of the extended system, we may apply the results from (\romannumeral 2) to it and receive its corresponding representation on the center manifold $\left( (\tilde{v}_\sigma(\lambda))_{\sigma\in\Sigma}, \lambda \right) \subset M^c$ with the same cell-by-cell asymptotics. As in the proof of part (\romannumeral 2), we compute its $v$-component to be
	\begin{equation}
		\label{eq:brcmfgker-iquad}
		\tilde{v}_\tau(\lambda) = y_\tau(\lambda) + \uu{\psi}_v ((y_\sigma (\lambda))_{\sigma\in\Sigma}, \lambda) = \lambda^{2^{-\mu_\tau}} \cdot \eta_\tau + \varkappa_\tau (\lambda),
	\end{equation}
	where $\eta_\tau \in V$ is as in \eqref{eq:brcmfgker-iproj} and $\|\varkappa_\tau (\lambda)\| = \OO \left( |\lambda|^{2^{-(\mu_\tau-1)}} \right)$. Furthermore, $\left. \uu{P}^c \right|_{M^c}$ and $\left( \mathbbm{1}_{\uu{\XX}^c} + \uu{\psi} \right)$ are mutually inverse bijections between $M^c$ and $\uu{\XX}^c$. Hence,
	\begin{align*}
		\left( (v_\sigma(\lambda))_{\sigma\in\Sigma}, \lambda \right) &= \left( \left( \mathbbm{1}_{\uu{\XX}^c} + \uu{\psi} \right) \circ {\uu{P}^c}_{| M^c} \right) \left( (v_\sigma(\lambda))_{\sigma\in\Sigma}, \lambda \right) \\
			&= \left( \mathbbm{1}_{\uu{\XX}^c} + \uu{\psi} \right) \left( (y_\sigma(\lambda))_{\sigma\in\Sigma}, \lambda \right) \\
			&= \left( (\tilde{v}_\sigma(\lambda))_{\sigma\in\Sigma}, \lambda \right).
	\end{align*} 
	Using \eqref{eq:brcmfgker-iquad}, we obtain
	\[ v_\tau(\lambda) = \lambda^{2^{-\mu_\tau}} \cdot \vartheta_\tau + R_\tau (\lambda) = \lambda^{2^{-\mu_\tau}} \cdot \eta_\tau + \varkappa_\tau (\lambda) = \tilde{v}_\tau(\lambda). \]
	In particular, $\eta_\tau = \vartheta_\tau \ne 0$, which proves $y_\tau \sim \lambda^{2^{-\mu_\sigma}}$.
	
	On the other hand, $\mu_\tau=0$ implies $\mu_\sigma=0$ and therefore $\|v_\tau(\lambda)\|= \OO(|\lambda|)$ for all $\sigma \noge \tau$ (which includes $\tau$). Thus we compute
	\[ \|y_\tau(\lambda)\| = \left\| \ell_\tau \left( (v_\sigma(\lambda) \mid \sigma \noge \tau ), \lambda \right) \right\| = \OO(|\lambda|), \]
	which shows $y_\tau \sim \lambda$, completing the proof of (\romannumeral 1).
\end{proof}
\begin{remk}
	\label{rem:relaxed3}
	The proof for \Cref{lem:3} allows for a slightly more precise statement when relaxing the ordering assumption on the cell-by-cell asymptotics \eqref{eq:asorder}. That is, for every cell $\tau$ for which the ordering assumption is fulfilled, $\mu_\sigma \le \mu_\tau$ for all $\sigma \noge \tau$, we obtain that the $\tau$-component exhibits the same asymptotics in the full representation of the solution branch and in its representation in the generalized kernel of the extended system, i.e., $v_\tau \sim \lambda^{2^{-\mu_\tau}}$ and $y_\tau \sim \lambda^{2^{-\mu_\tau}}$. This does not require the square root orders of \emph{all} cells to be ordered.
	\hspace*{\fill}$\triangle$
\end{remk}

In the previous lemma we have shown that, under the ordering assumption \eqref{eq:asorder}, cell-by-cell asymptotics agree in the representation of a steady state branch on the center manifold with those in the representation in the generalized kernel of the extended system. A similar result holds true when comparing representations in $\uu{\XX}^c$ and in $\gker{D_v \Gamma_f(0,0)}\times \RR$ even without the ordering assumption, i.e., $y_\tau \sim z_\tau$ with $y_\tau$ and $z_\tau$ as in \eqref{eq:repr}. This follows from the properties of the mutually inverse isomorphisms $P'$ and $Q'$.
\begin{lemma}
	\label{lem:4}
	\begin{enumerate}[label=(\roman*)]
		\item Let $\left( (y_\sigma (\lambda))_{\sigma\in\Sigma}, \lambda \right) \subset \uu{\XX}^c$ be the representation of a branch of steady states of \eqref{eq:ffext} in the generalized kernel of the extended system. Assume cell-by-cell asymptotics $y_\sigma \sim \lambda^{2^{-\mu_\sigma}}$. Then the representation
		\[ \left( (z_\sigma(\lambda))_{\sigma\in\Sigma}, \lambda \right) = \uu{P'} \left( (y_\sigma (\lambda))_{\sigma\in\Sigma}, \lambda \right) \subset \gker{D_v \Gamma_f(0,0)} \times \RR \]
		exhibits the same cell-by-cell asymptotics, i.e.,
		\[ z_\sigma \sim \lambda^{2^{-\mu_\sigma}} \quad \text{for all} \quad \sigma \in \Sigma. \]
		\item Let $\left( (z_\sigma(\lambda))_{\sigma\in\Sigma}, \lambda \right) \subset \gker{D_v \Gamma_f(0,0)} \times \RR$ be the representation of a branch of steady states of \eqref{eq:ffext} in $\gker{D_v \Gamma_f(0,0)}$. Assume cell-by-cell asymptotics $z_\sigma \sim \lambda^{2^{-\mu_\sigma}}$. Then the representation in the generalized kernel of the extended system
		\[ \left( (y_\sigma (\lambda))_{\sigma\in\Sigma}, \lambda \right) = Q' \left( (z_\sigma (\lambda))_{\sigma\in\Sigma}, \lambda \right) \subset \uu{\XX}^c \]
		exhibits the same cell-by-cell asymptotics, i.e.,
		\[ y_\sigma \sim \lambda^{2^{-\mu_\sigma}}	\quad \text{for all} \quad \sigma \in \Sigma. \]
	\end{enumerate}
\end{lemma}
\begin{proof}
	We prove (\romannumeral 1) and (\romannumeral 2) similarly to the proof of \Cref{lem:3}. Recall that the maps \mbox{$P'\colon (y,\lambda) \mapsto (P^c(y), \lambda)$} as in \eqref{eq:Pprimeff} and $Q'\colon (z,\lambda) \mapsto (z + \lambda v^h, \lambda)$ as in \eqref{eq:Qprimeff} form equivariant isomorphisms between $\uu{\XX}^c$ and $\gker{D_v \Gamma_f(0,0)} \times \RR$. In particular, due to \Cref{lem:1} we know that $P^c$ is equivariant with respect to the non-extended representation $\sigma \mapsto A_\sigma$. Hence, it is a linear, admissible map for the fundamental network and we obtain
	\begin{equation*}
		\left( P^c \left( (y_\sigma)_{\sigma\in\Sigma} \right) \right)_\tau = \ell_\tau \left( y_\sigma \mid \sigma \noge \tau \right),
	\end{equation*}
	in which $\ell_\tau$ is linear.
	
    Let $\left( (y_\sigma (\lambda))_{\sigma\in\Sigma}, \lambda \right) \subset \uu{\XX}^c$ be the representation of a branch of steady states of \eqref{eq:ffext} in the generalized kernel of the extended system and assume cell-by-cell asymptotics $y_\sigma \sim \lambda^{2^{-\mu_\sigma}}$. For a cell $\tau$ with $\mu_\tau \ge 1$, this is equivalent to
	\[ y_\tau (\lambda) = \lambda^{2^{-\mu_\tau}} \cdot \vartheta_\tau + R_\tau(\lambda), \]
	where $\vartheta_\tau \in V\setminus\{0\}$ and $\| R_\tau (\lambda) \| = \OO \left( |\lambda|^{2^{-(\mu_\tau-1)}} \right)$. As $\left( (z_\sigma(\lambda))_{\sigma\in\Sigma}, \lambda \right) = \left( P^c \left( (y_\sigma(\lambda))_{\sigma\in\Sigma} \right) , \lambda \right)$, we compute
	\[ z_\tau(\lambda) = \ell_\tau \left( y_\sigma(\lambda) \mid \sigma \noge \tau \right) = \lambda^{2^{-\xi_\tau}} \cdot \eta_\tau + \rho_\tau(\lambda), \]
	where $\eta_\tau \in V, \xi_\tau = \max_{\sigma \noge \tau} \mu_\sigma$ and $\| \rho_\tau (\lambda) \| = \OO \left( |\lambda|^{2^{-(\xi_\tau-1)}} \right)$. Furthermore, note that
	\[ Q' \left( (z_\sigma(\lambda))_{\sigma\in\Sigma}, \lambda \right) = \left( (z_\sigma(\lambda))_{\sigma\in\Sigma} + \lambda \cdot (v^h_\sigma)_{\sigma\in\Sigma}, \lambda \right). \]
	Therein
	\[ z_\tau (\lambda) + \lambda \cdot v^h_\tau = \lambda^{2^{-\xi_\tau}} \cdot \eta_\tau + \lambda \cdot v^h_\tau + \rho_\tau(\lambda) = \lambda^{2^{-\xi_\tau}} \cdot \eta_\tau + \varkappa_\tau(\lambda), \]
	where $\| \varkappa_\tau (\lambda) \| = \OO \left( |\lambda|^{2^{-(\xi_\tau-1)}} \right)$, since $\mu_\tau \ge 1$. As $P'$ and $Q'$ are mutually inverse, we obtain
	\[ y_\tau(\lambda) = \lambda^{2^{-\mu_\tau}} \cdot \vartheta_\tau + R_\tau(\lambda) = z_\tau (\lambda) + \lambda \cdot v^h_\tau = \lambda^{2^{-\xi_\tau}} \cdot \eta_\tau + \varkappa_\tau(\lambda). \]
	Thus, $\eta_\tau=\vartheta_\tau\ne 0$ and $\xi_\tau=\mu_\tau$ so that $z_\tau \sim \lambda^{2^{-\mu_\tau}}$. If on the other hand $\mu_\tau = 0$, we obtain $\mu_\sigma=0$ and therefore $y_\sigma \sim \lambda$ for all $\sigma \noge \tau$, including $\tau$. We compute
	\[ \| z_\tau (\lambda) \| = \left\| \ell_\tau \left( (y_\sigma(\lambda) \mid \sigma \noge \tau), \lambda \right) \right\| = \OO (|\lambda|), \]
	which shows $z_\tau \sim \lambda$, completing the proof of (\romannumeral 1).
	
	Conversely, consider the representation of a branch of steady states in \mbox{$\gker{D_v \Gamma_f(0,0)} \times \RR$} given by $\left( (z_\sigma(\lambda))_{\sigma\in\Sigma}, \lambda \right)$ with cell-by-cell asymptotics $z_\sigma \sim \lambda^{2^{-\mu_\sigma}}$. Fix a cell $\tau \in \Sigma$ with $\mu_\tau \ge 1$ first. Then
	\[ z_\tau(\lambda) = \lambda^{2^{-\mu_\tau}} \cdot \vartheta_\tau + R_\tau (\lambda) \]	
	with $\vartheta_\tau \in V\setminus \{0\}$ and $\|R_\tau(\lambda)\| = \OO \left( |\lambda|^{2^{-\mu_\tau-1}} \right)$. As before, we compute the representation in the generalized kernel of the extended system to be
	\[ y_\tau(\lambda) = z_\tau(\lambda) + \lambda \cdot v^h_\tau = \lambda^{2^{-\mu_\tau}} \vartheta_\tau + \lambda \cdot v^h_\tau + R_\tau(\lambda). \] 
	As $\vartheta_\tau \ne 0$, this implies $y_\tau \sim \lambda^{2^{-\mu_\tau}}$. On the other hand, assume $\mu_\tau = 0$, which implies \mbox{$\|z_\tau (\lambda)\| = \OO(|\lambda|)$} as before. In particular, we compute
	\[ \|y_\tau(\lambda)\| = \left\|z_\tau(\lambda) + \lambda\cdot v^h_\tau \right\| = \OO(|\lambda|), \]
	which shows $y_\tau \sim \lambda$ and therefore completes the proof of (\romannumeral 2).
\end{proof}
\begin{cor}
	\label{cor:asymptotics}
	In feedforward fundamental networks, cell-by-cell asymptotics of the full representation of branches of steady states of \eqref{eq:ffext}, agree with those of the representation in $\gker{D_v \Gamma_f(0,0)} \times \RR$, if they respect the ordering relation \eqref{eq:asorder}. That is the cell-by-cell asymptotics are invariant under the bijective conjugation $P = P' \circ \uu{P}^c \colon M^c \to \gker{D_v \Gamma_f(0,0)} \times \RR$.
\end{cor}
\begin{remk}
	\label{rem:asymptotics}
	Similar to \Cref{rem:relaxed3} a slightly more precise statement relaxing the ordering assumption on the cell-by-cell asymptotics \eqref{eq:asorder} is possible. For every cell $\tau$ for which the ordering assumption is satisfied, $\mu_\sigma \le \mu_\tau$ for all $\sigma \noge \tau$, we obtain that the $\tau$-component exhibits the same asymptotics in every representation of the solution branch. Here we do not require the square root orders of \emph{all} cells to be ordered as in \eqref{eq:asorder}.
	\hspace*{\fill}$\triangle$
\end{remk}

These technical results put us in the position to prove the main theorems of the investigation of $1$-parameter steady state bifurcations in fundamental feedforward networks with high-dimensional internal dynamics. For bifurcation problems in which the generalized kernel at the bifurcation point fulfills an additional condition, cell-by-cell asymptotics turn out to be the same as in the case \od. This follows almost directly from \Cref{cor:asymptotics}. In an arbitrary generic bifurcation problem only a slightly weaker statement holds. That is, the amplification effect as in \eqref{eq:defmu} is the same in both cases. The precise square root order of each cell, however, is not necessarily the same in both cases. In order to properly formulate the results we make use of the tensor notation introduced in \Cref{sec:hd}. In particular, we use indices and superscripts $\mathbf{\mathsf{1}}$ and $\mathbf{\mathsf{D}}$ with objects that have been defined before to indicate one-dimensional or $d$-dimensional internal dynamics, respectively, without explicitly defining them. For example $\Gamo_\phi$ is a fundamental network vector field with internal phase space $\Vo$ and internal dynamics $\phi$, while $\Gamd_f$ is an admissible vector field of the same network with internal phase space $\Vd$ and internal dynamics $f$.
\begin{theorem}
	\label{thm:ffddspecific}
	Consider a feedforward fundamental network with one-dimensional internal dynamics and a generic $1$-parameter family of admissible vector fields satisfying the bifurcation assumption \eqref{eq:critical} formulated in the beginning of \Cref{subsec:ffod} with the absolutely indecomposable subrepresentation $Y \subset \No = \bigoplus_{\sigma\in\Sigma}\Vo$ as the generalized kernel at the bifurcation point. Denote the set of all $1$-parameter families of admissible vector fields satisfying the bifurcation assumption \eqref{eq:critical} with $Y$ as the generalized kernel at the bifurcation point by
	\[ \FF = \left\{ \Gamo_\phi\ \left|\ \Gamo_\phi \text{ generic and } \gker{D_v \Gamo_\phi (0,0)} = Y \right. \right\}. \]
	Furthermore, consider a generic $1$\nobreakdash-parameter family $\Gamd_f \colon \Nd \times \RR \to \Nd$ of admissible vector fields on the same network with internal phase space $\Vd \cong \RR^d$ with $d >1$ satisfying the bifurcation assumption \eqref{eq:critical} and the additional condition
	\begin{equation}
		\label{eq:ffddspecific}
		\gker{D_\omega \Gamd_f (0,0)} = Y \ox \langle w \rangle
	\end{equation}
	on the kernel in tensor notation, where $w \in \Vd$ (see \Cref{thm:dechd}).\footnotemark\ Then there is a generic $\Gamo_\phi \in \FF$ such that $\Gamd_f$ exhibits the same pattern of local branches of steady states with the same cell-by-cell asymptotics. 
	
	More precisely, the branches of steady states of $\Gamo_\phi$ are known from \Cref{prop:maxcrit,prop:maxnoncrit}. Denote them by $(x_\sigma(\lambda))_{\sigma\in\Sigma}$, so that \mbox{$x_\sigma \sim \lambda^{2^{-\mu_\sigma}}$} with $\mu_\sigma$ as in \eqref{eq:defmu}. Then each branch of steady states $(w_\sigma(\lambda))_{\sigma\in\Sigma}$ of $\Gamd_f$ uniquely corresponds to a branch $(x_\sigma(\lambda))_{\sigma\in\Sigma}$ of $\Gamo_\phi$ and the square root orders are the same in both cases, i.e., $w_\sigma \sim \lambda^{2^{-\mu_\sigma}}$ for all $\sigma \in \Sigma$.
\end{theorem}
\footnotetext{In the original coordinates this means for any basis $(b^1_\sigma)_{\sigma\in\Sigma},\dotsc,(b^k_\sigma)_{\sigma\in\Sigma}$ of $Y$ the elements $(b^i_\sigma w)_{\sigma\in\Sigma}$ span $\gker{\Gamd_f(0,0)}$.}
\begin{proof}
	Let $\Gamd_f \colon \Nd \times \RR \to \Vd$ be a generic $1$-parameter family of admissible vector fields with $\Gamd_f(0,0)=0$ and
	\[ \gker{D_\omega\Gamd_f(0,0)} = Y \ox \langle w \rangle, \]
	where $Y \subset \No$ is an absolutely indecomposable subrepresentation. Due to the center manifold reduction, the dynamics of $\Gamd_f$ on its center manifold $M^c$ is bijectively conjugate to those of a generic vector field
	\[ F \colon \gker{D_\omega \Gamd_f(0,0)} \times \RR \to \gker{D_\omega \Gamd_f(0,0)}. \]
	In particular, all bifurcating branches of steady states of $\Gamd_f$ can uniquely be represented as branches of steady states of $F$ in $\gker{D_\omega \Gamd_f(0,0)} \times \RR$.	As the center manifold reduction preserves symmetry, the generic steady state bifurcations of the reduced system are entirely classified by $\Sigma$-equivariance.
	
	The same observation holds true for a generic $1$-parameter steady state bifurcation in the case \od. Let $\Gamo_\phi \colon \No \times \RR \to \Vo$ be a generic $1$-parameter family of admissible vector fields satisfying the bifurcation assumption \eqref{eq:critical}. Its dynamics on the center manifold is bijectively conjugate to that of a generic equivariant system 
	\[ G \colon \gker{D_x \Gamo_\phi(0,0)} \times \RR \to \gker{D_x \Gamo_\phi(0,0)} \]
	In particular, if $\Gamo_\phi \in \FF$ the generalized kernel satisfies 
	\[ \gker{D_x \Gamo_\phi(0,0)} = Y \]
	Once again, dynamics on $Y$ is classified by $\Sigma$-symmetry. As a result, there is a generic choice $\Gamo_\phi\in\FF$ such that
	\begin{equation}
		\label{eq:ffddspecificFG}
		F (y \ox w, \lambda) = G (y, \lambda) \ox w
	\end{equation}
	(recall that every element in $Y \ox \langle w \rangle$ can be represented as a pure tensor $y \ox w$ from the proof of \Cref{lem:subrepDD}).
	
	On the other hand, from \Cref{prop:maxcrit,prop:maxnoncrit}, we know all branching solutions for a generic $\Gamo_\phi$. In particular, for any branch $(x_\sigma(\lambda))_{\sigma\in\Sigma} \subset \No$ we know that $x_\sigma \sim \lambda^{2^{-\mu_\sigma}}$, where the $\mu_\sigma$ satisfy the ordering \eqref{eq:asorder}. \Cref{cor:asymptotics} shows that the representation $(y_\sigma(\lambda))_{\sigma\in\Sigma}$ on $Y$ has the same cell-by-cell asymptotics $y_\sigma \sim \lambda^{2^{-\mu_\sigma}}$. For this representation we have $G((y_\sigma(\lambda))_{\sigma\in\Sigma}, \lambda) = 0$ for all $\lambda$ with small absolute value.
	
	Applying \eqref{eq:ffddspecificFG} this directly implies
	\[ F((y_\sigma(\lambda))_{\sigma\in\Sigma} \ox w, \lambda) = G((y_\sigma(\lambda))_{\sigma\in\Sigma}, \lambda) \ox w = 0. \]
	Hence, we obtain a branch of steady states of the generic equivariant vector field $F$ on $Y\ox \langle w \rangle$ by attaching the vector $w\in\Vd\setminus\{0\}$ to each coordinate of the branch $y(\lambda)$:
	\[ (y_\sigma(\lambda))_{\sigma\in\Sigma} \ox w \subset \gker{D_\omega \Gamd_f(0,0)}. \]
	Furthermore, as the dynamics on $Y \times \RR$ and $Y \ox \langle w \rangle \times \RR$ are determined by symmetry, we obtain all branches of steady states of $F$ in this way (the representation on $\Nd$ via $\left\{ A_\sigma \ox \idvd \right\}_{\sigma\in\Sigma}$ is trivial in the $\langle w \rangle$-component). In the original coordinates this branch is denoted by
	\[ (z_\sigma(\lambda))_{\sigma\in\Sigma} = (y_\sigma(\lambda) \cdot w)_{\sigma\in\Sigma} \]
	(see \eqref{eq:tensorelt}). In particular $z_\sigma \sim \lambda^{2^{-\mu_\sigma}}$. Due to the center manifold reduction, this branch uniquely corresponds to a branch $(w_\sigma(\lambda))_{\sigma\in\Sigma}$ of the full system governed by $\Gamd_f$. As $(y_\sigma(\lambda))_{\sigma\in\Sigma}$ satisfies the ordering \eqref{eq:asorder} on its cell-by-cell asymptotics, \Cref{cor:asymptotics} implies $w_\sigma \sim \lambda^{2^{-\mu_\sigma}}$. Note that all steps in this proof require the application of bijective maps. Hence, we indeed have a one-to-one correspondence between branches in the \od-case and branches in the \dd-case.
\end{proof}
\begin{remk}
	Note that condition \eqref{eq:ffddspecific} imposes a restriction of generality. In general, due to \Cref{thm:dechd}, the generalized kernel is isomorphic -- as a subrepresentation -- but not equal to $Y \ox \langle w \rangle$.
	\hspace*{\fill}$\triangle$
\end{remk}
\begin{figure}[h]
	\begin{center}
		\resizebox{.4\linewidth}{!}{	
			\centering
\begin{tikzpicture}[->,
	>=stealth',
	shorten >=1pt,
	auto,
	node distance=1cm,
	main node/.style={line width=1.5pt, circle, scale = 3, draw, font=\sffamily\tiny, inner sep=1pt}]
	\node[main node] (1) {$1$};
	\node[main node, left of=1] (2) {$2$};
	\node[main node, left of=2] (3) {$3$};
	\path[every node/.style={font=\sffamily\small}, line width =1.5pt]
	(2) edge [color = {red}] node {} (1)
	(3) edge [color = {red}, bend left = 5] node {} (2)
	(3) edge [color = {red}, in = 190, out = 170, looseness = 8] node {} (3)
	(3) edge [color = {blue}, bend left = 20] node {} (1)
	(3) edge [color = {blue}, bend right = 5] node {} (2)
	(3) edge [color = {blue}, in = 200, out = 160, looseness = 8] node {} (3)
;
\end{tikzpicture}%
		}%
	\end{center}%
	\caption{A $3$-cell homogeneous feedforward chain.}
	\label{fig:hdff_3cellff}
\end{figure}%
\begin{ex}
	We illustrate the previous remark by recalling \Cref{ex:int} from the introduction.
	The network in \Cref{fig:hdff_3cellff} is the network in \Cref{fig:3cellff} after its input maps have been completed to a monoid. It is easy to see, that it is also a feedforward fundamental network. Similar to before, we observe that the linearization of a $1$-parameter family of admissible vector fields satisfying the bifurcation assumption \eqref{eq:critical} is of the form
	\[ L = \begin{pmatrix} 
		A & \textcolor{red}{B}	& \textcolor{blue}{C} \\
		0 & A 					& \textcolor{red}{B} + \textcolor{blue}{C} \\
		0 & 0 					& A+\textcolor{red}{B} + \textcolor{blue}{C}
	\end{pmatrix}, \]
	where $A,\textcolor{red}{B},\textcolor{blue}{C} \in \gl{V}$. Under the assumption of non-critical maximal cells, we have a (generically simple) eigenvalue $0$ of $A$ -- i.e., $A=0$ in the case \od\ --, while generically $\textcolor{red}{B}$ and $\textcolor{blue}{C}$ are invertible and we compute
	\begin{equation*}
		\operatorname{ker}_0^{\mathbf{\mathsf{1}}} (L^{\mathbf{\mathsf{1}}}) = \left\langle \begin{pmatrix} 1 \\ 0 \\ 0  \end{pmatrix}, \begin{pmatrix} 0 \\ \frac{1}{\textcolor{red}{B}} \\ 0 \end{pmatrix} \right\rangle, \qquad \operatorname{ker}_0^{\mathbf{\mathsf{D}}} (L^{\mathbf{\mathsf{D}}})	= \left\langle \begin{pmatrix} Y \\ 0 \\ 0  \end{pmatrix}, \begin{pmatrix} Y' \\ Y \\ 0 \end{pmatrix} \right\rangle,
	\end{equation*}
	where $Y, Y' \in \Vd$ are such that $AY=0$ and $AY' = \left( \idv - \textcolor{red}{B} \right) Y$. Hence, we see
	\[ \operatorname{ker}_0^{\mathbf{\mathsf{D}}} (L^{\mathbf{\mathsf{D}}}) \ne \left\langle \begin{pmatrix} Y \\ 0 \\ 0  \end{pmatrix}, \begin{pmatrix} 0 \\ \textcolor{red}{B}^{-1}Y \\ 0 \end{pmatrix} \right\rangle \cong \operatorname{ker}_0^{\mathbf{\mathsf{1}}} (L^{\mathbf{\mathsf{1}}}) \ox \langle Y \rangle. \]
	Nevertheless, it can easily be verified that
	\begin{alignat*}{4}
		& A_{\sigma_2}^{\mathbf{\mathsf{D}}}  \begin{pmatrix} Y' \\ Y \\ 0 \end{pmatrix} &&= \begin{pmatrix} Y \\ 0 \\ 0 \end{pmatrix}, \quad && A_{\sigma_3}^{\mathbf{\mathsf{D}}}  \begin{pmatrix} Y' \\ Y \\ 0 \end{pmatrix} &&= 0 \\
		& A_{\sigma_2}^{\mathbf{\mathsf{D}}}  \begin{pmatrix} 0 \\ \textcolor{red}{B}^{-1}Y \\ 0 \end{pmatrix} &&= \begin{pmatrix} \textcolor{red}{B}^{-1}Y \\ 0 \\ 0 \end{pmatrix}, \quad && A_{\sigma_3}^{\mathbf{\mathsf{D}}}  \begin{pmatrix} 0 \\ \textcolor{red}{B}^{-1}Y \\ 0 \end{pmatrix} &&= 0.
	\end{alignat*}
	In particular the actions of $\sigma_2$ and $\sigma_3$ on $(Y',Y,0)^T$ and on $(0,\textcolor{red}{B}^{-1}Y,0)^T$ are conjugate. Hence, we have indeed
	\[ \operatorname{ker}_0^{\mathbf{\mathsf{D}}} (L^{\mathbf{\mathsf{D}}}) \cong \left\langle \begin{pmatrix} Y \\ 0 \\ 0  \end{pmatrix}, \begin{pmatrix} 0 \\ \textcolor{red}{B}^{-1}Y \\ 0 \end{pmatrix} \right\rangle \cong \operatorname{ker}_0^{\mathbf{\mathsf{1}}} (L^{\mathbf{\mathsf{1}}}) \ox \langle Y \rangle. \]
	\hspace*{\fill}$\triangle$
\end{ex}

\Cref{thm:ffddspecific} describes generic steady state bifurcations in a $\dd$-feedforward network only in the special case that the generalized kernel of the vector field is essentially equal to one that also occurs generically in the \od-case. For a generic feedforward fundamental network with high-dimensional internal dynamics the result is slightly weaker. Here we only obtain the same cell asymptotics as in the \od-case for cells $\sigma\in\Sigma$ for which $\mu_\sigma>\mu_\tau$ for all $\tau \nog \sigma$. In particular, this is the case for cells that are critical in the \od-setting. In cells that are non-critical the square root order in the case \dd is less than or equal to that in the case \od. However, as amplification in the case \od is only visible in the critical cells -- $\mu_\sigma$ increases only in critical cells (see \eqref{eq:defmu}) --, this can be interpreted as the same amplification effect in steady state bifurcations in networks with high-dimensional internal dynamics. In order to make the notion of genericity precise, we need to make the following remark first.
\begin{remk}
	\label{rem:hdff_subrepkernel}
	\label{rem:subrepkernel}
	Consider a generic $1$-parameter family of admissible vector fields for a feedforward fundamental network $\Gamd_f \colon \Nd \times \RR \to \Nd$ satisfying the bifurcation assumption \eqref{eq:critical} in the case \dd. As generic steady state bifurcations occur along absolutely indecomposable subrepresentations, we know that $\gker{D_\omega \Gamd_f(0,0)}$ is absolutely indecomposable. Furthermore, due to \Cref{thm:dechd}, there is a decomposition into indecomposable subrepresentations
	\[ \No = Y_1 \oplus \dotsb \oplus Y_s \]
	such that
	\[ \gker{D_\omega \Gamd_f(0,0)} \cong Y_j \ox \langle w \rangle \cong Y_j \]
	for some $1 \le j \le s$. Any subrepresentation $Y_i$ can occur in a generic $1$-parameter steady state bifurcation in the case \od. Furthermore, the subrepresentations are all absolutely indecomposable and in one-to-one correspondence with the eigenvalues of a generic linear admissible map (see \textcite{vonderGracht.2022}). This can readily be proven by exploiting the structure of upper triangular matrices using methods which have been formalized in \cite{Nijholt.2021}.
    That is, there is a generic $1$-parameter family of admissible vector fields for the same fundamental feedforward network $\Gamo_\phi \colon \No \times \RR \to \No$ satisfying the bifurcation assumption \eqref{eq:critical} in the case \od such that
	\[ \gker{D_x \Gamo_\phi(0,0)} \cong Y_j \cong \gker{D_\omega \Gamd_f(0,0)}, \]
	where the isomorphisms are isomorphisms of subrepresentations.
	\hspace*{\fill}$\triangle$
\end{remk}
\begin{theorem}
	\label{thm:ffddgeneric}
	Consider a feedforward fundamental network with internal phase space $\Vd \cong \RR^d$ with $d >1$. Let $\Gamd_f\colon\Nd \times \RR \to \Nd$ be a generic $1$-parameter family of admissible vector fields satisfying the bifurcation assumption \eqref{eq:critical}. Then $\Gamd_f$ exhibits the same amplification effect in its branching steady state solutions as a generic $1$-parameter family of admissible vector fields for the same network with one-dimensional internal dynamics $\Gamo_\phi\colon\No\times\RR \to \No$.
	
	That is, any branch $(w_\sigma(\lambda))_{\sigma\in\Sigma}$ of steady states of $\Gamd_f$ uniquely corresponds to a branch $(x_\sigma(\lambda))_{\sigma\in\Sigma}$ of\ $\Gamo_\phi$. These are known from \Cref{prop:maxcrit,prop:maxnoncrit}. In particular, each cell $\sigma\in\Sigma$ is of square root order $\mu_\sigma$, i.e., \mbox{$v_\sigma \sim \lambda^{2^{-\mu_\sigma}}$} with $\mu_\sigma$ as in \eqref{eq:defmu}. Then
	\[ w_\sigma \sim \lambda^{2^{-\xi_\sigma}} \quad \text{for all} \quad \sigma \in \Sigma, \]
	where
	\begin{equation*}
		\xi_\sigma \quad \begin{cases}
			= \mu_\sigma								&\text{for} \quad \sigma \text{ critical},\\[5pt]
			\in \left\{ 0, \dotsc, \mu_\sigma \right\}	&\text{for} \quad \sigma \text{ non-critical}.
		\end{cases}
	\end{equation*}
	Here, criticality is to be understood with respect to $\Gamo_\phi$  -- the definition of criticality depends on a non-invertible linearization of a $1$-parameter family of admissible vector fields. 
\end{theorem}
\begin{proof}
	The general idea of the proof is similar to that of \Cref{thm:ffddspecific}. Therefore, we omit some of the details and focus on the difficulties that arise when condition \eqref{eq:ffddspecific} is violated. Once again, dynamics restricted to the center manifold of a generic parameter dependent system \mbox{$\Gamd_f \colon \Nd \times \RR \to \Nd$} satisfying the bifurcation assumption \eqref{eq:critical} is bijectively conjugate to that given by a generic equivariant vector field 
	\[ F \colon \gker{D_\omega \Gamd_f(0,0)} \times \RR \to \gker{D_\omega \Gamd_f (0,0)}. \]
	Due to \Cref{thm:dechd}, we know that
	\begin{equation}
		\label{eq:ffddgeneric-ddker}
		\Nd \supset\gker{D_\omega \Gamd_f(0,0)} \cong Y \ox \langle w \rangle,
	\end{equation}
	where $w \in \Vd\setminus\{0\}$ and $Y \subset \No$ is an absolutely indecomposable subrepresentation. \Cref{rem:hdff_subrepkernel} implies, that there is a generic $1$-parameter family of admissible vector fields $\Gamo_\phi \colon \No \times \RR \to \No$ for the same network with one-dimensional internal dynamics such that
	\begin{equation}
		\label{eq:ffddgeneric-odker}
		\gker{D_x \Gamo_\phi (0,0)} \cong Y \subset \No.
	\end{equation}
	
	Generic dynamics -- most importantly steady state bifurcations -- on $Y$ is uniquely determined by symmetry. On the other hand, the same holds for generic dynamics on $\gker{\Dsub_x \Gamo_\phi(0,0)}$. As in the proof of \Cref{thm:ffddspecific}, these are also determined by the reduction of the generic system $\Gamo_\phi$ to its generalized kernel by the center manifold reduction. In particular, using \Cref{prop:maxcrit,prop:maxnoncrit} and \Cref{cor:asymptotics}, we obtain a unique representation $(y_\sigma(\lambda))_{\sigma\in\Sigma} \subset \gker{D_x \Gamo_\phi(0,0)}$ of each branch of $\Gamo_\phi$ such that $y_\sigma \sim \lambda^{2^{-\mu_\sigma}}$, where $\mu_\sigma$ is defined as in \eqref{eq:defmu}. In the same manner, we obtain all branching steady state solutions of a generic $\Ao_\sigma \ox \mathbbm{1}_\Vd$-equivariant system
	\[ G \colon \gker{D_x \Gamo_\phi(0,0)} \ox \langle w \rangle \times \RR \to \gker{D_x \Gamo_\phi(0,0)} \ox \langle w \rangle \]
	as
	\[ (y_\sigma(\lambda))_{\sigma\in\Sigma} \ox w \subset \gker{D_x \Gamo_\phi(0,0)} \ox \langle w \rangle. \]
	This is due to the obvious observation
	\[ \gker{D_x \Gamo_\phi(0,0)} \ox \langle w \rangle \cong \gker{D_x \Gamo_\phi(0,0)} \]
	as representations with respect to $\sigma \mapsto \Ao_\sigma \ox \mathbbm{1}_\Vd$ and $\sigma\mapsto\Ao_\sigma$, respectively. In the remainder of this proof, it is convenient not to use the tensor notation but to rely on the original coordinates $(y_\sigma w)_{\sigma\in\Sigma} \in \gker{D_x \Gamo_\phi(0,0)} \ox \langle w \rangle$. In particular, the generic branch in $\gker{D_x \Gamo_\phi(0,0)} \ox \langle w \rangle$ becomes $(y_\sigma(\lambda)w)_{\sigma\in\Sigma}$. It satisfies
	\[ y_\sigma(\lambda)w \sim \lambda^{2^{-\mu_\sigma}} \]
	for all $\sigma\in\Sigma$ as $w$ is fixed.
	
	Using \eqref{eq:ffddgeneric-ddker} and \eqref{eq:ffddgeneric-odker} we see that there is an isomorphism of subrepresentations
	\[ \Psi \colon \gker{D_\omega \Gamd_f(0,0)} \to \gker{D_x \Gamo_\phi(0,0)} \ox \langle w \rangle. \]
	with inverse $\Xi$. This isomorphism conjugates generic dynamics on the two spaces as it respects symmetry. In particular, all branches of steady states of the generic reduced system $F$ are of the form
	\[ (z_\sigma(\lambda))_{\sigma\in\Sigma} = \Xi \left( (y_\sigma(\lambda))_{\sigma\in\Sigma} w \right) \subset \gker{D_\omega \Gamd_f(0,0)}. \]
	Similar to the proof of \Cref{thm:ffddspecific} we have found a relation of the generic branches of steady states in the case \dd to those in the case \od. However, due to the characterization of $(z_\sigma(\lambda))_{\sigma\in\Sigma}$ using the isomorphism $\Xi$ it is not obvious what the square root orders of each cell are in this representation. Hence, we cannot apply \Cref{cor:asymptotics} directly to obtain square root orders of the corresponding branches of $\Gamd_f$. Hence, we begin by investigating the effect of the isomorphism $\Psi$ on cell-by-cell square root orders similarly to \Cref{lem:3,lem:4}. To that end, we extend $\Psi$ and $\Xi$ trivially to the full space, i.e., $\Psi,\Xi \colon \Nd \to \Nd$. Hence, both maps are equivariant maps with respect to $\sigma \mapsto \Ad_\sigma$ (they are, however, only invertible when restricted to the respective subrepresentations). Due to the equivalence of equivariance and admissibility, they are therefore linear admissible maps for the feedforward fundamental network with internal phase space $\Vd$. In particular, both maps respect the partial order $\nole$ so that we may write
	\begin{subequations}
		\label{eq:ffddgeneric-PsiXi}
		\begin{alignat}{3}
			& \left( \Psi((Z_\sigma)_{\sigma\in\Sigma}) \right)_\tau &	&= \ell_\tau (Z_\sigma \mid \sigma \noge \tau) & &= \ell_\tau^y(Z_\sigma \mid \sigma \noge \tau) \cdot w, \label{eq:ffddgeneric-Psi}\\
			& \left( \Xi((Y_\sigma w)_{\sigma\in\Sigma}) \right)_\tau &	&= \varell_\tau (Y_\sigma w \mid \sigma \noge \tau) & &= \sum_{\sigma \noge \tau} Y_\sigma \varell_\tau^\sigma (w), \label{eq:ffddgeneric-Xi}
		\end{alignat}
	\end{subequations}
	where $\ell_\tau \colon \Vd^{\# \{\sigma\noge\tau\}} \to \Vd$, $\ell_\tau^y\colon \Vd^{\# \{\sigma\noge\tau\}} \to \RR$ as well as $\varell_\tau \colon \colon \Vd^{\# \{\sigma\noge\tau\}} \to \Vd$ and $\varell_\tau^\sigma \colon \Vd \to \Vd$ are suitable linear maps.
	
	Similar to the proofs of \Cref{lem:3,lem:4}, we investigate the implications of this structure of the linear maps $\Psi$ and $\Xi$ on the cell-by-cell asymptotics of $(z_\sigma(\lambda))_{\sigma\in\Sigma} = \Xi \left( (y_\sigma(\lambda))_{\sigma\in\Sigma} w \right)$. By definition
	\[ y_\sigma(\lambda) = \begin{cases}
	\alpha_\sigma \cdot \lambda^{2^{-\mu_\sigma}} + \OO\left( |\lambda|^{2^{-(\mu_\sigma - 1)}} \right)	&\text{if} \quad \mu_\sigma \ge 1, \\
	\OO(|\lambda|)																						&\text{if} \quad \mu_\sigma = 0,
	\end{cases} \]
	where $\alpha_\sigma \in \RR\setminus\{0\}$. Consider a cell $\tau \in \Sigma$ with $\mu_\tau=0$ first. By definition $\mu_\sigma=0$ for all $\sigma \noge \tau$. Hence, using \eqref{eq:ffddgeneric-Xi} we obtain
	\[ \|z_\tau (\lambda)\| = \left\| \sum_{\sigma \noge \tau} y_\sigma(\lambda) \varell_\tau^\sigma (w) \right\| = \OO(|\lambda|) \]
	and, therefore, $z_\tau \sim \lambda$. In particular $\tau$ is of square root order $\iota_\tau=0$ in that representation.
	
	Next, consider a cell $\tau\in\Sigma$ with $\mu_\tau\ge 1$. Then \eqref{eq:ffddgeneric-Xi} gives
	\[ z_\tau (\lambda) = \sum_{\sigma \noge \tau} y_\sigma(\lambda) \varell_\tau^\sigma (w) = \lambda^{2^{-\mu_\tau}} \cdot \eta_\tau + \varkappa_\tau(\lambda), \]
	where $\eta_\tau \in \Vd$ and $\|\varkappa_\tau(\lambda)\|=\OO \left( \lambda^{2^{-(\mu_\tau-1)}} \right)$. In particular, $z_\tau \sim \lambda^{2^{-\iota_\tau}}$ with
	\begin{equation}
		\label{eq:ffddgeneric-iota}
		\iota_\tau \quad \begin{cases}
			= \mu_\tau							&\text{if} \quad \eta_\tau \ne 0, \\
			\in \{ 0, \dotsc, \mu_\tau -1 \}	&\text{if} \quad \eta_\tau = 0.			
		\end{cases}
	\end{equation}
	The second case can be avoided, if additionally $\tau$ is a critical cell with respect to $\Gamo_\phi$. By definition, this implies $\mu_\sigma < \mu_\tau$ for all $\sigma \nog \tau$. Hence, \eqref{eq:ffddgeneric-iota} gives 
	\begin{equation}
		\label{eq:asorderkernel}
		\iota_\sigma \le \mu_\sigma < \mu_\tau \quad \text{for all} \quad \sigma \nog \tau.
	\end{equation}
	As $\Psi$ and $\Xi$ are mutually inverse on the respective subrepresentations, we obtain
	\begin{align*}
		y_\tau(\lambda) w &= \left( \alpha_\tau \cdot \lambda^{2^{-\mu_\tau}} + \OO\left( |\lambda|^{2^{-(\mu_\tau - 1)}} \right) \right) \cdot w \\
		&= \left( \Psi((z_\sigma(\lambda))_{\sigma\in\Sigma}) \right)_\tau \\
		&= \ell_\tau^y(z_\sigma(\lambda) \mid \sigma \noge \tau) \cdot w \\
		&= \left( \beta_\tau (\eta_\tau) \cdot \lambda^{2^{-\mu_\tau}} + \OO\left( |\lambda|^{2^{-(\mu_\tau - 1)}} \right) \right) \cdot w,
	\end{align*}
	where $\beta_\tau \colon \Vd \to \RR$ linear. Therein, the last equation holds due to \eqref{eq:asorderkernel}. Hence, $\beta_\tau(\eta_\tau) = \alpha_\tau \ne 0$ which implies $\eta_\tau \ne 0$. Hence, due to \eqref{eq:ffddgeneric-iota}, $\iota_\tau=\mu_\tau$ so that $z_\tau \sim \lambda^{2^{-\mu_\tau}}$.
	
	Summarizing, we have shown that the representation $(z_\sigma(\lambda))_{\sigma\in\Sigma}$ of the generic branch of steady states $(y_\sigma(\lambda)w)_{\sigma\in\Sigma}$ in $\gker{D_\omega \Gamd_f(0,0)}$ has cell-by-cell asymptotics $z_\sigma \sim \lambda^{2^{-\iota_\sigma}}$. These do not satisfy the ordering \eqref{eq:asorder} on all cells. However, for a cell $\tau$ that is critical with respect to $\Gamo_\phi$, we have that $\iota_\sigma < \iota_\tau$ for all $\sigma \nog \tau$. This follows from \eqref{eq:asorderkernel} and the fact that $\iota_\tau=\mu_\tau$ for such $\tau$. Due to the bijective conjugation between the dynamics of $F$ and those of the generic original system $\Gamd_f$ restricted to its center manifold, each branch $(z_\sigma(\lambda))_{\sigma\in\Sigma}$ of steady states of $F$ uniquely corresponds to a branch of steady states $(w_\sigma(\lambda))_{\sigma\in\Sigma} \subset \Nd$ of $\Gamd_f$. It can be computed as in \eqref{eq:cmf}
	\begin{equation}
		\label{eq:ffddgeneric-w}
		\left( (w_\sigma(\lambda))_{\sigma\in\Sigma} ,\lambda \right) = \left( \mathbbm{1}_{\uu{\XX}^c} + \uu{\psi} \right) \left( Q' \left( (z_\sigma(\lambda))_{\sigma\in\Sigma}, \lambda \right) \right)
	\end{equation}
	using the map whose graph is the center manifold $M^c$. It remains to be investigated, how the square root orders translate to this representation. From \Cref{cor:asymptotics} and \Cref{rem:asymptotics} we immediately obtain
	\begin{equation}
		\label{eq:ffddgeneric-wcrit}
		w_\tau \sim \lambda^{2^{-\iota_\tau}} = \lambda^{2^{-\mu_\tau}}
	\end{equation}
	for all cells $\tau\in\Sigma$ that are critical with respect to $\Gamo_\phi$, as the ordering \eqref{eq:asorder} is satisfied for these cells as well as $\iota_\tau=\mu_\tau$.
	
	We cannot apply \Cref{cor:asymptotics} to cells that are not critical with respect to $\Gamo_\phi$. However, using \Cref{lem:4}, we may take an intermediate step in \eqref{eq:ffddgeneric-w} as
	\[ \left( (X_\sigma(\lambda))_{\sigma\in\Sigma}, \lambda \right) = Q' \left( (z_\sigma (\lambda))_{\sigma\in\Sigma}, \lambda \right) \subset \uu{\XX}^c \subset \Nd \times \RR, \]
	where $ \uu{\XX}^c$ is the generalized kernel of the linearization of the extended system at the bifurcation point $D \uu{\Gamma}_f^\mathbf{\mathsf{D}}(0,0)$ and $Q' \colon \gker{D_{\omega} \Gamd_f(0,0)} \times \RR \to \uu{\XX}^c$ is the isomorphism of $\Sigma$-representations given in \eqref{eq:Qprimeff}. Then $((X_\sigma(\lambda))_{\sigma\in\Sigma}, \lambda)$ is the unique representation of $((z_\sigma(\lambda))_{\sigma\in\Sigma}, \lambda)$ in $\uu{\XX}^c$. Due to \Cref{lem:4}, the $X$-component of this representation exhibits the same cell-by-cell asymptotics as $(z_\sigma(\lambda))_{\sigma\in\Sigma}$, i.e.,
	\[ X_\sigma \sim \lambda^{2^{-\iota_\sigma}} \]
	for all $\sigma \in \Sigma$.
	
	Finally, as in \eqref{eq:ffddgeneric-w} we have
	\[ \left( (w_\sigma(\lambda))_{\sigma\in\Sigma} ,\lambda \right) = \left( \mathbbm{1}_{\uu{\XX}^c} + \uu{\psi} \right) \left( (X_\sigma(\lambda))_{\sigma\in\Sigma}, \lambda \right). \]
	Hence, it remains to investigate what impact the map $\left( \mathbbm{1}_{\uu{\XX}^c} + \uu{\psi} \right)$ has on the square root orders of each cell. \Cref{lem:3} does not apply to this situation, as the square root orders for the branch $(X_\sigma(\lambda))_{\sigma\in\Sigma}$ are given by $\iota_\sigma$ for all $\sigma\in\Sigma$ which do not respect the partial order $\nole$ as in \eqref{eq:asorder}. Nevertheless, a proof similar to that of \Cref{lem:3} allows for a characterization of the square root orders of the branch $(w_\sigma(\lambda))_{\sigma\in\Sigma}$. Recall from \eqref{eq:brcmfgker-quadform} that $\uu{\psi}$ is trivial in the $\lambda$-component, i.e., \mbox{$\uu{\psi} ( X, \lambda ) = ( \uu{\psi}_X (X,\lambda), 0 )$}, and $\uu{\psi}_X$ is a parameter-dependent admissible map of the fundamental network. That is,
	\[ \left( \uu{\psi}_X (X, \lambda) \right)_\tau = q_\tau \left( (X_\sigma \mid \sigma \noge \tau), \lambda \right). \]
	In particular, we compute
	\begin{equation}
		\label{eq:ffddgeneric-quadform}
		w_\tau(\lambda) = X_\tau(\lambda) + q_\tau ((X_\sigma(\lambda) \mid \sigma \noge \tau), \lambda),
	\end{equation}
	for all $\tau \in\Sigma$, where $q_\sigma(0,0)=0$ and $dq_\sigma(0,0)=0$.
	
	Fix a cell $\tau\in\Sigma$ with $\iota_\tau=0$. Then $\iota_\sigma=0$ for all $\sigma \nog \tau$. Hence, we have $\|X_\sigma(\lambda)\| = \OO(|\lambda|) $ for all $\sigma \noge \tau$ (which includes $\tau$). Therefore, from \eqref{eq:ffddgeneric-quadform} we obtain
	\[\|  w_\tau(\lambda) \| = \left\| X_\tau(\lambda) + q_\tau ((X_\sigma(\lambda) \mid \sigma \noge \tau), \lambda) \right\| = \OO(|\lambda|), \]
	which is equivalent to $w_\tau \sim \lambda$. That is $\tau$ has square root order $0 = \iota_\tau=\mu_\tau$.
	
	On the other hand, consider a cell $\tau\in\Sigma$ with $\iota_\tau \ge 1$. We may focus on the case that $\tau$ is non-critical with respect to $\Gamo_\phi$, as we have investigated the other case already in \eqref{eq:ffddgeneric-wcrit}. By definition
	\[ X_\sigma(\lambda) = \lambda^{2^{-\iota_\sigma}} \cdot \vartheta_\sigma + \rho_\sigma(\lambda), \]
	for all $\sigma\in\Sigma$, where $\vartheta_\sigma \in \Vd\setminus\{0\}$ and $\|\rho_\sigma(\lambda)\| = \OO\left( |\lambda|^{2^{-(\iota_\sigma-1)}} \right)$. Hence, using \eqref{eq:ffddgeneric-quadform} we obtain
	\[ w_\tau(\lambda) = X_\tau(\lambda) + q_\tau ((X_\sigma(\lambda) \mid \sigma \noge \tau), \lambda) = \lambda^{2^{-\iota_\sigma}} \cdot \vartheta_\sigma + \lambda^{2^{-\varUpsilon_\tau}} \cdot \varpi_\tau + \varOmega_\tau(\lambda), \]
	where $\varpi_\tau \in \Vd$, $\|\varOmega_\tau(\lambda)\| = \OO\left( |\lambda|^{2^{-(\iota_\tau-1)}} + |\lambda|^{2^{-(\varUpsilon_\tau-1)}} \right)$ and $\varUpsilon_\tau = \max_{\sigma \noge \tau} \iota_\sigma -1$. Recall that $\iota_\sigma \le \mu_\sigma$ for all $\sigma\in\Sigma$. Furthermore,
	\[ \mu_\tau = \max_{\sigma \nog \tau} \mu_\sigma, \]
	as $\tau$ is non-critical for $\Gamo_\phi$. Hence, we obtain
	\[ \varUpsilon_\tau = \max_{\sigma \noge \tau} \iota_\sigma -1 \le \max_{\sigma \noge \tau} \mu_\sigma -1 = \mu_\tau -1 < \mu_\tau. \]
	This proves $w_\tau \sim \lambda^{2^{-\xi_\tau}}$ where $\xi_\tau \le \iota_\tau \le \mu_\tau$.
\end{proof}
\begin{remk}
	\label{rem:critical}	
	In \Cref{thm:ffddspecific,thm:ffddgeneric} we show that the $1$-parameter family of admissible vector fields $\Gamd_f$ in the case \dd admits branching steady states that exhibit the same amplification effect as those of $\Gamo_\phi$ in the case \od. In particular, cells that are critical with respect to $\Gamo_\phi$ exhibit the same asymptotics in both cases -- the same holds for cells $\tau\in\Sigma$ with $\mu_\tau=0$.
	
	On the other hand, consider a cell $\tau$ with $\mu_\tau\ge 1$ that is non-critical with respect to $\Gamd_f$. Assume, furthermore, that for a specific branch of steady states $(w_\sigma(\lambda))_{\sigma\in\Sigma}$, we know $w_\sigma \sim \lambda^{2^{-\iota_\sigma}}$ with $\iota_\sigma\ge 0$ for all $\sigma \nog \tau$. In order to determine the $\tau$-coordinate of this branch, the equation to be solved is
	\[ 0 = f(w_{\sigma_1 \tau}, \dotsc, w_{\sigma_n \tau}, \lambda). \]
	Note that for all $\sigma\in\Sigma$ either $\sigma\tau=\tau$ or $\sigma\tau \nog \tau$. In particular, we may fill in the coordinates of the branch $w_\sigma(\lambda)$ for all $\sigma\nog\tau$. The $\tau$ equation reduces to a bifurcation equation on $\Vd$
	\[ 0 = g(w_\tau, \lambda). \]
	Due to the bifurcation assumption \eqref{eq:critical}, we have
	\[ D_{w_\tau} g(0,0) = D_{w_\tau} f(0,0) = \sum_{\sigma \in \LL_{\tau}} a_\sigma \]
	which is invertible as $\tau$ is assumed not to be critical. Therefore, we may apply the implicit function theorem to obtain
	\[ w_\tau (\lambda) = G((w_\sigma(\lambda) \mid \sigma \nog \tau), \lambda) \]
	such that
	\[ 0 = g(w_\tau(\lambda), \lambda) = f(w_{\sigma_1 \tau}(\lambda), \dotsc, w_{\sigma_n \tau}(\lambda), \lambda) \]
	for all $\lambda$ close to $0$. As $G$ is at least linear up to lowest order in all its arguments, we immediately obtain
	\[ w_\tau \sim \lambda^{2^{-\iota_\tau}} \]
	with $\iota_\tau \le \max_{\sigma \nog \tau} \iota_\sigma$. In particular, no amplification may occur in a cell that is non-critical with respect to $\Gamd_f$. As a result, all cells that provide an amplification must be critical. Since these are precisely the ones that are critical with respect to $\Gamo_\phi$, we obtain that the set of critical cells is the same in both cases.
	\hspace*{\fill}$\triangle$
\end{remk}
\begin{remk}
	The bifurcation results in the case \od in \textcite{vonderGracht.2022} are proven by bare-hands analysis exploiting the partial order in the network. The fact that the equation for cell $\sigma$ depends only on the state variables of cells $\tau \noge \sigma$, allows for inductive computations of solutions. A similar approach is possible in the case \dd as well. Computations in non-critical cells can be solved using the implicit function theorem (compare to \Cref{rem:critical}), critical cells require application of the Lyapunov\nobreakdash-Schmidt reduction. However, taking care of the numerous different cases of branching patterns, depending on a large number of parameters in the equations, and in particular investigating genericity is at least tedious if not factually impossible, due to the arbitrary dimension of the internal phase space.
	\hspace*{\fill}$\triangle$
\end{remk}

\section{Discussion}
\label{sec:discussion}
In the theoretical study of networked dynamical system, one often restricts to the case that the state of each unit evolves in a one-dimensional state space (typically a real vector space). However, many complex systems in applications do not inherently justify this restriction. For instance, coupled oscillators, unless reduced to phase dynamics, require at least a two-dimensional phase space. As a more specific example, the state of each Hodgkin-Huxley neuron in a network consists of four values describing the cell membrane's capacitance, as well as its conductance to three different ions. Similarly, the positional variables of multiple interacting mobile agents in realistic applications require two or three spatial dimensions. In simplifications to one spatial dimension, the collective typically behaves vastly differently.

In this paper, we justify the restriction to one-dimensional phase spaces from a theoretical point of view. Specifically, we study how the representation theory of monoids can be exploited to study bifurcations of coupled cell systems with high-dimensional cell dynamics. It is known that homogeneous coupled cell systems with asymmetric inputs are related to monoids and their representation on the total phase space. This representation uniquely decomposes into indecomposable subrepresentations. We show that the dimension of the internal phase space only affects the multiplicity with which the components occur.
Together with known transversality results for monoid equivariant bifurcation problems, this further leads to powerful statements about generic bifurcations in these networks. In particular, there is a minimal dimension for the internal phase space (depending on the number of parameters and the nature of the bifurcation) that is needed for all possible combinations of indecomposable subrepresentations that a given coupling structure allows to emerge in a generic bifurcation problem. Hence, all possible generic bifurcation problems for any dimension of internal phase space are \emph{equal} to one of those in the minimal dimension in the sense that there is an equivariant isomorphism identifying them. As such, the emerging bifurcation patterns are `qualitatively the same'. This argument is based on the Lyapunov-Schmidt reduction. If more detailed knowledge about the generic bifurcation behavior of a given coupling structure is known, the center manifold reduction can be exploited for a fine-grained relation between bifurcation behavior for different internal phase space dimensions. This is exemplified in this paper at the hand of feedforward networks.

We would like to mention one notable treatment of the question how the dimension of cell-dynamics influences the qualitative bifurcation behavior, which was performed in \cite{Gandhi.2020}. While not being the main focus of their study, the authors prove that $1$- and $2$-parameter bifurcations are effectively independent of the internal dimension (we believe this result to hold true for an arbitrary number of parameters). This apparent contradiction to our results in \Cref{sec:bi} stems from the fact that both works study entirely different classes of coupled cell systems. In fact, both classes are in some sense as different as they can be: While in this paper the network is assumed to be \emph{homogeneous} -- the dynamics of all cells are governed by the same function --, the networks in \cite{Gandhi.2020} are assumed to be \emph{fully inhomogeneous} -- the dynamics of all cells are governed by pairwise (generically strictly) different functions. Both approaches can be motivated by complex systems coming from applications in their own right. Their respective analyses however differ greatly. This can best be illustrated at the hand of a feedforward network in the sense of \Cref{sec:ff}. As we have seen, the homogeneity assumption leads to the phenomenon that multiple cells in the network can be critical at the same time which in turn causes higher multiplicities of critical eigenvalues in bifurcation problems. On the other hand, this effect can generically not occur in fully inhomogeneous networks, as the contribution of all cells to the spectrum is necessarily different. This in turn implies that the multiplicity of critical eigenvalues in a generic bifurcation problem is completely determined by the number of parameters.

Therefore, in our approach, we need to address high-dimensional center subspaces. To that end, we have the advantage that the homogeneity introduces monoid symmetries into the system whose representations and decompositions thereof can be exploited to combinatorically understand which components can make up the center subspace in a given generic bifurcation problem, as outlined above. This algebraic machinery is no longer available in the fully inhomogeneous context. In fact, fully inhomogeneous networks only have trivial monoid symmetries. However, the impact of high-dimensional cell dynamics can be circumvented in a more direct way. In fact, in a generic fully inhomogeneous network vector field, the impact of the different components of high-dimensional cell states on each other are all pairwise different. Hence, the state of a high-dimensional cell can instead be interpreted as the collection of states of multiple one-dimensional cells. This identification does not require any technical steps but is simply a change in perspective (or in notation for that matter). Hence, any bifurcation problem with high-dimensional cells can be identified with a bifurcation problem on a larger network with only one-dimensional cells. At this point, the authors of \cite{Gandhi.2020} go a step beyond what is provided in this paper: They fully classify all $1$- and $2$-parameter bifurcations in fully inhomogeneous networks with one-dimensional cells. By the argument outlined above, this result then immediately yields the same classification for high-dimensional cells.

\section*{Acknowledgement}
\noindent
Parts of this work originated in Sören von der Gracht's doctoral project and are contained in his thesis (``Genericity in Network Dynamics'', 2019 \cite{Schwenker.2019}), written under the primary supervision of Reiner Lauterbach (Universität Hamburg) and co-examined by Bob Rink (Vrije Universiteit Amsterdam) and Ana Paula Dias (Universidade do Porto). The author wishes to express his gratitude to the examiners for helpful comments, discussions, and support.
\medskip

\noindent
SvdG was partially funded by the Deutsche Forschungsgemeinschaft (DFG, German Research Foundation)---453112019.
\medskip

\noindent
EN acknowledges support from the São Paulo Research Foundation (FAPESP, grant no. 2024/00930-2).
\medskip

\noindent
\begin{minipage}{0.8\textwidth}
	This work is part of EN's research program \textit{Designing Network Dynamical Systems through Algebra}, which is financed by the Dutch Research Council (NWO).
\end{minipage}
\begin{minipage}{0.2\textwidth}
	\hfill \includegraphics[width=1.2cm]{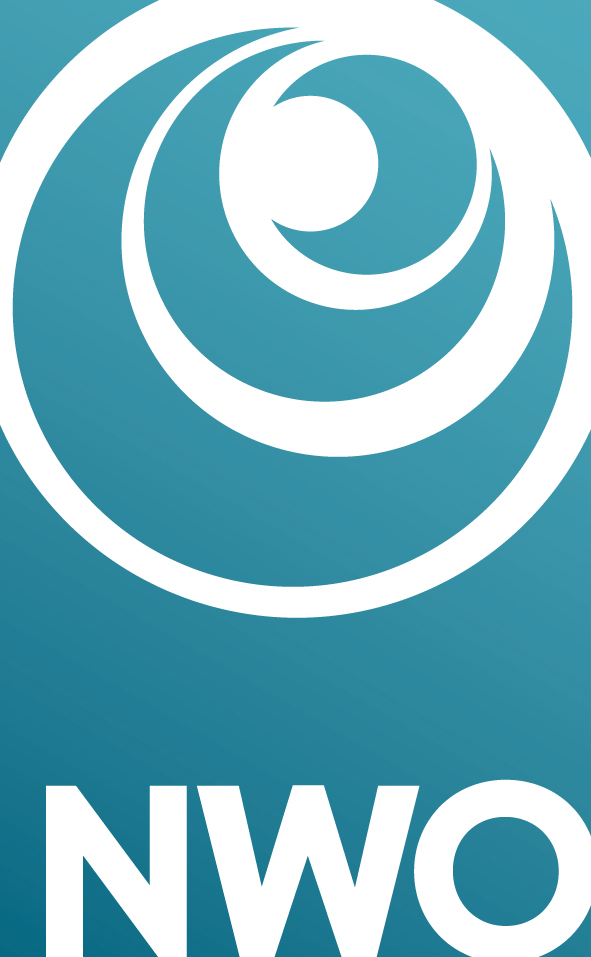}
\end{minipage}

\begingroup
\RaggedRight
\printbibliography
\endgroup
\phantomsection
\addcontentsline{toc}{section}{References}

\end{document}